\documentclass [a4paper,leqno,12pt]{article}
\usepackage[intlimits]{amsmath}
\usepackage{amsbsy}
\usepackage{amstext}
\usepackage{amsthm}
\usepackage{amssymb}
\usepackage[cp1251]{inputenc}
\usepackage[T2A]{fontenc}
\inputencoding{cp1251}
\usepackage{indentfirst}
\usepackage{a4wide}
\usepackage{rotating}
\usepackage{pstricks}
\usepackage[]{graphicx}
\newcommand{\beq}{\begin{equation}}
\newcommand{\eeq}{\end{equation}}
\newcommand{\bs}{\boldsymbol}
\newtheorem{theorem}{Theorem}
\newtheorem{proposition}[theorem]{Proposition}
\newtheorem{lemma}[theorem]{Lemma}
\newtheorem{cor}[theorem]{Corollary}
\theoremstyle{definition}
\newtheorem{remark}[theorem]{Remark}

\numberwithin{equation}{section} \numberwithin{figure}{section}
\numberwithin{theorem}{section} \numberwithin{table}{section}
\def\proof{{\bf \noindent Proof. }}
\def\endprf {\mbox{$\blacksquare$}}
\DeclareMathOperator*{\argmax}{argmax}

\title{Price dynamics in a strategic model of trade between two regions}
\author{Iordan V. Iordanov \footnote{Faculty of Mathematics and Informatics, Sofia University. E-mail: iordanov@fmi.uni-sofia.bg} \and Stoyan V. Stoyanov\footnote{Finanalytica. E-mail: stoyan.stoyanov@finanalytica.com} \and Andrey A. Vassilev\footnote{Faculty of Mathematics and Informatics, Sofia University. E-mail: avassilev@fmi.uni-sofia.bg}}
\date{}

\begin{document} \maketitle 
\begin{abstract}This paper develops a strategic model of trade between two regions in
which, depending on the relation among output, financial resources
and transportation costs, the adjustment of prices towards an
equilibrium is studied. We derive conditions on the relations
among output and financial resources which produce different types
of Nash equilibria. The paths obtained in the process of
converging toward a steady state for prices under discrete-time
and continuous-time dynamics are derived and compared. It turns
out that the results in the two cases differ substantially. Some
of the effects of random disturbances on the price dynamics in
continuous time are also studied.

\end{abstract}

\section{Introduction}

The present work develops a model of trade between two regions in
which, depending on the relation among output, financial resources
and transportation costs, the adjustment of prices towards a
steady state is studied. We assume that there is one type of
traded good and local producers can supply only a fixed amount of
this traded good, which cannot be stored for future consumption.
As usual, prices change to balance supply and demand. In the
chosen setup, the evolution of prices according to an exogenous
rule is studied, starting from pre-specified initial conditions.
More specifically, whenever there are unsold quantities left, the
price is decreased proportionally and when there are local
financial resources unspent, the price is increased
proportionally. This allows us to abstract away from producer
behaviour and focus exclusively on consumers' decisions. The
representative consumers in the two regions seek to maximize their
per-period utility in a strategic situation arising from the need
to compete for scarce resources. We utilize the concept of Nash
equilibrium to characterize optimal behaviour in the game
theoretic interaction. This equilibrium concept has the advantage
of delivering consistent predictions of the outcomes of a game,
assuming that each player takes into account the other players'
optimizing decisions (see Ch.1 in \cite{Fud91} for a more detailed
discussion of the concept).

Under the above setup we derive conditions on the relations among
quantities produced and financial resources, for which different
types of Nash equilibria arise. We also compute the paths obtained
in the process of prices converging toward a steady state. In
certain cases the laws governing price dynamics in discrete time
lead to a zero price in one of the regions, which can be
interpreted as a breakdown of economic activity in the region.
Such pathologies do not arise in the case of continuous-time price
dynamics, where the continuous nature of the adjustment process
provides a natural balancing mechanism against degenerate
stationary points for prices. The stability properties of the
stationary points in the continuous-time case are proved
analytically and illustrated through the behaviour of the phase
trajectories of the system in the presence of stochastic
disturbances.

The paper is organized as follows. Section \ref{sec:model}
introduces the model and key notational conventions. Section
\ref{sec:exist} shows the existence and form of Nash equilibria
for the model under discussion. Section \ref{sec:discprice}
calculates and compares the dynamics governing prices in discrete
time, while section \ref{sec:contprice} presents the counterpart
analysis in the continuous-time case. The proofs of the results
from section \ref{sec:discprice} are provided in the appendix.
Section \ref{sec:stochprice} contains the results of some
simulations for the continuous-time case with stochastic shocks.
Partial announcements of the results reported in this paper
appeared in \cite{Ior04} and \cite{Vas05}.

\section{The model}\label{sec:model}

We consider the consumption decisions of two economic agents
occupying distinct spatial locations, called region $I$ and $II$,
respectively. The consumer in region $I$ (or, shortly, consumer
$I$) exogenously receives money income $Y_1>0$ in each period.
Similarly, the consumer in region $II$ (consumer $II$) receives
money income $Y_2>0$. For each period $t$, in region $i$, $i=1,2$,
a fixed quantity $q_i>0$ of a certain good is supplied at a price
$p_{i,t}$. The consumers place orders for the desired quantities
in each region, observing their budget constraints and incurring
symmetric transportation costs $\rho>0$ per unit of shipment from
the ``foreign'' region. Each consumer attempts to maximize their
total consumption for the current period. Consumers can be
considered myopic in that they do not optimize their consumption
over a specified time horizon but their decisions are confined
only to the current period.

In cases when total orders for the respective region exceed the
quantity available, the following distribution rule is applied:
first, the order of the local consumer is executed to the extent
possible and then the remaining quantity, if any, is allocated to
the consumer from the other region. We sometimes refer to this
distribution scheme and its consequences as \emph{local
dominance}. It is clear then that the choice of orders to be
placed has a strategic element to it, since the actual quantity
received by the consumer depends on the choices made by the
counterpart in the other region. The agents are assumed to have
complete knowledge of all the relevant aspects of the situation
under discussion.

More formally, for each period $t$ we model the above situation as
a static noncooperative game of complete information. Denote by
$\alpha$ and $\beta$ the orders placed by consumer $I$ in region
$I$ and $II$, respectively. In an analogous manner, $\gamma$ and
$\delta$ stand for the orders of consumer $II$ in regions $I$ and
$II$, all orders obviously being nonnegative quantities. In period
$t$ consumer $I$'s strategy space $S_1$ is determined by the
budget constraint and the nonnegativity restrictions on the
orders: \beq \label{eq:ss1}S_1=\{(\alpha,\beta) \in
\mathbb{R}^2|\alpha p_{1,t}+ \beta (p_{2,t}+\rho)\leq Y_1,\quad
\alpha,~\beta \geq 0 \}.\eeq Consumer $II$'s strategy space in
period $t$ is \beq \label{eq:ss2} S_2=\{(\gamma,\delta) \in
\mathbb{R}^2|\gamma (p_{1,t}+\rho)+ \delta p_{2,t}\leq Y_2,\quad
\gamma,~\delta \geq 0 \}.\eeq Below we adopt the shorthand
$p_{1,t}':= p_{1,t}+\rho$ and $p_{2,t}':=p_{2,t}+\rho$. We also
omit the subscript $t$ whenever it is evident from the context or
irrelevant.

The payoff function for consumer $I$ is given by \beq
\begin{split} P_1(\alpha,\beta,\gamma,\delta)= &
\min (\alpha,~q_1)+\min(\beta,~q_2-\min (\delta,~q_2))\equiv \\
\equiv & \min (\alpha,~q_1)+\min (\beta,~\max
(0,~q_2-\delta))\end{split}\label{eq:payoff1} \eeq and that for
consumer $II$ by \beq
\begin{split}P_2(\alpha,\beta,\gamma,\delta)= & \min
(\gamma,~q_1-\min (\alpha,~q_1) )+\min(\delta,~q_2)\equiv \\
\equiv & \min(\delta,~q_2)+\min(\gamma,~\max(0,~q_1-\alpha)).
\end{split} \label{eq:payoff2}\eeq

Any unspent fraction of the current-period income is assumed to
perish and consequently the accumulation of stocks of savings is
not allowed in the model. Similarly, the goods available each
period cannot be stored for future consumption. Let $q_i^{cons}$
denote the total amount consumed in region $i$ and $Y_i^{res}$
stand for the part of the region $i$'s income not spent in the
other region. In other words, $q_1^{cons}:=\alpha_0+\gamma_0$,
$q_2^{cons}:=\beta_0+\delta_0$,
$Y_{1,t}^{res}:=Y_1-p'_{2,t}\beta_0$ and
$Y_{2,t}^{res}:=Y_2-p'_{1,t}\gamma_0$.

There are two mutually exclusive situations leading to an
adjustment in prices. First, if the quantity available in the
respective region has not been entirely consumed, prices are
adjusted downwards. In discrete time this is captured by the
equation \beq
\frac{p_{i,t}-p_{i,t+1}}{p_{i,t}}=\frac{q_i-q_{i,t}^{cons}}{q_i}\textrm{
or } p_{i,t+1}q_i=p_{i,t}q_{i,t}^{cons}. \label{eq:pricedown} \eeq
Clearly, if $q_i^{cons}=0$, then $p_{i,t+1}=0$. Second, if
$Y_i^{res}$ is not entirely exhausted in absorbing local supply,
which can be expressed in value terms as $p_iq_i$, then the price
$p_{i,t}$ is adjusted upwards to $p_{i,t+1}$ to ensure residual
income exhaustion: \beq
\frac{p_{i,t+1}-p_{i,t}}{p_{i,t}}=\frac{Y_{i,t}^{res}-p_{i,t}q_i}{p_{i,t}q_i}
\textrm{ or } p_{i,t+1}q_i=Y_{i,t}^{res}. \label{eq:priceup} \eeq
Obviously the fraction of $Y_i$ spent on the ``foreign'' market
cannot be attracted back for domestic consumption if local prices
are increasing. Later we formally prove the claim that the two
situations leading to price adjustment cannot occur
simultaneously. As usual, we consider prices in a steady
state\footnote{Equivalently, we say that the price adjustment
process has reached a \emph{stationary point}.} if the rules given
by equations \eqref{eq:pricedown} and \eqref{eq:priceup} do not
lead to a change in prices.

For the above model we are interested in two main questions.
First, it would be desirable to establish the existence of a Nash
equilibrium for the one-period game and specify it in closed form.
Second, one would like to be able to trace out the price dynamics
entailed by a sequence of one-period games for a given set of
initial conditions $p_{1,0},~p_{2,0},~q_1,~q_2,~Y_1,~Y_2$ and
$\rho$, and characterize their properties.

\section{Existence and form of equilibrium}\label{sec:exist}

In this section we study the existence and properties of the most
popular equilibrium concept -- that of Nash equilibrium -- for the
model specified above,for a fixed time period $t$. Our basic tool
for establishing existence is a theorem \cite[p. 72]{Fri90}
asserting that at least one Nash equilibrium exists for a game of
complete information for which: \begin{itemize}
  \item[(a)] the strategy spaces of all
players are compact and convex subsets of $\mathbb{R}^m$;
  \item[(b)] all payoff functions are defined, continuous and bounded over the
strategy space of the game, and
  \item[(c)] any payoff function is
quasiconcave in the player's own feasible strategies for a fixed
strategy profile of the opponents.
\end{itemize} We remind the reader that a function $f : X \rightarrow
\mathbb{R}$ is called quasiconcave if for any $x,y \in X$ we have
$f(\lambda x+ (1-\lambda) y) \geq \min \{ f(x), f(y) \}$ for all
$\lambda \in (0,1)$.

Properties (a) and (b) are immediately verified for our model when
the prices $p_i$ are positive. (If a price is zero, economically
plausible restrictions are imposed on the model in order to ensure
that the above properties hold in this case as well; see SR3.) To
establish property (c) note that the payoff function for each
consumer is separable in the consumer's orders and each component
of the sum in the payoff is a concave function in the respective
order. These observations entail the concavity and hence the
quasiconcavity of the payoffs.

Since all the hypotheses of the existence theorem are satisfied
for our model, it has at least one Nash equilibrium. We proceed to
compute the Nash equilibria for all possible configurations of
$Y_1,~Y_2,~q_1,~q_2,~p_1,~p_2$ and $\rho$. To this end, we derive
the best-reply correspondences (see \cite[pp. 69-75]{Fri90} ) for
the two consumers. We remind that these correspondences are
defined as follows. Let us denote
$$u_1(\alpha,\beta)=P_1(\alpha,\beta,\gamma,\delta)$$ for fixed
values of $(\gamma,\delta)\in S_2$ and
$$u_2(\gamma,\delta)=P_2(\alpha,\beta,\gamma,\delta)$$ for fixed
values of $(\alpha,\beta)\in S_1$. Let
$(\bar{\alpha}(\gamma,\delta),\bar{\beta}(\gamma,\delta))$ be
defined as $$\argmax_{(\alpha,\beta)\in S_1}u_1(\alpha,\beta)$$
and, similarly, let
$(\bar{\gamma}(\alpha,\beta),\bar{\delta}(\alpha,\beta))$ be
defined as $$\argmax_{(\gamma,\delta)\in S_2}u_2(\gamma,\delta).$$
The (multivalued) correspondence
$(\alpha,\beta,\gamma,\delta)\mapsto
(\bar{\alpha},\bar{\beta},\bar{\gamma},\bar{\delta})$ is called
best-reply correspondence for the problem, with
$r_1(\gamma,\delta):=(\bar{\alpha},\bar{\beta})$ and
$r_2(\alpha,\beta):=(\bar{\gamma},\bar{\delta})$ being the
best-reply functions for consumers $I$ and $II$, respectively.

The derivation of the best-reply correspondences is
straightforward and we omit the details, presenting only the
end-results. Table \ref{tab:breply1} presents the best-reply
correspondence for consumer $I$ and Table \ref{tab:breply2} shows
the best-reply correspondence for consumer II. For simplicity in
the tables we use $\alpha$ instead of $\bar{\alpha}$ etc. for the
equilibrium values. We note in advance that in the course of the
price adjustment process one of the prices can become zero, in
which case the best reply correspondences take a different form
(see Tables \ref{tab:breply1prime} and \ref{tab:breply2prime}
below).

\begin{sidewaystable}
  \centering 
  \begin{tabular}{|c|c|c|c|} \hline
     & \multicolumn{2}{c|}{$\mathbf{A:\frac{Y_1}{p_1}> q_1}$} &  $\mathbf{B:\frac{Y_1}{p_1}\leq q_1}$
     \\\hline
    $\mathbf{I.~q_2-\bs{\delta} \leq 0}$ & \multicolumn{2}{c|}{$q_1\leq\alpha\leq \frac{Y_1}{p_1},~0\leq\beta\leq \frac{Y_1-p_1\alpha}{p_2'}$, SR1: $\alpha=q_1$, $\beta=0$}   & $\alpha=\frac{Y_1}{p_1},~\beta=0$    \\\hline
    & $\mathbf{A_1:0<q_1<\frac{Y_1-p_2'(q_2-\bs{\delta})}{p_1}}$ & $\mathbf{A_2:\frac{Y_1-p_2'(q_2-\bs{\delta})}{p_1}\leq q_1<\frac{Y_1}{p_1}}$ &
      \\\cline{2-4}
     &   &  $\mathsf{(1)}:\alpha=q_1,~\beta=\frac{Y_1-p_1q_1}{p_2'}$ & $\mathsf{(1)}:\alpha=\frac{Y_1}{p_1},~\beta=0$ \\
     $\mathbf{II.~q_2-\bs{\delta} \in \left(0,\frac{Y_1}{p_2'}\right]}$& $q_1\leq\alpha\leq \frac{Y_1-p_2'(q_2-\delta)}{p_1}$ & $\mathsf{(2)}:\frac{Y_1-p_2'(q_2-\delta)}{p_1}\leq\alpha\leq q_1,~$  & $\mathsf{(2)}: \alpha=\frac{Y_1-p_2'\beta}{p_1},$ \\
      & $q_2-\delta \leq \beta \leq \frac{Y_1-p_1\alpha}{p_2'}$ & $\beta=\frac{Y_1-p_1\alpha}{p_2'}$, SR2: as in (1) & $0\leq\beta\leq\frac{p_2'(q_2-\delta)}{p_1}$, SR2: as in (1) \\
     & SR1: $\alpha=q_1$, $\beta=q_2-\delta$ & $\mathsf{(3)}:\alpha=\frac{Y_1-p_2'(q_2-\delta)}{p_1},$ & $\mathsf{(3)}:\alpha=\frac{Y_1-p_2'(q_2-\delta)}{p_1},$ \\
     &  & $\beta=q_2-\delta$ & $\beta=q_2-\delta$ \\\hline
     & \multicolumn{2}{c|}{$\mathsf{(1)}:\alpha=q_1,~\beta=\frac{Y_1-p_1 q_1}{p_2'}$} &  $\mathsf{(1)}:\alpha=\frac{Y_1}{p_1},~\beta=0$  \\
     $\mathbf{III.~q_2-\bs{\delta} > \frac{Y_1}{p_2'}}$ & \multicolumn{2}{c|}{$\mathsf{(2)}:0\leq\alpha\leq q_1,~\beta=\frac{Y_1-p_1 \alpha}{p_2'}$, SR2: as in (1)} &  $\mathsf{(2)}:0\leq\alpha\leq \frac{Y_1}{p_1},$  \\
     &  \multicolumn{2}{c|}{$\mathsf{(3)}:\alpha=0,~\beta=\frac{Y_1}{p_2'}$}&   $\beta=\frac{Y_1-p_1 \alpha}{p_2'}$, SR2: as in (1) \\
     & \multicolumn{2}{c|}{$\qquad \qquad \qquad$} &  $\mathsf{(3)}:\alpha=0,~\beta=\frac{Y_1}{p_2'}$ \\ \hline
     \multicolumn{4}{c}{Shorthand notation used: $\mathsf{(1)}$ for $\mathsf{p_1<p_2'}$, $\mathsf{(2)}$ for
     $\mathsf{p_1=p_2'}$ and $\mathsf{(3)}$ for
     $\mathsf{p_1>p_2'}$}\\
     \multicolumn{4}{c}{Whenever the shorthand notation is not employed, the result should be taken to apply to each of the three cases.}
  \end{tabular}\caption{Best-reply correspondence for consumer $I$. }\label{tab:breply1}
\end{sidewaystable}

\begin{sidewaystable}
 \centering 
  \begin{tabular}{|c|c|c|c|} \hline
    & \multicolumn{2}{c|}{$\mathbf{A:\frac{Y_2}{p_2}> q_2}$} & $\mathbf{B:\frac{Y_2}{p_2}\leq q_2}$ \\\hline
    $\mathbf{I.~q_1-\bs{\alpha} \leq 0}$ & \multicolumn{2}{c|}{$0\leq\gamma\leq \frac{Y_2-p_2q_2}{p_1'},~q_2\leq\delta\leq \frac{Y_2-p_1'\gamma}{p_2}$, SR1: $\gamma=0$, $\delta=q_2$}   & $\gamma=0,~\delta=\frac{Y_2}{p_2}$    \\\hline
    & $\mathbf{A_1:0<q_2<\frac{Y_2-p_1'(q_1-\bs{\alpha})}{p_2}}$ & $\mathbf{A_2:\frac{Y_2-p_1'(q_1-\bs{\alpha})}{p_2}\leq q_2<\frac{Y_2}{p_2}}$ &
      \\\cline{2-4}
     &   &  $\mathsf{(1)}:\gamma=\frac{Y_2-p_2q_2}{p_1'},~\delta=q_2$ & $\mathsf{(1)}:\gamma=0,~\delta=\frac{Y_2}{p_2}$ \\
     $\mathbf{II.~q_1-\bs{\alpha} \in \left(0,\frac{Y_2}{p_1'}\right]}$& $q_1-\alpha\leq\gamma\leq \frac{Y_2-p_2q_2}{p_1'}$ & $\mathsf{(2)}:\frac{Y_2-p_2q_2}{p_1'}\leq\gamma\leq q_1-\alpha,~$  & $\mathsf{(2)}:0\leq\gamma\leq q_1-\alpha,$ \\
      & $q_2\leq\delta \leq \frac{Y_2-p_1'\gamma}{p_2}$ & $\delta=\frac{Y_2-p_1'\gamma}{p_2}$, SR2: as in (1) & $\delta=\frac{Y_2-p_1'\gamma}{p_2}$, SR2: as in (1) \\
     & SR1: $\gamma=q_1-\alpha$, $\delta=q_2$ & $\mathsf{(3)}:\gamma=q_1-\alpha,$ & $\mathsf{(3)}:\gamma=q_1-\alpha,$ \\
     &  & $\delta=\frac{Y_2-p_1'\gamma}{p_2}$ & $\delta=\frac{Y_2-p_1'\gamma}{p_2}$ \\\hline
     & \multicolumn{2}{c|}{$\mathsf{(1)}:\gamma=\frac{Y_2-p_2 q_2}{p_1'},~\delta=q_2$} &  $\mathsf{(1)}:\gamma=0,~\delta=\frac{Y_2}{p_2}$  \\
     $\mathbf{III.~q_1-\bs{\alpha} > \frac{Y_2}{p_1'}}$ & \multicolumn{2}{c|}{$\mathsf{(2)}:\frac{Y_2-p_2q_2}{p_1'}\leq\gamma\leq \frac{Y_2}{p_1'},~\delta=\frac{Y_2-p_1' \gamma}{p_2}$, SR2: as in (1)} &  $\mathsf{(2)}:0\leq\gamma\leq \frac{Y_2}{p_1'},$  \\
     & \multicolumn{2}{c|}{$\mathsf{(3)}:\gamma=\frac{Y_2}{p_1'},~\delta=0$} &   $\delta=\frac{Y_2-p_1' \gamma}{p_2}$, SR2: as in (1) \\
     & \multicolumn{2}{c|}{$\qquad \qquad \qquad$} &  $\mathsf{(3)}:\gamma=\frac{Y_2}{p_1'},~\delta=0$ \\ \hline
     \multicolumn{4}{c}{Shorthand notation used: $\mathsf{(1)}$ for $\mathsf{p_1'>p_2}$, $\mathsf{(2)}$ for
     $\mathsf{p_1'=p_2}$ and $\mathsf{(3)}$ for
     $\mathsf{p_1'<p_2}$}\\
     \multicolumn{4}{c}{Whenever the shorthand notation is not employed, the result should be taken to apply to each of the three cases.}
  \end{tabular}\caption{Best-reply correspondence for consumer $II$. }\label{tab:breply2}
\end{sidewaystable}

Several comments are in order with respect to the Tables
\ref{tab:breply1} and \ref{tab:breply2}. Because of the presence
of the parameters $(q_i,Y_i,p_i,\rho)$, as well as the different
feasible values of the fixed variables, the procedure for
maximizing $u_i,~i=1,2,$ can be decomposed into different cases in
a natural manner. For example, to find $\max u_1$ we get rid of
$\max (0,q_2-\delta)$ by successively analyzing the two cases $q_2
\leq \delta$ and $q_2>\delta$. The second case in turn decomposes
into two subcases depending on whether the quantity $q_2-\delta$
is smaller or greater than the maximal feasible value $Y_1/p'_2$
of $\beta$. For $\alpha$, which we compare with $q_1$, the
geometry of the feasible set depends on, first, the size of the
difference $q_1-\alpha_0$, where $\alpha_0$ is determined from the
condition $Y_1=p_1\alpha_0+p'_2(q_2-\delta)$, assuming $q_2-\delta
\in (0,Y_1/p'_2)$ and, second, by the ratio between the size of
local supply $q_1$ and the maximum purchasing power of local
income, $Y_1/p_1$.

With the aid of the best-reply correspondences we can compute the
Nash equilibria for the game as solutions to a system of
equations. However, uniqueness is not guaranteed in this model and
we therefore have to resort to additional rules for equilibrium
selection in order to choose a single equilibrium. To this end we
define the following supplementary selection rules (SR), which we
deem logical from a practical viewpoint:
\begin{itemize} \item[SR1] \textbf{(Expenditure minimization)} For a set of Nash equilibria yielding the same utility we select the one minimizing the expenditures
made. (The expenditures made by the first consumer are
$p_1\alpha+p'_2\beta$ and those made by the second consumer are
$p'_1\gamma+p_2\delta$.)
\item[SR2] \textbf{(Home bias)} If more than one Nash equilibrium with the same utility can be obtained with the same (minimal) expenditure, then we select the one in which consumers receive the maximum amount possible in their own region in preference over the ``foreign''
consumer. (I.e. if for the first market we have
$p_1\alpha+p'_2\beta=const$ for more than one point
$(\alpha,\beta)$, we choose the point with the largest value of
$\alpha$. We proceed analogously for the second market.)
\item[SR3] \textbf{(Free disposal)} In the degenerate case when a price is equal to zero, we assume that the actual amount bought is equal to the quantity available in the respective region.  \end{itemize}

SR3 points to a modification in the best-reply correspondences
required in the degenerate case of a zero price. If $p_1=0$, i.e.
$p'_1=\rho$, Table \ref{tab:breply2} should be modified into a
table identical to the former with $p'_1=\rho$ and $p_2>0$. Table
\ref{tab:breply1} should be replaced by Table
\ref{tab:breply1prime}.

\begin{table}[ht]\centering
\begin{tabular}{|l|l|l|}
\hline
\textbf{Subcase} & $\bs{(\bar{\alpha},\bar{\beta})}$ & \textbf{Selection as per SR3} \\
\hline
$\mathbf{q_2-\bs{\delta} \leq 0}$ & $q_1 \leq \alpha$ & $(q_1,0)$ \\
 & $0 \leq \beta \leq \frac{Y_1}{p'_2}$ &  \\
\hline
$\mathbf{0<q_2-\bs{\delta} \leq \frac{Y_1}{p'_2}}$ & $q_1 \leq \alpha$ & $(q_1,q_2-\delta)$ \\
 & $q_2-\delta \leq \beta \leq \frac{Y_1}{p'_2}$ &  \\
\hline
$\mathbf{\frac{Y_1}{p'_2}<q_2-\bs{\delta}}$ & $q_1 \leq \alpha$ & $\left ( q_1,\frac{Y_1}{p'_2}  \right )$ \\
 & $\beta=\frac{Y_1}{p'_2}$ &  \\
\hline
\end{tabular} \caption{Modification of Table \ref{tab:breply1} for the degenerate
case $p_1=0$.}\label{tab:breply1prime}
\end{table}

\begin{table}[ht]\centering
\begin{tabular}{|l|l|l|}
\hline
\textbf{Subcase} & $\bs{(\bar{\gamma},\bar{\delta})}$ & \textbf{Selection as per SR3} \\
\hline
$\mathbf{q_1-\bs{\alpha} \leq 0}$ & $0 \leq \gamma \leq \frac{Y_2}{p'_1}$ & $(0,q_2)$ \\
 & $q_2 \leq \delta$ &  \\
\hline
$\mathbf{0<q_1-\bs{\alpha} \leq \frac{Y_2}{p'_1}}$ & $q_1-\alpha \leq \gamma \leq \frac{Y_2}{p'_1}$ & $(q_1-\alpha,q_2)$ \\
 & $q_2 \leq \delta$ &  \\
\hline
$\mathbf{\frac{Y_2}{p'_1}<q_1-\bs{\alpha}}$ & $\gamma=\frac{Y_2}{p'_1}$ & $\left (\frac{Y_2}{p'_1},q_2  \right )$ \\
 & $q_2 \leq \delta$ &  \\
\hline
\end{tabular}
\caption{Modification of Table \ref{tab:breply2} for the
degenerate case $p_2=0$.}\label{tab:breply2prime}
\end{table}
If $p_2=0$, i.e. $p'_2=\rho$, Table \ref{tab:breply1} should be
replaced by an identical table with $p'_2=\rho,~p_1=0$ and Table
\ref{tab:breply2} should be replaced by Table
\ref{tab:breply2prime}.

In an analogous manner, by using SR1 and SR2 we can dispose of the
multiplicity of solutions and arrive at a unique Nash equilibrium
$(\alpha,\beta,\gamma,\delta)$, as reflected in Tables
\ref{tab:breply1} and \ref{tab:breply2}. The latter can be
obtained as a solution of the system \beq \left \{
\begin{array}{ccc}
  \bar{\alpha}(\gamma,\delta)& = & \alpha \\
  \bar{\beta}(\gamma,\delta) & = & \beta \\
  \bar{\gamma}(\alpha,\beta) & = & \gamma \\
  \bar{\delta}(\alpha,\beta) & = & \delta ,
\end{array}\right . \label{eq:syssol}\eeq which we derive explicitly in section
\ref{sec:discprice} and in the appendix.

\section{Price dynamics in the discrete-time case}\label{sec:discprice}

By definition, a \emph{steady state} (\emph{point of equilibrium},
p.e.) for the prices $p_i$ is the value for which consumption (as
given by the \emph{Nash equilibrium}
$(\alpha,\beta,\gamma,\delta)$, NE for brevity) leads to a
complete depletion of both the available quantities $q_i$ of the
good  and the financial resources $Y_i$. In other words, for $q_i$
and $Y_i$ exogenously given, we have $\alpha+\gamma=q_1$,
$\beta+\delta=q_2$, $p_1\alpha+p'_2\beta = Y_1$, $p'_1 \gamma +
p_2 \delta = Y_2 $. As a result prices are not adjusted for the
next period but retain their current values.

Naturally, if the initial values of prices are not a p.e., they
are corrected prior to next period's consumption, as described
above. The present section aims to characterize the evolution of
prices for all possible values of $Y_i$ (for fixed $q_i,\rho>0$).
For this purpose, it is convenient to present the results in
$Y_1$-$Y_2$ space. In accordance with the different cases
presented in Tables \ref{tab:breply1} and \ref{tab:breply2}, we
partition (by means of a set of lines and additional restrictions)
the nonnegative quadrant of the plane into disjoint subsets of
points\footnote{Formally, we should consider the set of all
possible combinations of incomes for the two regions, which
coincides with the nonnegative quadrant. A generic point in this
set is denoted $(\tilde{Y}_1,\tilde{Y}_2)$, while the particular
income pair under consideration is $(Y_1,Y_2)$. Therefore, the
definitions of all the zones and lines below should be presented
in terms of $\tilde{Y}_1$ and $\tilde{Y}_2$. However, to simplify
the notation we sometimes depart from this convention when no
confusion can arise (e.g. when defining the various zones) and
write the objects simply in terms of $Y_1$ and $Y_2$.}
$(\tilde{Y}_1,\tilde{Y}_2)$ for which a unique NE exists. Each NE
corresponding to an element of this partition is presented as a
closed-form expression involving the exogenous parameters.
However, the partition itself (or, respectively, the set of
lines), crucially depends on the values of $p_i$. After prices
have been adjusted, the point $(Y_1,Y_2)$ may turn out to be in
another element of the partition and, possibly, require another
round of adjustment etc. This evolution of the prices will be
detailed below, where we list the p.e.s attained (after a finite
or infinite number of steps) for each initial point.

\begin{figure}[ht]
\centering
\input{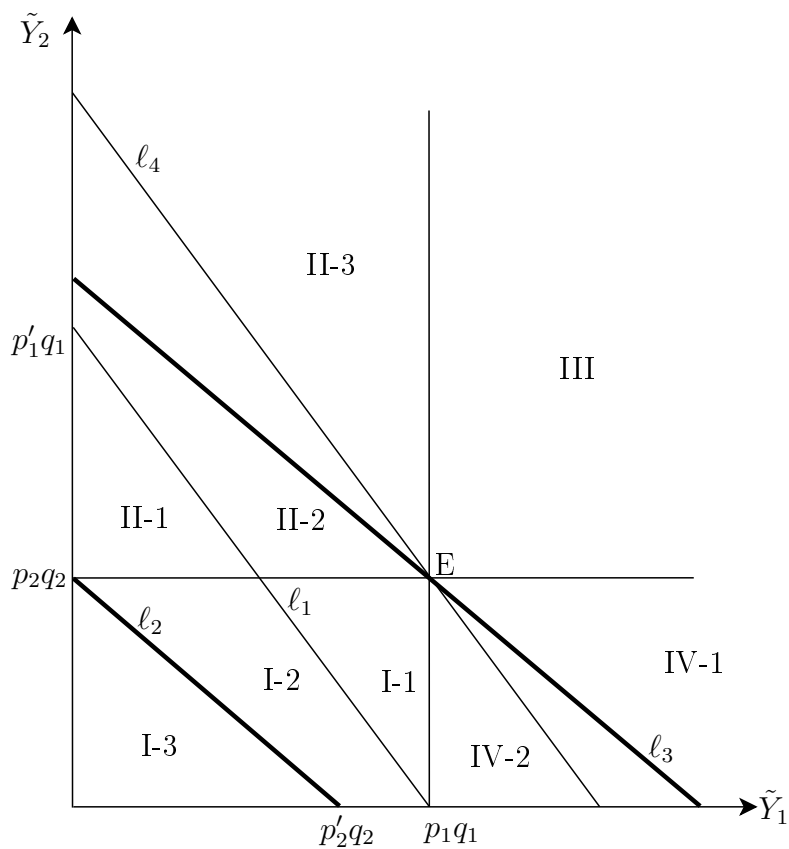}
\caption{Income space partition obtained for a fixed parameter
set.} \label{fig:income}
\end{figure}

The two main cases $A$ and $B$ in the tables define the lines
$\tilde{Y}_i=p_iq_i$. The lines divide the nonnegative quadrant
into four zones (see Figure \ref{fig:income}). We label these
zones in roman numerals: \begin{enumerate} \item[I)]
$Y_1<p_1q_1,~Y_2<p_2q_2$
\item[II)] $Y_1<p_1q_1,~Y_2 \geq p_2q_2$ \item[III)]$Y_1\geq p_1q_1,~Y_2 \geq
p_2q_2$ \item[IV)]$Y_1\geq p_1q_1,~Y_2 < p_2q_2$
\end{enumerate}

The situations obtained in II) and IV) are symmetric with the
roles of the consumers simply being swapped. In order to reduce
the number of cases explored, however, we shall break this
symmetry and make the assumption \beq p_2q_2<p_1q_1
\label{eq:assump}\eeq

The restrictions specified in the left-hand columns of Tables
\ref{tab:breply1} and \ref{tab:breply2} for maximal values of the
priority orders $\alpha$ and $\delta$ define the lines \beq
\ell_1:q_1=\frac{\tilde{Y}_1}{p_1}+\frac{\tilde{Y}_2}{p'_1},\textrm{
i.e. }p_1q_1=\tilde{Y}_1+(p_1/p'_1)\tilde{Y}_2,
\label{eq:ell1line}\eeq \beq
\ell_2:q_2=\frac{\tilde{Y}_1}{p'_2}+\frac{\tilde{Y}_2}{p_2},\textrm{
i.e. }p_2q_2=(p_2/p'_2)\tilde{Y}_1+\tilde{Y}_2.
\label{eq:ell2line}\eeq

\begin{figure}[ht]
\centering
\input{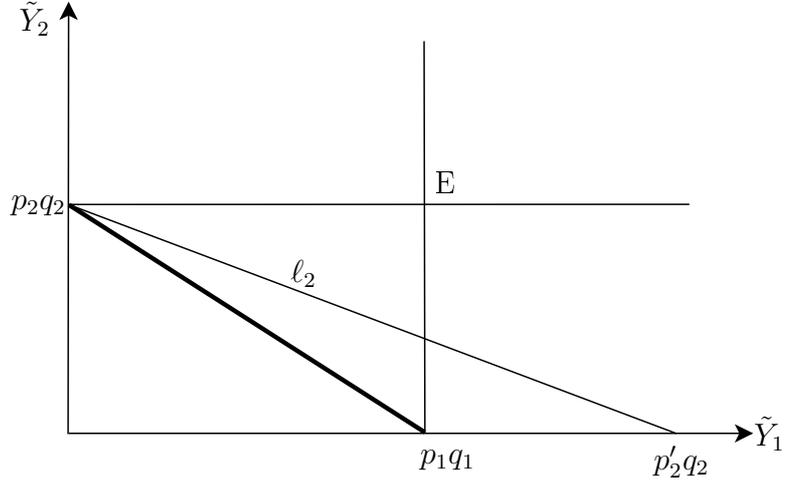}
\caption{The partition from Figure \ref{fig:income} when
$p'_2q_2>p_1q_1$.} \label{fig:incFig2}
\end{figure}

In view of \eqref{eq:assump} there are three possible cases:
\begin{subequations}\renewcommand{\theequation}{\theparentequation \roman{equation}}\begin{equation} p'_2q_2<p_1q_1 \quad \textrm{ (see Figure \ref{fig:income}),} \label{eq:4.4i} \end{equation}
\begin{equation} p'_2q_2=p_1q_1 \quad \textrm{ (the thick line in Figure \ref{fig:incFig2}),} \label{eq:4.4ii} \end{equation}
\begin{equation} p'_2q_2>p_1q_1 \quad \textrm{ (see Figure \ref{fig:incFig2}).} \label{eq:4.4iii} \end{equation}
\end{subequations}

The restrictions defining the cases $A_1$ and $A_2$ in the two
tables (again for maximal values of the variables $\alpha$ and
$\delta$) can be represented geometrically by the lines \beq
\ell_3:q_1=\frac{\tilde{Y}_1-p'_2\left(
q_2-\frac{\tilde{Y}_2}{p_2}\right)}{p_1},\textrm{ i.e.
}\frac{p_2}{p'_2}(p_1q_1+p'_2q_2)=\frac{p_2}{p'_2}\tilde{Y}_1+\tilde{Y}_2,
\label{eq:ell3line}\eeq \beq
\ell_4:q_2=\frac{\tilde{Y}_2-p'_1\left(
q_1-\frac{\tilde{Y}_1}{p_1} \right)}{p_2},\textrm{ i.e.
}\frac{p_1}{p'_1}(p'_1q_1+p_2q_2)=\tilde{Y}_1+\frac{p_1}{p'_1}\tilde{Y}_2
\label{eq:ell4line}\eeq

It is obvious that $\ell_1$ and $\ell_4$ are parallel, as are
$\ell_2$ and $\ell_3$. The point $E(p_1q_1,p_2q_2)$ lies on
$\ell_3$ and $\ell_4$. The role of these two lines is different in
zones II and IV, as will be described in more detail later (in
zone II it is only line $\ell_4$ that matters and line $\ell_3$
matters only in zone IV).

\begin{remark} To facilitate the verification of the statements for the
different cases, we start by advancing a comment on $\Delta
p_t:=p_{1,t}-p_{2,t}$, denoted for brevity simply as $\Delta p$.
Table \ref{tab:breply1} contains the cases \beq
\begin{array}{l l}
1) & \Delta p < \rho ,\\
2) & \Delta p = \rho ,\\
3) & \Delta p > \rho .\\
\end{array}\label{eq:A1}
\eeq

In an analogous manner, Table \ref{tab:breply2} contains the cases
defined by the conditions \beq
\begin{array}{l l}
1) & \Delta p > -\rho ,\\
2) & \Delta p = -\rho ,\\
3) & \Delta p < -\rho .\\
\end{array}\label{eq:A2}
\eeq

By combining case $i),~i=1,2,3,$ from \eqref{eq:A1} and case
$j),~j=1,2,3,$ from \eqref{eq:A2}, we obtain the following
subcases $(i,j)$: $$
  \begin{array}{l l}
    (1,1) & \Delta p \in (-\rho,\rho), \\
    (1,2) & \Delta p = -\rho, \\
    (1,3) & \Delta p < -\rho, \\
    (2,1) & \Delta p = \rho, \\
    (3,1) & \Delta p > \rho,
  \end{array}
$$ while subcases $(2,2),(2,3),(3,2)$ and $(3,3)$ are rendered
impossible by incompatible constraints. \endprf
\label{rem:Deltap}\end{remark}

We proceed to describe the evolution of prices for given initial
values $(p_{1,0},p_{2,0})$, for which financial resources
$(Y_1,Y_2)$ lie in zone III, i.e. \beq \frac{Y_1}{q_1}\geq
p_{1,0},~\frac{Y_2}{q_2}\geq p_{2,0}. \label{eq:zIIIinit}\eeq

Let $(\alpha,\beta,\gamma,\delta )$ be a NE for the chosen values
of parameters $Y_i,q_i,p_{i,0}$.

\textbf{1)} Suppose first that $q_2-\delta \leq 0$. Using Table
\ref{tab:breply1}, case I-$A$ (or 1-I-$A$ for brevity), we find
the possible range of values for $(\alpha,\beta)$. Then, by SR1 we
determine $\alpha=q_1,~\beta=0$. Now Table \ref{tab:breply2}, case
I-$A$ (2-I-$A$ for brevity) shows that $\delta=q_2$ and, according
to SR1, $\gamma=0$.

We note that for the given choice of the inequalities for the
different cases in the tables, from a formal standpoint sometimes
there arises the need to analyse other cases as well. (Case
1-I-$B$ is an example in the present situations.) In fact, we
could define the subcases in such a manner as to have only closed
sets. Our choice, however, has the advantage of simplifying the
exposition, while leading to the same results.

Thus the NE is $(q_1,0,0,q_2)$ and obviously the quantities
supplied are depleted. If \eqref{eq:zIIIinit} holds with
equalities, the financial resources are also depleted and no
adjustments in the prices are necessary as the point $E$ turns out
to be a p.e. with
$$p_{i,t}=p_{i,0}=\frac{Y_i}{q_i},~\forall t \geq 0.$$

If for some $i=1,2$ we have $$\frac{Y_i}{q_i}>p_{i,0},$$ then
according to \eqref{eq:priceup} the corresponding price is
adjusted upward to \beq
p_{i,1}=\frac{Y_i}{q_i}\label{eq:zIIIequi}\eeq and no further
adjustments are required. In other words, after one step prices
stabilize at the prices given by \eqref{eq:zIIIequi}:
$p_{i,t}=p_{i,1},~\forall t \geq 1.$

Next, we check for other NEs that may possibly be obtained in the
case given by \eqref{eq:zIIIinit}.

\textbf{2)} Suppose now that $0<q_2-\delta \leq Y_1/p'_2$. Then
using 1-II-$A_1$ we get $q_1-\alpha \leq 0$. Now 2-I-$A$ implies
$q_2-\delta \leq 0$, which is a contradiction under our
hypothesis.

The formula in 1-II-$A_2$ shows that $q_1-\alpha >0$ is possible
only for \beq p_{1,0}>p'_{2,0}.\label{eq:zIIIp1geqp'2}\eeq Then
2-II-$A_1$ implies $q_2-\delta \leq 0$, which is impossible, and
for 2-II-$A_2$ the only situation compatible with
\eqref{eq:zIIIp1geqp'2} is $p'_{1,0}>p_{2,0}$, which however
produces $\delta=q_2$, a contradiction with the assumption
$0<q_2-\delta \leq Y_1/p'_2$. It remains to check whether there is
an equilibrium for which $q_1-\alpha > Y_2/p'_1$. We see from
2-III-$A$ that when $p'_{1,0} \geq p_{2,0}$, then $q_2-\delta \leq
0$ is implied again. The case $p'_{1,0}<p_{2,0}$ is incompatible
with \eqref{eq:zIIIp1geqp'2}. This completes the analysis of the
case $0<q_2-\delta \leq Y_1/p'_2$.

\textbf{3)} Let us now assume that $$\frac{Y_1}{p'_2}<q_2-\delta
.$$ Then from 1-III-$A$ it follows that $q_1-\alpha>0$ only when
$p_{1,0}>p'_{2,0}$, since $\alpha=0$ then. One can verify as above
that the values from 2-II-$A_1$, 2-II-$A_2$ and 2-III-$A$ lead to
a
contradiction. 

In this case the two economies operate in autarky and reach
equilibrium without interacting with each other. This result is a
natural consequence of the fact that consumers in both economies
have local dominance (orders placed by the ``local'' consumer are
executed first), as well as sufficient financial resources to
absorb the entire local supply.

This proves \begin{proposition} For initial price values $p_{i,0}$
for which the financial resources $(Y_1,Y_2)$ are in zone III, as
defined in \eqref{eq:zIIIinit}, the unique NE is $(q_1,0,0,q_2)$
and after at most one (upward) price adjustment the equilibrium
point $E(p_{1,1}q_1,p_{2,1}q_2)$ is reached.
\label{prop:zoneIII}\end{proposition}

In what follows, we refer to the equilibrium described in
Proposition \ref{prop:zoneIII} as \textbf{equilibrium of type}
$\mathbf{E}$. This is the type of equilibrium arising in the case
of affluent economies that trade in conditions of ample financial
resource availability.

The analysis of the other cases is technically more complicated
and we list only the final results here, relegating sketches of
the proofs to the appendix.

We also note that it is possible for more than one NE to arise
depending on the relationship between prices and transportation
costs.

\begin{proposition} For initial prices $p_{i,0}$ for which $(Y_1,Y_2)$ is in zone
II, \beq Y_1 < p_{1,0}q_1,~p_{2,0}q_2 \leq Y_2 ,
\label{eq:zoneIIinit}\eeq which is divided into the following
subzones:

\noindent i) Zone II-1: \beq \left\{
  \begin{array}{l}
    Y_1< p_{1,0}q_1 ,~p_{2,0} q_2 \leq Y_2 , \\
    q_1 > \frac{Y_1}{p_{1,0}}+\frac{Y_2}{p'_{1,0}}\quad ( \textrm{strictly below
    }\ell_1),
  \end{array}
   \right.\label{eq:zoneII-1init}\eeq

\noindent ii) Zone II-2: \beq \left\{
  \begin{array}{l}
    Y_1 < p_{1,0}q_1,~p_{2,0}q_2\leq Y_2 , \\
    q_2 > \frac{1}{p_{2,0}}\left[ Y_2-p'_{1,0}\left(q_1-\frac{Y_1}{p_{1,0}}\right) \right]\quad ( \textrm{strictly below
    }\ell_4), \\
    \frac{Y_1}{p_{1,0}}+\frac{Y_2}{p'_{1,0}}\geq q_1 \quad (
    \textrm{on or above
    }\ell_1),
  \end{array}
 \right.
\label{eq:zoneII-2init}\eeq

\noindent iii) Zone II-3: \beq \left\{
  \begin{array}{l}
    Y_1< p_{1,0}q_1 ,~p_{2,0} q_2 < Y_2 , \\
    q_2 \leq \frac{1}{p_{2,0}}\left[ Y_2-p'_{1,0}\left(q_1-\frac{Y_1}{p_{1,0}}\right) \right]\quad ( \textrm{on or
    above
    }\ell_4),
  \end{array}
   \right. \label{eq:zoneII-3init}\eeq we have respectively:

\noindent a) in zone II-3 there exists a unique NE
$(\frac{Y_1}{p_1},0,q_1-\frac{Y_1}{p_1},q_2)$, for which either
$(Y_1,Y_2)\in\ell_{4,0}$, which is a p.e., or $(Y_1,Y_2)$ is
strictly above $\ell_{4,0}$ and we obtain $(Y_1,Y_2)\in\ell_{4,1}$
after one upward adjustment in $p_{2,0}$. (We refer to the p.e.s
of this type as \textbf{$\bs{\ell_4}$-equilibria}.)

\noindent b) in zone II-2 there exist two types of NE:

I) for $\Delta p < -\rho$: NE
$\left(\frac{Y_1}{p_{1,0}},0,q_1-\frac{Y_1}{p_{1,0}},\frac{Y_2-p'_{1,0}(q_1-Y_1/p_{1,0})}{p_{2,0}}
\right)$,

II) for $\Delta p \geq -\rho$: NE $\left(
\frac{Y_1}{p_{1,0}},0,\frac{Y_2-p_{2,0}q_2}{p'_{1,0}},q_2
\right)$.

For the two types of NE the following price adjustment patterns
obtain:

- \textbf{in case I)}: 1) when $(Y_1,Y_2)$ is strictly above
$\ell_{1,0}$, after one downward adjustment step in $p_{2,0}$ we
reach an $\ell_4$-equilibrium;

2) when $(Y_1,Y_2)$ lies on $\ell_{1,0}$, the price $p_{2,0}$ is
reduced to $p_{2,1}=0$ and we reach a \textbf{degenerate
$\bs{\ell_1}$-equilibrium}.

- \textbf{in case II)}: 1) when $Y_2=p_{2,0}q_2$, after one
downward adjustment in $p_{1,0}$ we reach a p.e. of type $E$;

2) when $Y_2>p_{2,0}q_2$, let $$k=\frac{1}{2q_1}\left[ \sqrt{(\rho
q_1-Y_1-Y_2+p_{2,0}q_2)^2 + 4q_1\rho Y_1} - (\rho
q_1-Y_1-Y_2+p_{2,0}q_2) \right]$$ and then, depending on whether

2.1) $p_{2,0}-\rho \leq k < p_{1,0}$

or

2.2) $k < p_{2,0}-\rho < p_{1,0}$

we have, respectively, in:

2.1) an infinite downward adjustment process in $p_{1,t}$ for
which $\lim_{t\rightarrow \infty}p_{1,t}=k$. (In this case the
system of two economies tends to a \textbf{degenerate
$\bs{\ell_4}$-equilibrium}, where the limit line is defined with
the aid of the number $k$ as $\ell_{4,\infty} : \tilde{Y}_1 +
\frac{k}{k+\rho}\tilde{Y}_2=k q_1+\frac{k}{k+\rho} p_{2,0}q_2$.)

2.2) after $s$ downward adjustments of $p_{1,t}, ~t=0,\ldots,s$,
we reach the situation described in case I). (Here the number
$s\in \mathbb{N}$ is determined by the condition
$$p_{2,0}-\rho \in \left[ g^s(p_{1,0}),g^{s-1}(p_{1,0}) \right],$$
where $g(x)=\frac{1}{q_1}\left[ Y_1 +
(Y_2-p_{2,0}q_2)\frac{x}{x+\rho}\right].$)

\noindent c) in zone II-1 there are two types of NE:

I) for $\Delta p < -\rho$: NE $\left(
\frac{Y_1}{p_{1,0}},0,\frac{Y_2}{p'_{1,0}},0 \right)$,

II) for $\Delta p \geq -\rho$: NE $\left(
\frac{Y_1}{p_{1,0}},0,\frac{Y_2-p_{2,0}q_2}{p'_{1,0}},q_2
\right)$.

For the two types of NE the following price adjustment patterns
obtain:

- \textbf{in case I)} we have an infinite downward adjustment
process in $p_{1,t}$, under which it tends monotonically to
$$p_{1,\infty}=\frac{1}{2q_2}\left( Y_1+Y_2-\rho q_1 + \sqrt{(Y_1+Y_2-\rho q_1)^2+4\rho q_1 Y_1}
\right),$$while $p_{2,t}=0,~\forall t \geq 1$. (In this case the
system of two economies tends to a degenerate
$\ell_1$-equilibrium, where the limit line is defined with the aid
of the number $p_{1,\infty}$ -- see \eqref{eq:A_ell1infty}.)

- \textbf{in case II)}: see case II) in b). \label{prop:zoneII}
\end{proposition}

Proposition \ref{prop:zoneII} deals with the interaction of an
affluent economy with abundant financial resources (region $II$)
and a relatively poor one (region $I$). In our setup ``affluence''
is defined in terms of the financial ability of consumers to
absorb the local (and, potentially, foreign) supply and is
unrelated to the production side of the economy. This allows for a
rich variety of situations in zone II. For instance, with very
high financial resources in region $II$, cheap output in region
$I$ and transportation costs that are not prohibitively high,
local consumers in region $I$ buy all they can afford, so that
consumers from region $II$ can absorb the residual supply in
region $I$, as well as the entire supply in their own region, and
still have income unspent. This naturally leads to a price
increase in the rich region, while prices in the poorer region are
unaffected by virtue of the pricing mechanism (see zone II-3 and
case a) above). As another example, if there is very abundant and
cheap supply in region $I$ (accounting for transportation costs in
the case of region $II$), the financial resources of the two
economies are entirely spent there and yet there remain unrealized
quantities, which keeps driving down the price in region $I$ (zone
II-1, case c)-I)). At the same time, the market in region $II$
becomes redundant and stops functioning, with a zero price
obtaining there and the entire amount of the good available being
consumed by the local consumer for free (off the market).

\begin{proposition} For initial prices $p_{i,0}$ for which $(Y_1,Y_2)$ is in zone
I, \beq 0<Y_1 < p_{1,0}q_1,~0< Y_2 <p_{2,0}q_2 ,
\label{eq:zoneIinit}\eeq which is divided into the following
subzones:

\noindent i) Zone I-1: \beq \left\{
  \begin{array}{l}
    0<Y_1< p_{1,0}q_1 ,~0< Y_2 <p_{2,0} q_2  , \\
    q_1 \leq \frac{Y_1}{p_{1,0}}+\frac{Y_2}{p'_{1,0}}\quad (
    \textrm{above or on
    }\ell_{1,0}),
  \end{array}
   \right.\label{eq:zoneI-1init}\eeq

\noindent ii) Zone I-2: \beq \left\{
  \begin{array}{l}
    0<Y_1 < p_{1,0}q_1,~0< Y_2 <p_{2,0}q_2 , \\
    \frac{Y_1}{p_{1,0}}+\frac{Y_2}{p'_{1,0}}<q_1  \quad ( \textrm{strictly below
    }\ell_{1,0}), \\
    q_2 \leq \frac{Y_1}{p'_{2,0}}+\frac{Y_2}{p_{2,0}} \quad (
    \textrm{on or above
    }\ell_{2,0}),
  \end{array}
 \right.
\label{eq:zoneI-2init}\eeq

\noindent iii) Zone I-3: \beq \left\{
  \begin{array}{l}
    0<Y_1< p_{1,0}q_1 ,~0< Y_2 <p_{2,0} q_2 , \\
     q_2 > \frac{Y_1}{p'_{2,0}}+\frac{Y_2}{p_{2,0}} \quad (
     \textrm{strictly below
    }\ell_{2,0}),
  \end{array}
\right. \label{eq:zoneI-3init}\eeq

we have respectively:

\noindent a) in zone I-1 there exist three types of NE:

I) for $\Delta p \in [-\rho,\rho]$: NE $\left(
\frac{Y_1}{p_{1,0}},0,0,\frac{Y_2}{p_{2,0}} \right)$,

II) for $\Delta p < -\rho$: NE $\left(
\frac{Y_1}{p_{1,0}},0,q_1-\frac{Y_1}{p_{1,0}},\frac{Y_2-p'_{1,0}\left(
q_1-\frac{Y_1}{p_{1,0}} \right)}{p_{2,0}} \right)$,

III) for $\Delta p > \rho$: NE $\left(
\frac{Y_1-p'_{2,0}\left(q_2-\frac{Y_2}{p_{2,0}}\right)}{p_{1,0}},q_2-\frac{Y_2}{p_{2,0}},0,\frac{Y_2}{p_{2,0}}
\right).$

For the three types of NE the following price adjustment patterns
obtain:

- \textbf{in case I)}: after one downward price adjustment we
reach a type $E$ equilibrium;

- \textbf{in case II)}: 1) when $(Y_1,Y_2)$ is strictly above
$\ell_1$, after one downward adjustment in $p_{2,0}$ we reach an
$\ell_4$-equilibrium;

2) when $(Y_1,Y_2)$ lies on $\ell_1$, we reach a degenerate
$\ell_1$-equilibrium ($p_{2,t}=0,~\forall t \geq 1$);

- \textbf{in case III)}: after one downward adjustment in
$p_{1,0}$, we reach an $\ell_3$-equilibrium (see Proposition
\ref{prop:zoneIV});

\noindent b) in zone I-2 there exist three types of NE with the
corresponding adjustment patterns:

I) for $\Delta p \in [-\rho,\rho]$: see I) in zone I-1,

II) for $\Delta p < -\rho$: see I) in zone II-1,

III) for $\Delta p > \rho$: see III) in zone I-1.

\noindent c) in zone I-3 there exist three types of NE with the
corresponding adjustment patterns:

I) for $\Delta p \in [-\rho,\rho]$: see I) in zone I-1,

II) for $\Delta p < -\rho$: see II) in zone I-2,

III) for $\Delta p > \rho$: NE $\left( 0,
\frac{Y_1}{p'_{2,0}},0,\frac{Y_2}{p_{2,0}}\right)$, in which case
after an infinite downward adjustment process for $p_{2,0}$, we
reach a \textbf{degenerate $\bs{\ell_2}$-equilibrium}.
\label{prop:zoneI}
\end{proposition}

Proposition \ref{prop:zoneI} analyzes the interaction of two
regions that are relatively poor in terms of initial wealth.
Naturally, the low purchasing power of the consumers in the two
regions results in deflationary developments, while the exact
distribution of consumption across regions also depends on the
size of transportation costs.

\begin{proposition} For initial prices $p_{i,0}$ for which $(Y_1,Y_2)$ is in zone
IV, \beq p_{1,0}q_1\leq Y_1 ,~0< Y_2 <p_{2,0}q_2 ,
\label{eq:zoneIVinit}\eeq which is divided into the following
subzones:

\noindent i) Zone IV-1: \beq \left\{
  \begin{array}{l}
    p_{1,0}q_1 \leq Y_1,~0< Y_2 <p_{2,0} q_2  , \\
    q_1 \leq \frac{Y_1-p'_{2,0}\left(q_2-\frac{Y_2}{p_{2,0}}\right)}{p_{1,0}}\quad (
    \textrm{above or on
    }\ell_{3,0}),
  \end{array}
   \right.\label{eq:zoneIV-1init}\eeq

\noindent ii) Zone IV-2: \beq \left\{
  \begin{array}{l}
    p_{1,0}q_1 \leq Y_1 ,~0< Y_2 <p_{2,0}q_2 , \\
    q_1 > \frac{Y_1-p'_{2,0}\left(q_2-\frac{Y_2}{p_{2,0}}\right)}{p_{1,0}}  \quad ( \textrm{strictly below
    }\ell_{3,0}),
  \end{array}
 \right.
\label{eq:zoneIV-2init}\eeq

we have respectively:

\noindent a) in zone IV-1 there exists a unique NE
$\left(q_1,q_2-\frac{Y_2}{p_{2,0}},0,\frac{Y_2}{p_{2,0}} \right)$,
which is symmetric (as regards a change of roles of the two
economies) to the NE from zone II-3 (see Proposition
\ref{prop:zoneII}, a)). The p.e. obtained in this case will be
referred to as an \textbf{$\bs{\ell_3}$-equilibrium}.

\noindent b) in zone IV-2 there exist two types of NE:

I) for $\Delta p \leq \rho$: NE $\left( q_1,
\frac{Y_1-p_{1,0}q_1}{p'_{2,0}},0,\frac{Y_2}{p_{2,0}}\right)$,

II) for $\Delta p > \rho$: NE $\left(\frac{Y_1-p'_{2,0}\left(
q_2-\frac{Y_2}{p_{2,0}} \right)}{p_{1,0}},
q_2-\frac{Y_2}{p_{2,0}},0,\frac{Y_2}{p_{2,0}}\right)$, \newline
which are symmetric (in the above sense) to cases II) and I) for
zone II-2 (see Proposition \ref{prop:zoneII}, b)).

\label{prop:zoneIV}
\end{proposition}

The results obtained for zone IV are symmetric to those for zone
II as regards a change of roles of the two economies. In this
case, region $I$ is the ``rich'' region and has the potential to
absorb a part of the supply in region $II$, while in the
``poorer'' region $II$ consumption is satisfied out of local
supply only.

\begin{proposition} For initial prices $p_{i,0}$ for which $(Y_1,Y_2)$ is in zone
I, defined by \eqref{eq:zoneIinit} under the condition
\eqref{eq:4.4iii} (see Figure \ref{fig:A5}), which is divided into
the following subzones:

\noindent i) zone 1-1: \beq \left\{
  \begin{array}{l}
    0< Y_1 <p_{1,0} q_1,~0< Y_2 <p_{2,0} q_2  , \\
    \frac{Y_1}{p_{1,0}}+\frac{Y_2}{p'_{1,0}} < q_1\quad (
    \textrm{strictly below }\ell_{1,0}), \\
    \frac{Y_1}{p'_{2,0}}+\frac{Y_2}{p_{2,0}} < q_2\quad (
    \textrm{strictly below }\ell_{2,0})
  \end{array}
   \right.\label{eq:zone1-1init}\eeq

\noindent ii) zone 1-2: \beq \left\{
  \begin{array}{l}
    0< Y_1 <p_{1,0} q_1,~0< Y_2 <p_{2,0} q_2  , \\
    \frac{Y_1}{p_{1,0}}+\frac{Y_2}{p'_{1,0}} \geq q_1\quad (
    \textrm{on or above }\ell_{1,0}), \\
    \frac{Y_1}{p'_{2,0}}+\frac{Y_2}{p_{2,0}} < q_2\quad (
    \textrm{strictly below }\ell_{2,0})
  \end{array}
   \right.\label{eq:zone1-2init}\eeq

\noindent iii) zone 1-3: \beq \left\{
  \begin{array}{l}
    0< Y_1 <p_{1,0} q_1,~0< Y_2 <p_{2,0} q_2  , \\
    \frac{Y_1}{p_{1,0}}+\frac{Y_2}{p'_{1,0}} < q_1\quad (
    \textrm{strictly below }\ell_{1,0}), \\
    \frac{Y_1}{p'_{2,0}}+\frac{Y_2}{p_{2,0}} \geq q_2\quad (
    \textrm{on or above }\ell_{2,0})
  \end{array}
   \right.\label{eq:zone1-3init}\eeq

\noindent iv) zone 1-4: \beq \left\{
  \begin{array}{l}
    0< Y_1 <p_{1,0} q_1,~0< Y_2 <p_{2,0} q_2  , \\
    \frac{Y_1}{p_{1,0}}+\frac{Y_2}{p'_{1,0}} \geq q_1\quad (
    \textrm{on or above }\ell_{1,0}), \\
    \frac{Y_1}{p'_{2,0}}+\frac{Y_2}{p_{2,0}} \geq q_2\quad (
    \textrm{on or above }\ell_{2,0})
  \end{array}
   \right.\label{eq:zone1-4init}\eeq

we have respectively:

\noindent a) in zone 1-4 the initial NEs and the respective price
adjustment processes coincide with those from zone I-1 (basic
case).

\noindent b) in zone 1-3 the initial NEs and the respective price
adjustment processes coincide respectively:

- for $\Delta p \geq -\rho$ - with those from zone I-1 (basic
case);

- for $\Delta p < -\rho$ - with those from zone I-2 (basic case).

\noindent c) in zone 1-2 the initial NEs and the respective price
adjustment processes coincide respectively:

- for $\Delta p \leq \rho$ - with those from zone 1-4;

- for $\Delta p > \rho$ - with those from zone I-3 (basic case).

\noindent d) in zone 1-1 the initial NEs and the respective price
adjustment processes coincide with those from zone I-3 (basic
case). \label{prop:zoneIagain}
\end{proposition}

Proposition \ref{prop:zoneIagain} revisits the analysis of the
interaction of two relatively poor regions in the special case
when the consumer in region $I$ needs more financial resources in
order to buy the entire supply in region $II$ than the resources
needed to entirely absorb local supply (condition
\eqref{eq:4.4iii}). Unsurprisingly, the results obtained replicate
the set of results from the basic case for zone I (Proposition
\ref{prop:zoneI}). The differences that arise are a natural
consequence of the different partitioning of zone 1 into subzones
due to the fact that the line $\ell_2$ now intersects the line
$\ell_1$ at the point $G$ (see Figure \ref{fig:A5}).

We conclude section \ref{sec:discprice} by formulating Theorem
\ref{thm:DiscPriceSummary}, which summarizes the results from
Propositions \ref{prop:zoneIII}-\ref{prop:zoneIagain}. Since we
have already stated the final results for the price dynamics
entailed by the model in discrete time, strictly accounting for
the various combinations of parameters possible, we now state the
theorem in a way that emphasizes the economic interpretation of
the results. For this purpose we introduce appropriate terms that
help illustrate the claims ($i=1,2$):
\begin{itemize}
  \item $RLS_i := q_i$ --- real local supply in the market in
  region $i$;
  \item $NLS_i := p_iq_i$ --- nominal local supply in region $i$ (valued at local
  prices);
  \item $NFR_i := Y_i$ --- nominal financial resources in region
  $i$;
  \item $TRFR_1 := \frac{Y_1}{p_1}+\frac{Y_2}{p'_1}$ --- total real
  financial resources, valued at region $I$'s prices;
  \item $TRFR_2 := \frac{Y_1}{p'_2}+\frac{Y_2}{p_2}$ --- total real
  financial resources, valued at region $II$'s prices;
  \item $T_1 := \frac{p_2q_2}{p'_1}+q_1$ --- total real
  supply in region $I$;
  \item $T_2 := \frac{p_1q_1}{p'_2}+q_2$ --- total real
  supply in region $II$.
\end{itemize}

With the help of the above quantities we can provide equivalent
formulations for the terms used in the different propositions:

\begin{flushleft}
  $\qquad$$\qquad$ $\left\{
  \begin{array}{l}
    \textrm{above} \\
    \textrm{on} \\
    \textrm{below}
  \end{array}
\right\} \ell_i \Leftrightarrow RLS_i\left\{
  \begin{array}{l}
    < \\
    = \\
    >
  \end{array}
\right\} TRFR_i,~i=1,2,$ \\

$\qquad$$\qquad$ $\left\{
  \begin{array}{l}
    \textrm{above} \\
    \textrm{on} \\
    \textrm{below}
  \end{array}
\right\} \ell_3 \Leftrightarrow T_2\left\{
  \begin{array}{l}
    < \\
    = \\
    >
  \end{array}
\right\} TRFR_2,$  \\

$\qquad$$\qquad$ $\left\{
  \begin{array}{l}
    \textrm{above} \\
    \textrm{on} \\
    \textrm{below}
  \end{array}
\right\} \ell_4 \Leftrightarrow T_1\left\{
  \begin{array}{l}
    < \\
    = \\
    >
  \end{array}
\right\} TRFR_1 .$
\end{flushleft}

\begin{theorem} I) When:
\begin{enumerate}
  \item[1)] $NFR_i > NLS_i,~i=1,2$ or
  \item[2)] $NFR_i \leq NLS_i,~i=1,2$ and at the same time $\Delta p \in
  [-\rho,\rho]$,
\end{enumerate} then with a one-time increase (case 1)) or decrease (case
2)) in prices we reach an equilibrium of type E (see zones I and
III).

II) When $NFR$ for one of the economies is less than the
respective $NLS$ but at the same time $TRFR$ valued at the local
price for this economy is not less than $T$, then we reach an
$\ell_3$-equilibrium or $\ell_4$-equilibrium after a one-time
\emph{increase} of the price in the \emph{other} economy (see
zones IV-1 and II-3).

III-a) When the second requirement in II) is violated (i.e. $T >
TRFR$) but we have \begin{enumerate}
  \item[i)] $RLS < TRFR < T$ and
  \item[ii)] the local price, adjusted for transportation costs,
  is strictly smaller than the price on the other market,
\end{enumerate} then the same result as in case II) obtains
through a \emph{decrease} of the latter price. (See zone II-3 for
the $\ell_4$-equilibrium, and zones IV-2 and I-1 for the
$\ell_3$-equilibria).

III-b) When condition ii) in III-a) is replaced by the opposite
condition, there are two situations:

III-b-1) The difference between the price in the other region and
the transportation costs does not exceed a critical threshold (the
number $k$ in Proposition \ref{prop:zoneII} for zone II-2);

III-b-2) The above difference is strictly greater than the
critical threshold.

Then, in case III-b-1) the system of the two economies tends to a
degenerate $\ell_4$-equilibrium through an infinite adjustment
process and in case III-b-2) after a finite number of steps a
regular $\ell_4$-equilibrium is reached.

IV) When in economy $i$ we have $RLS_i=TRFR_i$, we reach a
degenerate equilibrium in which the price in the other economy
immediately falls to zero (see the degenerate $\ell_1$-equilibrium
in the case in Figure \ref{fig:income} when \eqref{eq:4.4i}
holds).

V) When for economy $i$ we have $TRFR_i < q_i$ and condition ii)
from III-a) holds, then the two economies tend to a degenerate
equilibrium through a gradual decrease of the price which does not
automatically become zero (see zones II-1, I-2 and I-3 for the
$\ell_1$-equilibria, and zone I-3 for the $\ell_2$-equilibrium).
\label{thm:DiscPriceSummary}
\end{theorem}

\section{Price dynamics in the continuous-time case}\label{sec:contprice}

This section investigates the counterpart of the model in section
\ref{sec:model} under continuous-time price dynamics.  Similarly
to the setup described above, we formulate a price adjustment rule
on the basis of the residual income left unspent or the quantity
of the good not consumed at each instant $t$. This takes the form
of a system of ordinary differential equations, whose properties
are studied and compared to those of their discrete time
counterpart \eqref{eq:pricedown}-\eqref{eq:priceup}.

The problem setup and notation employed are identical to the ones
in section \ref{sec:model}. The static games played and all their
properties are the same as before, with the obvious difference
that the games are indexed by a set with the cardinality of the
continuum. To distinguish the continuous-time nature of the
present setup, we write the two prices as $p_{i}(t)$, $i=1,2$.

Thus, at time $t$ consumer $I$'s strategy space $S_1$ is
determined by the budget constraint and the nonnegativity
restrictions on the orders: \beq S_1=\{(\alpha,\beta)\in
\mathbb{R}^2_+~|~\alpha p_{1}(t)+ \beta (p_{2}(t)+\rho)\leq
Y_1\}.\label{eq:ss1cont}\eeq Consumer $II$'s strategy space in
period $t$ is \beq S_2=\{(\gamma,\delta)\in
\mathbb{R}^2_+~|~\gamma (p_{1}(t)+\rho)+ \delta p_{2}(t)\leq Y_2
\}.\label{eq:ss2cont}\eeq As before, we adopt the shorthand
$p'_{1}(t):= p_{1}(t)+\rho$ and $p'_{2}(t):=p_{2}(t)+\rho$. We
also omit the argument $t$ whenever it is evident from the context
or irrelevant.

The payoff (or utility) functions for consumers $I$ and $II$ are
denoted $U_1(\alpha,~\beta,~\gamma,~\delta)$ and
$U_2(\alpha,~\beta,~\gamma,~\delta)$, and are defined as in
\eqref{eq:payoff1} and \eqref{eq:payoff2}. Apart from the familiar
notation $q_i^{cons}$ and $Y_i^{res}$, $i=1,2$, we also define the
part of the instantaneous income flow for consumer $I$ that has
been spent as $Y_1^{cons}:=p_1(t)\alpha_0+p'_2(t)\beta_0\leq Y_1$.
The respective variable for consumer $II$ is
$Y_2^{cons}:=p'_1(t)\gamma_0+p_2(t)\delta_0 \leq Y_2$.

We first establish that at any moment in time we can have exactly
one of the two situations described in the previous paragraph
(Lemma \ref{lem:2.2}). We then show that
$Y_i^{res}-p_iq_i=Y_i-Y_i^{cons}$ (Corollary \ref{cor:2.3}).
\begin{lemma} \label{lem:2.1} Let $(\alpha_0,\beta_0,\gamma_0,\delta_0)$ be a Nash equilibrium as above. Then \beq \label{eq:2.3}\alpha_0+\gamma_0\leq q_1\quad\textrm{and}\quad \beta_0+\delta_0\leq q_2. \eeq  \end{lemma}
\proof We first observe that \beq \label{eq:2.4} \alpha_0\leq q_1
\quad \textrm{and}\quad \delta_0\leq q_2. \eeq To see this,
assume, for instance, that $\alpha_0>q_1$. The latter implies
$U_1(\alpha_0,\beta_0,\gamma_0,\delta_0)=U_1(q_1,\beta_0,\gamma_0,\delta_0)$.
When $p_1>0$, this contradicts SR1. If $p_1=0$, the claim follows
from SR3.

Taking into account \eqref{eq:2.4}, we obtain
\begin{equation*}\begin{split}
U_1(\alpha_0,\beta_0,\gamma_0,\delta_0)= &
\alpha_0+\min(\beta_0,q_2-\delta_0)\\
U_2(\alpha_0,\beta_0,\gamma_0,\delta_0)= &
\min(\gamma_0,q_1-\alpha_0)+\delta_0. \end{split}\end{equation*}
Now the first part in \eqref{eq:2.3} becomes obvious, as the
assumption $\gamma_0>q_1-\alpha_0$ contradicts SR1, applied to
$U_2$.
\endprf

\begin{lemma} It is impossible to have simultaneously \beq q_i^{cons}<q_i\quad \textrm{and}
\quad Y_i^{res}>p_iq_i,~i=1,2. \label{eq:2.5}\eeq
\label{lem:2.2}\end{lemma}

\proof Fix, for instance, $i=1$ and suppose the converse is true.
Then $$\alpha_0+\gamma_0=q_1^{cons}<q_1\quad \Rightarrow \quad
q_1-\alpha_0-\gamma_0>0.$$ Keeping $\beta_0,\gamma_0$ and
$\delta_0$ fixed, we increase $\alpha_0$ to
$\bar{\alpha}:=\alpha_0+\varepsilon,~\varepsilon \in
(0,q_1-\alpha_0-\gamma_0)$. This implies that
$$\alpha_0+\varepsilon<q_1-\gamma_0\quad \Rightarrow \quad p_1(\alpha_0+\varepsilon)<p_1q_1-p_1\gamma_0\leq p_1q_1<Y_1^{res},$$ which establishes the feasibility of $(\bar{\alpha},\beta_0,\gamma_0,\delta_0)$. Then
$U_1(\bar{\alpha},\beta_0,\gamma_0,\delta_0)>U_1(\alpha_0,\beta_0,\gamma_0,\delta_0)$,
which contradicts the assumption that
$(\alpha_0,\beta_0,\gamma_0,\delta_0)$ is a Nash equilibrium.
\endprf

\begin{cor}For consumer $I$ exactly one of the following alternatives is possible (with analogous results holding for consumer $II$):\begin{enumerate}
  \item[ i)] $\alpha_0=q_1$ and then
  $Y_1^{res}-p_1q_1=Y_1-Y_1^{cons}$;
  \item[ ii)] $\alpha_0<q_1$ and then $Y_1=Y_1^{cons}.$
\end{enumerate}
\label{cor:2.3}\end{cor}

\proof Since i) is obvious, we take up the case $\alpha_0<q_1$. By
Lemma \ref{lem:2.1} we have $\delta_0\leq q_2$ and we analyse two
cases:
\begin{enumerate}
  \item[I)] If $\delta_0=q_2$, then $\beta_0=0$. By definition
  $(\alpha_0,\beta_0)$ is the solution to \begin{equation*}\begin{split} & \max_{\alpha,\beta}U_1(\alpha,\beta,\gamma_0,\delta_0)\quad s.t. \\ & p_1\alpha+p'_2\beta\leq Y_1 \\ & \alpha,\beta\geq 0
  \end{split}\end{equation*}and since $\alpha_0<q_1$, we can focus
  on finding
  $\max_{\alpha,\beta}[\alpha+\min(\beta,0)]=\max_{\alpha}[\alpha]$
  subject to $\alpha\leq q_1$. Two subcases are possible:
\begin{enumerate}
  \item[a)]$\frac{Y_1}{p_1}\leq q_1$;
  \item[b)]$q_1<\frac{Y_1}{p_1}$.
\end{enumerate}

For subcase a) it is evident that $\alpha_0=\frac{Y_1}{p_1}$ and,
combined with $\beta_0=0$, this gives us
$Y_1^{cons}=p_1\alpha_0+p'_2\beta_0=Y_1$, as asserted. For subcase
b) there is no solution.
  \item[II)] If $\delta_0<q_2$, then Lemma \ref{lem:2.1} implies
  $\beta_0\leq q_2-\delta_0$. Again two subcases are possible:
  \begin{enumerate}
  \item[a)]$\frac{Y_1}{p'_2}\leq q_2-\delta_0$;
  \item[b)]$ q_2-\delta_0<\frac{Y_1}{p'_2}$.
\end{enumerate}
In maximizing $U_1(\alpha,\beta,\gamma,\delta)$ over $\{
(\alpha,\beta)|\alpha,\beta\geq 0,~p_1\alpha+p'_2\beta\leq Y_1\}$
we can restrict our attention to the intersection of this feasible
set with the set defined by $\alpha\leq q_1$ and $\beta \leq
q_2-\delta_0$, as the extremum $(\alpha_0,\beta_0)$ satisfies
these constraints as well. Then it is easily seen that either the
maximum is attained at a point along the budget constraint and
therefore $Y_1^{cons}=p_1\alpha_0+p'_2\beta_0=Y_1$, or (depending
on the magnitudes of $q_1$ and $Y_1/p_1$, and the slope of the
budget constraint line) we get $\alpha_0=q_1$, which violates the
initial assumptions.  \endprf
\end{enumerate}

\begin{remark}The claim in Corollary \ref{cor:2.3} can also be established through direct verification
by using the specific form of the NEs (see section
\ref{sec:exist}). \label{rem:DirectChk}\end{remark}

The price adjustment rule in discrete time is of the form (see
\eqref{eq:pricedown} and \eqref{eq:priceup}) \beq \label{eq:2.5a}
\frac{p_{t+1}-p_t}{p_t}=A(t)\times 1,\eeq i.e. $A(t)$ is related
to the change in the price for one time period. If we assume that
in the continuous-time case $A(t)$ does not change substantially
over a short time interval $[t,t+\Delta t]$, the counterpart of
the discrete-time adjustment rule will be
$$\frac{\Delta p(t)}{p(t)}=A(t)\Delta t,\quad \Delta
p(t):=p(t+\Delta t)-p(t).$$ Taking the limit in the above as
$\Delta t \rightarrow 0$, we obtain the differential equation
$$\frac{\dot{p}(t)}{p(t)}=A(t).$$

More precisely, the continuous-time counterpart of the adjustment
rules defined in \eqref{eq:pricedown} and \eqref{eq:priceup} is
given by the differential equation system \beq
\label{eq:2.6}\frac{\dot{p}_i}{p_i}=-\frac{q_i-q_i^{cons}}{q_i}+\frac{Y_i-Y_i^{cons}}{p_iq_i},
\quad i=1,2. \eeq For brevity we will employ the shorthand $Q_i$
for the right-hand side of equation \eqref{eq:2.6}. By virtue of
the results established above at most one of the terms on the
right-hand side of \eqref{eq:2.6} will be nonzero. To allow for
the possibility of the prices taking zero values, we rewrite the
above system as \beq \label{eq:2.7} \dot{p}_i=Q_ip_i,\quad
i=1,2.\eeq

For the problem at hand it is more convenient to switch to a
$p_1$-$p_2$ coordinate system instead of the $Y_1$-$Y_2$ system
used so far.

We obtain the following results: \begin{itemize}
\item[a)] The lines $Y_i=p_i q_i$ are transformed into the lines
$p_i=\frac{Y_1}{q_i}$, $i=1,2$, and we have the the conditions
$$Y_i \left\{
  \begin{array}{l}
    < \\
    = \\
    >
  \end{array}
 \right\} p_i q_i \quad \Longleftrightarrow \quad p_i \left\{
  \begin{array}{l}
    > \\
    = \\
    <
  \end{array}
 \right\} \frac{Y_i}{q_i} .$$

\item[b)] The equation $$\ell_1 :
q_1=\frac{Y_1}{p_1}+\frac{Y_2}{p'_1}$$ can be written in
equivalent form as $$q_1 p_1^2 +(q_1 \rho -Y_1-Y_2)p_1 - Y_1\rho =
0 $$ and, in view of the fact that $p_1$ is a price, we can take
only the positive root $$p_1^*=\frac{Y_1+Y_2-\rho
q_1+\sqrt{(Y_1+Y_2-\rho q_1)^2+4\rho q_1Y_1}}{2q_1}$$ and write
the last equation as
$$p_1=p^*_1 .$$ We also have $$(Y_1,Y_2)\textrm{ is } \left\{
  \begin{array}{l}
    \textrm{above} \\
    \textrm{on} \\
    \textrm{below}
  \end{array}
 \right\} \ell_1 \quad \Longleftrightarrow \quad
 p_1 \left\{
  \begin{array}{l}
    < \\
    = \\
    >
  \end{array}
 \right\} p^*_1 .$$

\item[c)] Analogously, the line $$\ell_2 :
q_2=\frac{Y_1}{p'_2}+\frac{Y_2}{p_2}$$ is transformed into
$$p_2=p^*_2, $$ where $$p_2^*=\frac{Y_1+Y_2-\rho
q_2+\sqrt{(Y_1+Y_2-\rho q_2)^2+4\rho q_2Y_2}}{2q_2}.$$ Moreover,
$$(Y_1,Y_2)\textrm{ is } \left\{
  \begin{array}{l}
    \textrm{above} \\
    \textrm{on} \\
    \textrm{below}
  \end{array}
 \right\} \ell_2 \quad \Longleftrightarrow \quad
 p_2 \left\{
  \begin{array}{l}
    < \\
    = \\
    >
  \end{array}
 \right\} p^*_2 .$$ Note that \beq p^*_1>\frac{Y_1}{q_1}\textrm{  and  } p^*_2>\frac{Y_2}{q_2}.
 \label{eq:pstarineq}\eeq For instance, it is easily verified that the first inequality in
 \eqref{eq:pstarineq} is equivalent to $$\sqrt{(Y_1+Y_2-\rho q_1)^2+4\rho q_1Y_1} > Y_1-Y_2+\rho
q_1 .$$ It is evidently true when the right-hand side is
non-positive. When the right-hand side is positive, we can square
the inequality and check that it is equivalent to $Y_2>0$.

\item[d)] Solving the equation for $\ell_3$ with respect to $p_1$ (for fixed
$Y_i$, $q_i$), we obtain the hyperbola $$ p_1=\frac{1}{q_1}\left(
Y_1+Y_2-\rho q_2-p_2q_2+\frac{\rho}{p_2}Y_2\right)=:h_3(p_2). $$
We note that

\item[d1)]$(Y_1,Y_2)\textrm{ is } \left\{
  \begin{array}{l}
    \textrm{above} \\
    \textrm{on} \\
    \textrm{below}
  \end{array}
 \right\} \ell_3  \quad \Longleftrightarrow \quad
 p_1 \left\{
  \begin{array}{l}
    < \\
    = \\
    >
  \end{array}
 \right\} h_3(p_2) \quad \Longleftrightarrow \quad $

 $\Longleftrightarrow \quad (p_1,p_2) \textrm{ is } \left\{
  \begin{array}{l}
    \textrm{below} \\
    \textrm{on} \\
    \textrm{above}
  \end{array}
 \right\} \textrm{ the graph of } p_1=h_3(p_2). $

\item[d2)] The hyperbola $p_1=h_3(p_2)$ crosses the $p_2$ axis at
the point $(0,p^*_2)$ and therefore has the form shown in Figure
\ref{fig:h3hyperb}.

\begin{figure}[ht]
\centering
\input{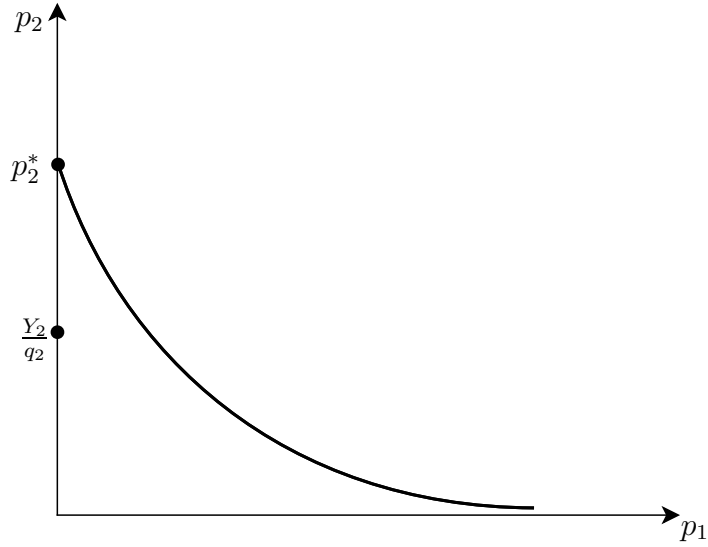}
\caption{The hyperbola $p_1=h_3(p_2)$.} \label{fig:h3hyperb}
\end{figure}

\item[e)] In a similar manner, $\ell_4$ is transformed into the hyperbola $$p_2=\frac{1}{q_2}\left( Y_1+Y_2-\rho
q_1-p_1q_1+\frac{\rho}{p_1}Y_1 \right) =: h_4(p_1). $$

Then, $$(Y_1,Y_2)\textrm{ is } \left\{
  \begin{array}{l}
    \textrm{above} \\
    \textrm{on} \\
    \textrm{below}
  \end{array}
 \right\} \ell_4  \quad \Longleftrightarrow \quad
 p_2 \left\{
  \begin{array}{l}
    < \\
    = \\
    >
  \end{array}
 \right\} h_4(p_1) \quad \Longleftrightarrow \quad $$  $$\Longleftrightarrow \quad (p_1,p_2) \textrm{ is } \left\{
  \begin{array}{l}
    \textrm{below} \\
    \textrm{on} \\
    \textrm{above}
  \end{array}
 \right\} \textrm{ the graph of } p_2=h_4(p_1). $$

Moreover, the hyperbola $p_2=h_4(p_1)$ crosses the $p_1$ axis at
the point $(p^*_1,0)$ and therefore has the form shown in Figure
\ref{fig:h4hyperb}.

\begin{figure}[ht]
\centering
\input{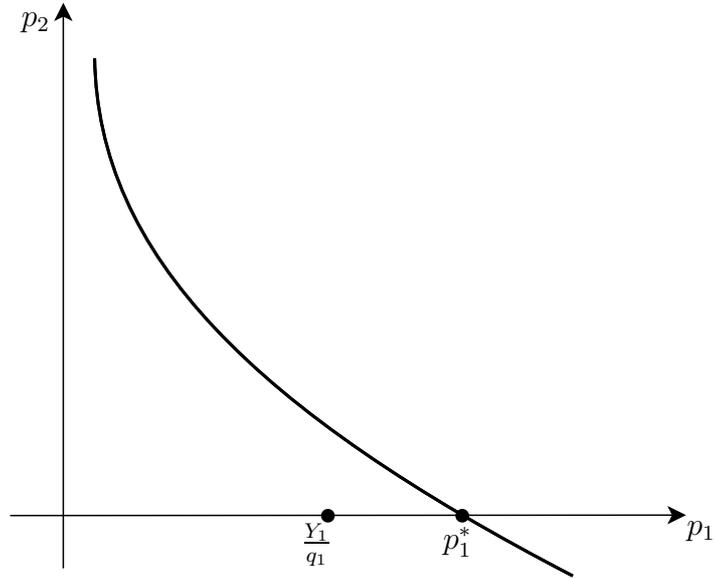}
\caption{The hyperbola $p_2=h_4(p_1)$.} \label{fig:h4hyperb}
\end{figure}

\end{itemize}

\begin{figure}[ht]
\centering
\input{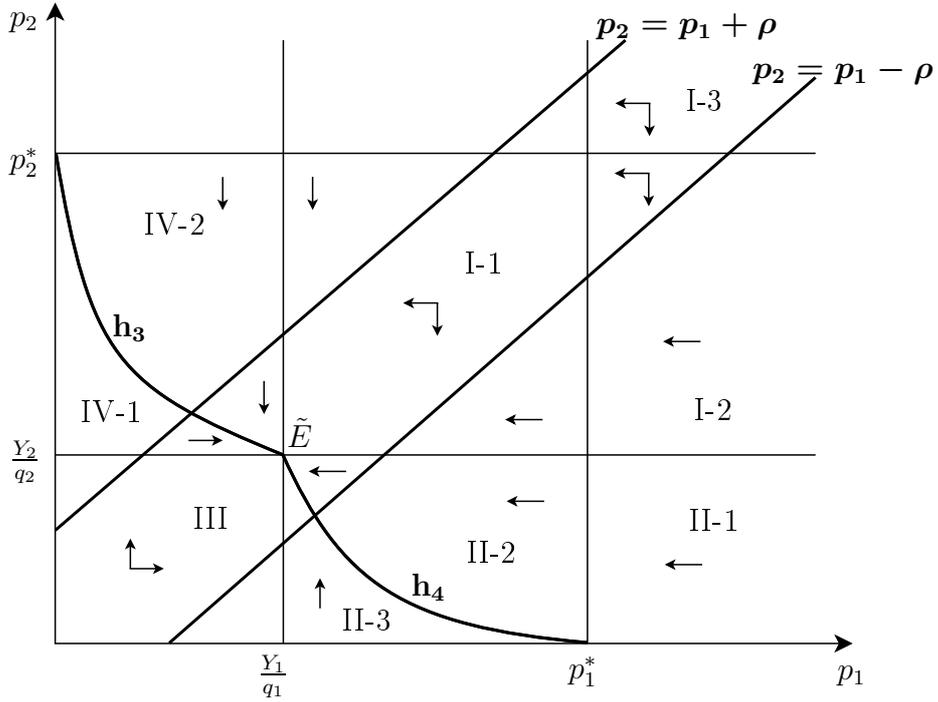}
\caption{The partition from Figure \ref{fig:income} in $p_1$-$p_2$ space.} \label{fig:pricespace1}
\end{figure}

Thus, we obtain Figure \ref{fig:pricespace1}, which is the
equivalent of Figure \ref{fig:income} in $p_1$-$p_2$ space.

For convenience we list the form of the NEs as described in the
discrete-time case (see Tables \ref{tab:Nash1} and
\ref{tab:Nash2}).

\begin{table}[hp]\centering
\begin{tabular}{|l|c|c|}
\hline
{\bf Zone} & {\bf Price relation} & {\bf Nash equilibrium }$\bs{(\alpha_0,\beta_0,\gamma_0,\delta_0)}$ \\
\hline
III & \multicolumn{1}{l|}{-} & $(q_1,0,0,q_2)$ \\
\hline
II-3 & \multicolumn{1}{l|}{-} & $\left(\frac{Y_1}{p_1},0,q_1-\frac{Y_1}{p_1},q_2\right)$ \\
\hline
II-2 & \multicolumn{1}{l|}{$\Delta p \geq -\rho$} & $\left(\frac{Y_1}{p_1},0,\frac{Y_2-p_2q_2}{p'_1},q_2\right)$ \\
\cline{2-3}
 & \multicolumn{1}{l|}{$\Delta p < -\rho$} & $\left(\frac{Y_1}{p_1},0,q_1-\frac{Y_1}{p_1},\frac{Y_2-p'_1\left(q_1-\frac{Y_1}{p_1}\right)}{p_2}\right)$ \\
\hline
II-1 & \multicolumn{1}{l|}{$\Delta p \geq -\rho$} & $\left(\frac{Y_1}{p_1},0,\frac{Y_2-p_2q_2}{p'_1},q_2\right)$ \\
\cline{2-3}
 & \multicolumn{1}{l|}{$\Delta p < -\rho$} & $\left(\frac{Y_1}{p_1},0,\frac{Y_2}{p'_1},0\right)$ \\
\hline
I-1 & \multicolumn{1}{l|}{$\Delta p \in [-\rho,\rho]$} & $\left(\frac{Y_1}{p_1},0,0,\frac{Y_2}{p_2}\right)$ \\
\cline{2-3}
 & \multicolumn{1}{l|}{$\Delta p < -\rho$} & $\left(\frac{Y_1}{p_1},0,q_1-\frac{Y_1}{p_1},\frac{Y_2-p'_1\left(q_1-\frac{Y_1}{p_1}\right)}{p_2}\right)$ \\
\cline{2-3}
 & \multicolumn{1}{l|}{$\Delta p >\rho$} & $\left(\frac{Y_1-p'_2\left(q_2-\frac{Y_2}{p_2}\right)}{p_1},q_2-\frac{Y_2}{p_2},0,\frac{Y_2}{p_2}\right)$ \\
\hline
I-2 & \multicolumn{1}{l|}{$\Delta p \in [-\rho,\rho]$} & $\left(\frac{Y_1}{p_1},0,0,\frac{Y_2}{p_2}\right)$ \\
\cline{2-3}
 & \multicolumn{1}{l|}{$\Delta p < -\rho$} & $\left(\frac{Y_1}{p_1},0,\frac{Y_2}{p'_1},0\right)$ \\
\cline{2-3}
 & \multicolumn{1}{l|}{$\Delta p >\rho$} & $\left(\frac{Y_1-p'_2\left(q_2-\frac{Y_2}{p_2}\right)}{p_1},q_2-\frac{Y_2}{p_2},0,\frac{Y_2}{p_2}\right)$ \\
\hline
I-3 & \multicolumn{1}{l|}{$\Delta p \in [-\rho,\rho]$} & $\left(\frac{Y_1}{p_1},0,0,\frac{Y_2}{p_2}\right)$ \\
\cline{2-3}
 & \multicolumn{1}{l|}{$\Delta p < -\rho$} & $\left(\frac{Y_1}{p_1},0,\frac{Y_2}{p'_1},0\right)$ \\
\cline{2-3}
 & \multicolumn{1}{l|}{$\Delta p >\rho$} & $\left(0,\frac{Y_1}{p'_2},0,\frac{Y_2}{p_2}\right)$ \\
\hline
IV-1 & \multicolumn{1}{l|}{-} & $\left(q_1,q_2-\frac{Y_2}{p_2},0,\frac{Y_2}{p_2}\right)$ \\
\hline
IV-2 & \multicolumn{1}{l|}{$\Delta p \leq \rho$} & $\left(q_1,\frac{Y_1-p_1q_1}{p'_2},0,\frac{Y_2}{p_2}\right)$ \\
\cline{2-3}
 & \multicolumn{1}{l|}{$\Delta p > \rho$} & $\left(\frac{Y_1-p'_2\left(q_2-\frac{Y_2}{p_2}\right)}{p_1},q_2-\frac{Y_2}{p_2},0,\frac{Y_2}{p_2}\right)$ \\
\hline
\end{tabular}
\caption{Nash equilibria for the case $p'_2q_2<p_1q_1$.}
\label{tab:Nash1}
\end{table}

\begin{table}[hp]
  \centering\begin{tabular}{|l|c|c|}
\hline
{\bf Zone} & {\bf Price relation} & {\bf Nash equilibrium }$\bs{(\alpha_0,\beta_0,\gamma_0,\delta_0)}$ \\
\hline
III & \multicolumn{1}{l|}{-} & $(q_1,0,0,q_2)$ \\
\hline
II-3 & \multicolumn{1}{l|}{-} & $\left(\frac{Y_1}{p_1},0,q_1-\frac{Y_1}{p_1},q_2\right)$ \\
\hline
II-2 & \multicolumn{1}{l|}{$\Delta p \geq -\rho$} & $\left(\frac{Y_1}{p_1},0,\frac{Y_2-p_2q_2}{p'_1},q_2\right)$ \\
\cline{2-3}
 & \multicolumn{1}{l|}{$\Delta p < -\rho$} & $\left(\frac{Y_1}{p_1},0,q_1-\frac{Y_1}{p_1},\frac{Y_2-p'_1\left(q_1-\frac{Y_1}{p_1}\right)}{p_2}\right)$ \\
\hline
II-1 & \multicolumn{1}{l|}{$\Delta p \geq -\rho$} & $\left(\frac{Y_1}{p_1},0,\frac{Y_2-p_2q_2}{p'_1},q_2\right)$ \\
\cline{2-3}
 & \multicolumn{1}{l|}{$\Delta p < -\rho$} & $\left(\frac{Y_1}{p_1},0,\frac{Y_2}{p'_1},0\right)$ \\
\hline
1-1 & \multicolumn{1}{l|}{$\Delta p \in [-\rho,\rho]$} & $\left(\frac{Y_1}{p_1},0,0,\frac{Y_2}{p_2}\right)$ \\
\cline{2-3}
 & \multicolumn{1}{l|}{$\Delta p < -\rho$} & $\left(\frac{Y_1}{p_1},0,\frac{Y_2}{p'_1},0\right)$ \\
\cline{2-3}
 & \multicolumn{1}{l|}{$\Delta p >\rho$} & $\left(0,\frac{Y_1}{p'_2},0,\frac{Y_2}{p_2}\right)$ \\
\hline
1-2 & \multicolumn{1}{l|}{$\Delta p \in [-\rho,\rho]$} & $\left(\frac{Y_1}{p_1},0,0,\frac{Y_2}{p_2}\right)$ \\
\cline{2-3}
 & \multicolumn{1}{l|}{$\Delta p < -\rho$} & $\left(\frac{Y_1}{p_1},0,q_1-\frac{Y_1}{p_1},\frac{Y_2-p'_1\left(q_1-\frac{Y_1}{p_1}\right)}{p_2}\right)$ \\
\cline{2-3}
 & \multicolumn{1}{l|}{$\Delta p >\rho$} & $\left(0,\frac{Y_1}{p'_2},0,\frac{Y_2}{p_2}\right)$ \\
\hline
1-3 & \multicolumn{1}{l|}{$\Delta p \in [-\rho,\rho]$} & $\left(\frac{Y_1}{p_1},0,0,\frac{Y_2}{p_2}\right)$ \\
\cline{2-3}
 & \multicolumn{1}{l|}{$\Delta p < -\rho$} & $\left(\frac{Y_1}{p_1},0,\frac{Y_2}{p'_1},0\right)$ \\
\cline{2-3}
 & \multicolumn{1}{l|}{$\Delta p >\rho$} & $\left(\frac{Y_1-p'_2\left(q_2-\frac{Y_2}{p_2}\right)}{p_1},q_2-\frac{Y_2}{p_2},0,\frac{Y_2}{p_2}\right)$ \\
\hline
1-4 & \multicolumn{1}{l|}{$\Delta p \in [-\rho,\rho]$} & $\left(\frac{Y_1}{p_1},0,0,\frac{Y_2}{p_2}\right)$ \\
\cline{2-3}
 & \multicolumn{1}{l|}{$\Delta p < -\rho$} & $\left(\frac{Y_1}{p_1},0,q_1-\frac{Y_1}{p_1},\frac{Y_2-p'_1\left(q_1-\frac{Y_1}{p_1}\right)}{p_2}\right)$ \\
\cline{2-3}
 & \multicolumn{1}{l|}{$\Delta p >\rho$} & $\left(\frac{Y_1-p'_2\left(q_2-\frac{Y_2}{p_2}\right)}{p_1},q_2-\frac{Y_2}{p_2},0,\frac{Y_2}{p_2}\right)$ \\
\hline
IV-1 & \multicolumn{1}{l|}{-} & $\left(q_1,q_2-\frac{Y_2}{p_2},0,\frac{Y_2}{p_2}\right)$ \\
\hline
IV-2 & \multicolumn{1}{l|}{$\Delta p \leq \rho$} & $\left(q_1,\frac{Y_1-p_1q_1}{p'_2},0,\frac{Y_2}{p_2}\right)$ \\
\cline{2-3}
 & \multicolumn{1}{l|}{$\Delta p > \rho$} & $\left(\frac{Y_1-p'_2\left(q_2-\frac{Y_2}{p_2}\right)}{p_1},q_2-\frac{Y_2}{p_2},0,\frac{Y_2}{p_2}\right)$ \\
\hline
\end{tabular}
\caption{Nash equilibria for the case
$p'_2q_2>p_1q_1$.}\label{tab:Nash2}
\end{table}

Now we can compute the right-hand sides $p_i Q_i$ of equations
\eqref{eq:2.7}, which are shown in Table \ref{tab:rhs2.7}.

\begin{table}[hp]\small
  \centering
\begin{tabular}{|l|c|c|l|c|l|}
\hline
{\bf Zone} & {\bf $\mathbf{p}$-$\bs{\rho}$ relation} & $\mathbf{p_1Q_1}$ & \multicolumn{1}{c|}{{\bf Sign}} & $\mathbf{p_2Q_2}$ & \multicolumn{1}{c|}{{\bf Sign}} \\
\hline
III & \multicolumn{1}{l|}{-} & $\frac{Y_1-p_1q_1}{q_1}$ & \multicolumn{1}{c|}{+} & $\frac{Y_2-p_2q_2}{q_2}$ & \multicolumn{1}{c|}{+} \\
\hline
II-3 & \multicolumn{1}{l|}{-} & $0$ & \multicolumn{1}{c|}{} & $\frac{Y_2-p_2q_2-p'_1\left ( q_1-\frac{Y_1}{p_1} \right )}{q_2}$ & \multicolumn{1}{c|}{+ below $h_4$} \\
\hline
II-2 & \multicolumn{1}{l|}{$\Delta p \geq -\rho$} & $\frac{Y_1}{q_1}+\frac{(Y_2-p_2q_2)p_1}{q_1p'_1}-p_1$ & \multicolumn{1}{c|}{- above $h_4$} & $0$ & \multicolumn{1}{c|}{} \\
\cline{2-6}
 & \multicolumn{1}{l|}{$\Delta p < -\rho$} & $0$ & \multicolumn{1}{c|}{} & $\frac{ Y_2-p'_1\left (q_1-\frac{Y_1}{p_1} \right)}{q_2}-p_2 $ & \multicolumn{1}{c|}{- above $h_4$} \\
\hline
II-1 & \multicolumn{1}{l|}{$\Delta p \geq -\rho$} & $\frac{Y_1}{q_1}+\frac{(Y_2-p_2q_2)p_1}{q_1p'_1}-p_1$ & \multicolumn{1}{c|}{- above $h_4$} & $0$ & \multicolumn{1}{c|}{} \\
\cline{2-6}
 & \multicolumn{1}{l|}{$\Delta p < -\rho$} & $\left ( \frac{Y_1}{p_1}+\frac{Y_2}{p'_1}-q_1 \right )\frac{p_1}{q_1}$ & \multicolumn{1}{c|}{- for $p_1>p_1^*$} & $-p_2$ & \multicolumn{1}{c|}{-} \\
\hline
I-1 & \multicolumn{1}{l|}{$\Delta p \in [-\rho,\rho]$} & $\frac{Y_1-p_1q_1}{q_1}$ & \multicolumn{1}{c|}{-} & $\frac{Y_2-p_2q_2}{q_2}$ & \multicolumn{1}{c|}{-} \\
\cline{2-6}
 & \multicolumn{1}{l|}{$\Delta p < -\rho$} & $0$ & \multicolumn{1}{c|}{} & $ \frac{ Y_2-p'_1\left (q_1-\frac{Y_1}{p_1} \right)}{q_2}-p_2 $ & \multicolumn{1}{c|}{- above $h_4$} \\
\cline{2-6}
 & \multicolumn{1}{l|}{$\Delta p >\rho$} & $\frac{ Y_1-p'_2\left (q_2-\frac{Y_2}{p_2} \right)}{q_1}-p_1 $ & \multicolumn{1}{c|}{- above $h_3$} & $0$ & \multicolumn{1}{c|}{} \\
\hline
I-2 & \multicolumn{1}{l|}{$\Delta p \in [-\rho,\rho]$} & $\frac{Y_1-p_1q_1}{q_1}$ & \multicolumn{1}{c|}{-} & $\frac{Y_2-p_2q_2}{q_2}$ & \multicolumn{1}{c|}{-} \\
\cline{2-6}
 & \multicolumn{1}{l|}{$\Delta p < -\rho$} & $\left ( \frac{Y_1}{p_1}+\frac{Y_2}{p'_1}-q_1 \right )\frac{p_1}{q_1}$ & \multicolumn{1}{c|}{- for $p_1>p_1^*$} & $-p_2$ & \multicolumn{1}{c|}{-} \\
\cline{2-6}
 & \multicolumn{1}{l|}{$\Delta p >\rho$} & $\frac{ Y_1-p'_2\left (q_2-\frac{Y_2}{p_2} \right)}{q_1}-p_1 $ & \multicolumn{1}{c|}{- above $h_3$} & $0$ & \multicolumn{1}{c|}{} \\
\hline
I-3 & \multicolumn{1}{l|}{$\Delta p \in [-\rho,\rho]$} & $\frac{Y_1-p_1q_1}{q_1}$ & \multicolumn{1}{c|}{-} & $\frac{Y_2-p_2q_2}{q_2}$ & \multicolumn{1}{c|}{-} \\
\cline{2-6}
 & \multicolumn{1}{l|}{$\Delta p < -\rho$} & $\left ( \frac{Y_1}{p_1}+\frac{Y_2}{p'_1}-q_1 \right )\frac{p_1}{q_1}$ & \multicolumn{1}{c|}{- for $p_1>p_1^*$} & $-p_2$ & \multicolumn{1}{c|}{-} \\
\cline{2-6}
 & \multicolumn{1}{l|}{$\Delta p >\rho$} & $-p_1$ & \multicolumn{1}{c|}{-} & $\left ( \frac{Y_1}{p'_2}+\frac{Y_2}{p_2}-q_2 \right )\frac{p_2}{q_2}$ & \multicolumn{1}{c|}{- for $p_2>p_2^*$} \\
\hline
IV-1 & \multicolumn{1}{l|}{-} & $\frac{Y_1-p_1q_1-p'_2\left (q_2-\frac{Y_2}{p_2} \right)}{q_1}$ & \multicolumn{1}{c|}{+ below $h_3$} & $0$ & \multicolumn{1}{c|}{} \\
\hline
IV-2 & \multicolumn{1}{l|}{$\Delta p \leq \rho$} & $0$ & \multicolumn{1}{c|}{} & $\frac{(Y_1-p_1q_1)p_2}{q_2p'_2}+\frac{Y_2}{q_2}-p_2 $ & \multicolumn{1}{c|}{- above $h_3$} \\
\cline{2-6}
 & \multicolumn{1}{l|}{$\Delta p > \rho$} & $\frac{ Y_1-p'_2\left (q_2-\frac{Y_2}{p_2} \right)}{q_1}-p_1 $ & \multicolumn{1}{c|}{- above $h_3$} & $0$ & \multicolumn{1}{c|}{} \\
\hline
\end{tabular}
\caption{Right-hand side expressions and signs for equations
\eqref{eq:2.7}.}\label{tab:rhs2.7}
\end{table}

\begin{figure}[ht]
\centering
\input{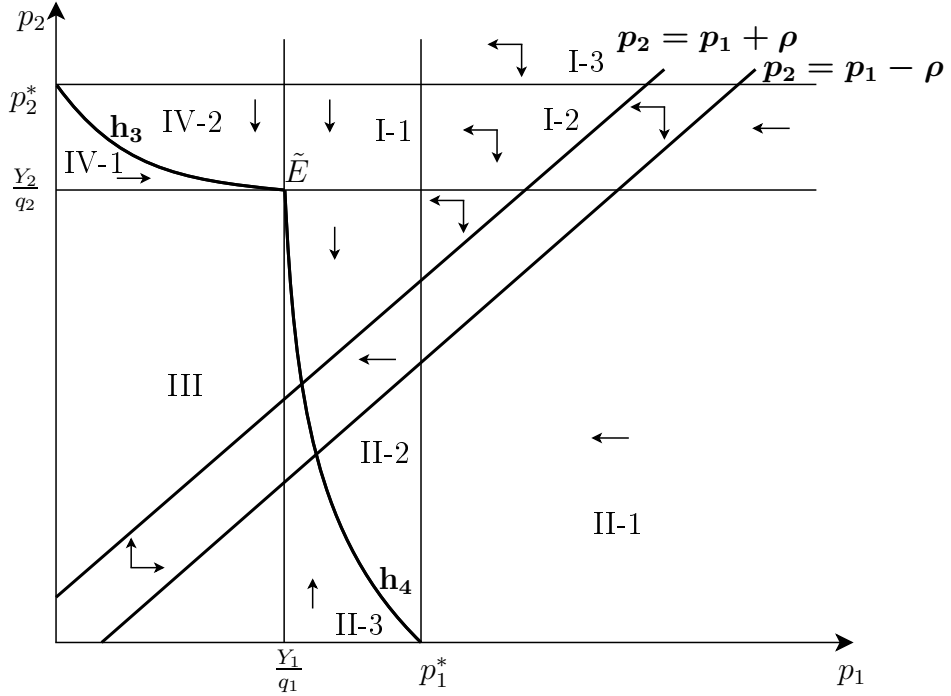}
\caption{The partition from Figure \ref{fig:income} in $p_1$-$p_2$
space for a different position of the point $\tilde{E}$ with
respect to the lines $p_2=p_1\pm \rho$.} \label{fig:pricespace2}
\end{figure}

The results in Table \ref{tab:rhs2.7} allow us to find the
direction of the phase flows shown in Figures
\ref{fig:pricespace1} and \ref{fig:pricespace2}. It is evident
that the position of the lines $p_2=p_1+\rho$ and $p_2=p_1-\rho$
relative to the partition in Figure \ref{fig:pricespace1} plays a
special role for the type of phase portrait obtained. For
instance, when the point $\tilde{E}$ lies in the set
$\{(p_1,p_2)|~ p_1-\rho<p_2<p_1+\rho \}$, we have a phase portrait
of the type shown in Fig. \ref{fig:pricespace1}. When $\tilde{E}$
is above the set $\{(p_1,p_2)|~ p_1-\rho<p_2<p_1+\rho \}$, the
situation shown in Figure \ref{fig:pricespace2} obtains. The
reader can easily produce phase portraits of this kind for various
assumptions about $Y_1,Y_2,q_1,q_2$ and $\rho$ with the help of
Table \ref{tab:rhs2.7}.

Table \ref{tab:rhs2.7} shows that, depending on the transportation
costs $\rho$ and the coordinates
$\left(\frac{Y_1}{q_1},\frac{Y_2}{q_2}\right)$ of point
$\tilde{E}$ (or, equivalently, depending on the position of this
point with respect to the set $\{(p_1,p_2)|~ p_1-\rho<p_2<p_1+\rho
\}$), it is possible to have discontinuities in the right-hand
side of \eqref{eq:2.7}. For example, in zone IV-2 it is possible
to have a discontinuity along the line $p_2=p_1-\rho$, provided
that $\tilde{E}$ lies below it, i.e.
$\frac{Y_2}{q_2}<\frac{Y_1}{q_1}-\rho$. In these cases we obtain
the situation in \cite[ pp. 41-42]{Fil85}, where the phase
trajectory, after hitting the surface of the discontinuity, stays
on it\footnote{For more details see \cite[p. 64 and pp.
82-83]{Fil85}}.

\begin{figure}[ht]
\centering
\input{Stab1.tex}
\caption{A neighborhood of the point $(p_1^0,p_2^0)$ on $h_3$ for
$(\bar{p}_1^0,\bar{p}_2^0)$ in zone IV-1.} \label{fig:stab1}
\end{figure}

The only equilibrium points for the system \eqref{eq:2.7} are
those on the hyperbolae $h_3$ and $h_4$, including point
$\tilde{E}$ (see Figure \ref{fig:pricespace1}). The direction of
the phase flow as presented in the figure makes it clear that the
fixed points we consider are Lyapunov stable but not
asymptotically stable. Let us take, for example, a point
$(p_1^0,p_2^0)$ on the graph of $h_3$, assuming that $\tilde{E}$
lies above the line $p_2=p_1-\rho$ (i.e. $\Delta p < \rho$). In
zone IV-1 the system takes the form \beq
\left\{\begin{array}{l}\dot{p_1}=\frac{1}{q_1}\left(Y_1+Y_2-\rho
q_2+\frac{\rho}{p_2}Y_2-p_2q_2-p_1q_1\right)\\\dot{p_2}=0\end{array}
\right. \label{eq:DS1}\eeq Fix a neighborhood $V$ of
$(p_1^0,p_2^0)$ and let $(\bar{p}_1^0,\bar{p}_2^0)$ be a point in
the intersection of $V$ and zone IV-1 (see Figure
\ref{fig:stab1}). In other words, this point is below the
hyperbola $p_1=h_3(p_2)$. Then the scalar
$$C\equiv \frac{1}{q_1}\left( Y_1+Y_2-\rho
q_2+\frac{\rho}{\bar{p}_2^0}Y_2-\bar{p}_2^0q_2
\right)-\bar{p}_1^0=h_3(\bar{p}_2^0)-\bar{p}_1^0>0$$ can be made
arbitrarily small if we shrink $V$ in an appropriate manner, since
$h_3(\cdot)$ is continuous and $h_3(p_2^0)-p^0_1=0$.

The second equation in \eqref{eq:DS1} implies $p_2=\bar{p}_2^0$
and therefore the first equation takes the form
$$\dot{p}_1=-p_1+C+\bar{p}_1^0,~p_1(0)=\bar{p}^0_1,$$ the
corresponding solution being $$p_1(t)=-Ce^{-t}+C+\bar{p}_1^0.$$

For any $t>0$ the distance between $(p_1(t),p_2(t))$ and
$(p_1^0,p_2^0)$ is bounded above by
$2C+|\bar{p}_1^0-p_1^0|+|\bar{p}_2^0-p_2^0|$. As $t \rightarrow
+\infty$, the solution $(p_1(t),p_2(t))$ tends to
$(C+\bar{p}_1^0,\bar{p}_2^0)=(h_3(\bar{p}_2^0),\bar{p}_2^0)\neq
(h_3(p_2^0),p_2^0)=(p_1^0,p_2^0)$, except in the special case when
$\bar{p}_2^0=p_2^0$.

In zone IV-2 (i.e. above the hyperbola $h_3$) and under the
condition $\Delta p \leq \rho$ we have the system \beq
\left\{\begin{array}{l}\dot{p_1}=0\\
\dot{p_2}=\frac{1}{q_2}(Y_1+Y_2-p_1 q_1+\rho q_2) +\rho\frac{p_1
q_1 -Y_1}{q_2}.\frac{1}{p'_2}-p'_2\equiv \tilde{H}(p_1,p_2).
\end{array} \right. \label{eq:DS2}\eeq

\begin{figure}[ht]
\centering
\input{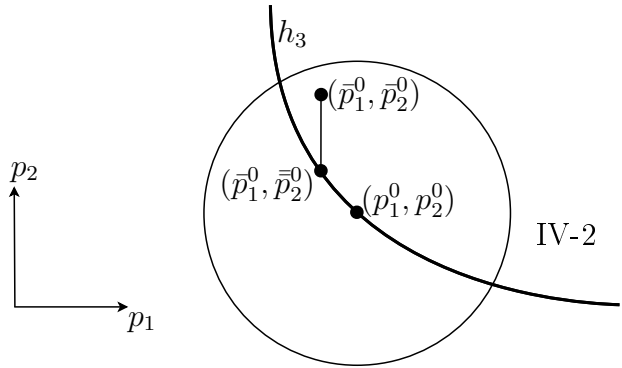}
\caption{A neighborhood of the point $(p_1^0,p_2^0)$ on $h_3$ for
$(\bar{p}_1^0,\bar{p}_2^0)$ in zone IV-2.} \label{fig:stab2}
\end{figure}

Fix a point $(\bar{p}_1^0,\bar{p}_2^0)$ in the intersection of $V$
and zone IV-2. In view of the first equation in \eqref{eq:DS2},
$p_1=\bar{p}_1^0$. Let $(\bar{p}_1^0,\bar{\bar{p}}_2^0)$ be a
point on the hyperbola $h_3$, i.e.
$\tilde{H}(\bar{p}_1^0,\bar{\bar{p}}_2^0)=0$ (see Figure
\ref{fig:stab2}). Then the second equation in \eqref{eq:DS2}
becomes
$$\dot{p}_2=\tilde{H}(\bar{p}_1^0,p_2)-\tilde{H}(\bar{p}_1^0,\bar{\bar{p_2}}).$$

We note that the derivative $$\frac{\partial}{\partial p_2
}\tilde{H}(p_1,p_2)=-\left[ \frac{\rho (p_1q_1 -Y_1)}{q_2
p^{\prime ~ 2}_2}+1 \right]$$ is strictly negative in a small
neighborhood of $(p_1^0,p_2^0)$. Indeed, as $(p_1^0,p_2^0)$ lies
on the graph of $h_3$, we have
$$p_1^0q_1+p_2^{0\prime}q_2-Y_1=\frac{p_2^{0\prime}}{p_2^0}Y_2>0$$ and,
consequently,
$$\frac{\rho}{p_2^{0\prime}}(Y_1-p_1^0q_1)<Y_1-p_1^0q_1<q_2p_2^{0\prime}, \textrm{ i.e. }\frac{\partial}{\partial p_2
}\tilde{H}(p_1^0,p_2^0)<0.$$ The claim follows from the latter
observation as $\frac{\partial}{\partial p_2 }\tilde{H}(p_1,p_2)$
is continuous. If we further contract the neighborhood $V$ so as
to ensure that $\frac{\partial}{\partial p_2
}\tilde{H}(p_1,p_2)<0$ in it, the equation under consideration
becomes $$\dot{p}_2(t)=(p_2(t)-\bar{\bar{p}}_2^0).\int_0^1
\frac{\partial}{\partial p_2
}\tilde{H}(\bar{p}_1^0,\bar{\bar{p}}_2^0+s(p_2(t)-\bar{\bar{p}}_2^0))ds.$$
Then
$$\frac{d}{dt}(p_2(t)-\bar{\bar{p}}_2^0)^2=2(p_2(t)-\bar{\bar{p}}_2^0)\dot{p}_2=2(p_2(t)-\bar{\bar{p}}_2^0)^2\int_0^1
\frac{\partial}{\partial p_2
}\tilde{H}(\bar{p}_1^0,\bar{\bar{p}}_2^0+s(p_2(t)-\bar{\bar{p}}_2^0))ds<0.$$
In words, for $(\bar{p}_1^0,\bar{p}_2^0)\in V$, the expression
$|(p_2(t)-\bar{\bar{p}}_2^0)|$ does not increase as $t \rightarrow
+\infty $ and stability is established since
$|(p_2^0-\bar{\bar{p}}_2^0)|$ is small.

One obtains analogous results for the points on the graph of
$h_4$. (For determinacy, we shall consider the setup in Figure
\ref{fig:pricespace1}, when $\tilde{E}$ is below the line
$p_2=p_1+\rho$, i.e. $\Delta p > -\rho$.)

Also, it is easy to verify that for initial data
$(\bar{p}_1^0,\bar{p}_2^0)$ in zones III or I-1, the phase
trajectories for $t\rightarrow +\infty$ tend to $\tilde{E} \left(
\frac{Y_1}{q_1},\frac{Y_2}{q_2} \right)$. In zone III the
differential equations system for the prices has the form
$$\left\{
  \begin{array}{l}
    \dot{p}_i=\frac{Y_i}{q_i}-p_i,~i=1,2  \\
    p_i(0)=\bar{p}_i^0
  \end{array}
 \right.$$ and its solution is $$p_i(t)=\left( \bar{p}_i^0 -\frac{Y_i}{q_i}
 \right)e^{-t}+\frac{Y_i}{q_i},$$which makes the claim obvious.

For initial data in zone I-1 (again in the setup from Figure
\ref{fig:pricespace1}, i.e. for $\Delta p \in [-\rho,\rho]$), the
differential system for the prices coincides with that for zone
III, which was just described.

To conclude, asymptotic stability does not in general hold even
for the point $\tilde{E}$, regardless of the properties of initial
data from the above described zones, for which the solutions of
the system tend to $\tilde{E}$. This conclusion remains valid for
other relations between the quantities $Y_i/p_i$ and the
transportation costs $\rho$ (see Figure \ref{fig:pricespace2}).
These observations explain the effects under stochastic
perturbations of the prices, obtained in section
\ref{sec:stochprice}.

\section{Price dynamics with stochastic shocks}\label{sec:stochprice}

The model studied here is deterministic and the agents are assumed
to have complete information. Given that this model abstracts from
many real-world complications, it would be worthwhile to study its
behaviour with respect to perturbations in some of the exogenous
variables. In this section we look at the case of adding shocks to
the prices by means of incorporating a nuisance stochastic process
in the differential system describing their evolution.

\begin{remark} Before proceeding to develop the setup for the main stochastic
simulation, we note that, heuristically, it seems plausible to
expect that the stability properties of the dynamical system from
section \ref{sec:contprice} will, in some sense, be preserved in
the presence of well-behaved stochastic disturbances. In other
words, if the shocks disturbing the system are sufficiently
``regular'', one may expect the deterministic component to
dominate in the stochastic dynamical system. This intuition can be
illustrated graphically with the aid of computer simulations
featuring a series of one-sided positive or negative stochastic
shocks on the prices. A representative outcome of the simulations
is shown in Figure \ref{fig:posneg}. As the figure shows, the
one-sided disturbances cause the equilibrium outcome to drift
along the locus of fixed points of the (deterministic)
differential system. Moreover, for appropriate one-sided
disturbances and initial conditions, the equilibrium will drift
toward the point $\tilde{E}$, which was shown in the previous
section to enjoy somewhat stronger stability properties than the
other fixed points of the system. \endprf

\begin{figure}[ht]
\begin{center}
\includegraphics[width=7 cm, height=8 cm]{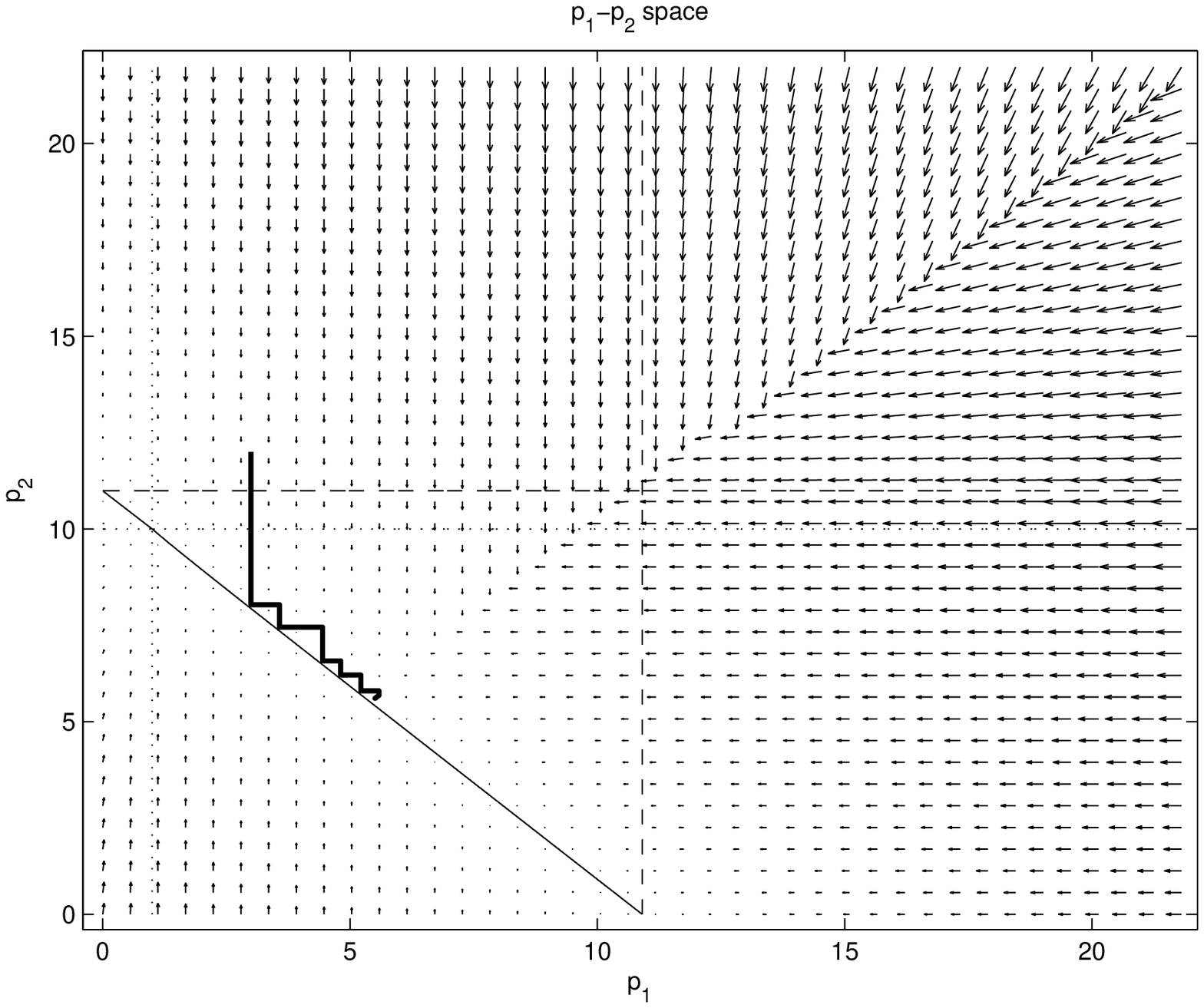}
\includegraphics[width=7 cm, height=8 cm]{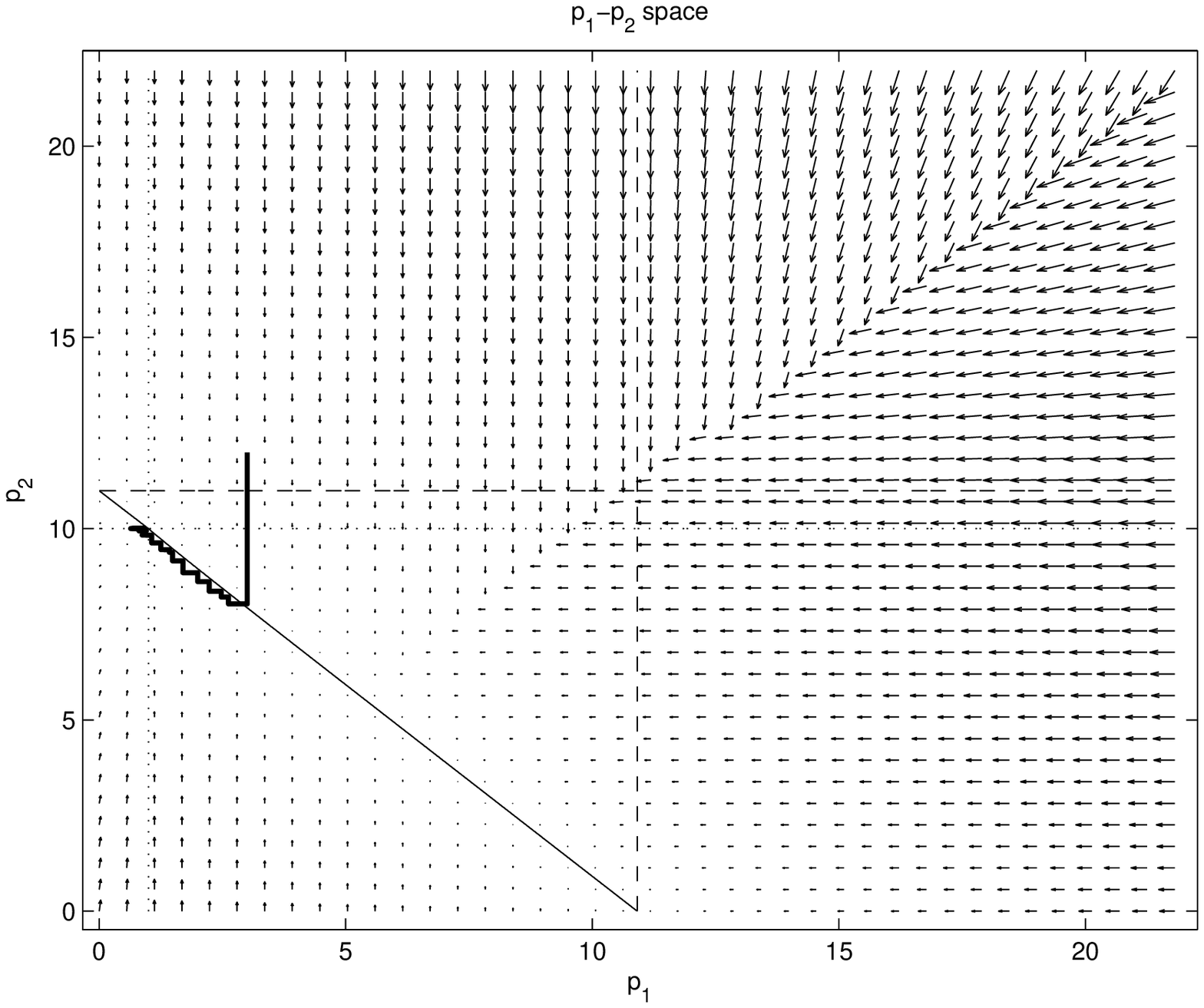}
\caption{Stability of the fixed points of the differential system
under one-sided positive (left-hand graph) and negative
(right-hand graph) stochastic shocks.} \label{fig:posneg}
\end{center}
\end{figure}\label{rem:onesidedshocks}\end{remark}

In the price equation \eqref{eq:2.5a} in discrete time we can
incorporate external random fluctuations by including a noise
variable $\Delta W_{t}$:
$$\frac{\Delta p_t}{p_t} = A(t)\Delta t + \sigma \Delta W_{t} \Delta t$$
Here $\sigma > 0$ is a coefficient characterizing the price
variability, $\Delta W_{t}$ are independent identically
distributed random variables that follow the standard Gaussian
distribution. If we rewrite the above equation as \beq
 \Delta p_t = A(t)p_t\Delta t + \sigma
p_t\Delta W_{t} \Delta t, \label{eq:Difference}\eeq under some
regularity conditions, at the limit $\Delta t \rightarrow 0$, the
solution of the difference equation \eqref{eq:Difference}
converges strongly to the solution of the stochastic differential
equation (SDE) \beq  d p(t) = A(t)p(t)dt + \sigma
p(t)dW_{t},\label{eq:Differential}\eeq where $dW_{t}$ is the It\^o
stochastic differential (for more details, see \cite[Theorem
9.6.2, p. 324]{Klo92}).

We choose a SDE of this type to govern the price dynamics in
continuous time. Let us consider a 2-dimensional Wiener process $W
= (W_{t}, t \in [0, T])$ with components $W_{t}^{1}$ and
$W_{t}^{2}$ which are independent scalar Wiener processes with
respect to a common family of $\sigma$-algebras
$\{\mathcal{A}_{t}, t \in [0, T]\}$. According to equation
\eqref{eq:2.7} we construct the following system
\begin{equation*}
\begin{split} d p_{t}^{1} &= Q_{1}(p_{t}^1,p_{t}^2)p_{t}^{1}dt +
\sigma_{1}p_{t}^{1}dW_{t}^{1}\\
d p_{t}^{2} &= Q_{2}(p_{t}^1,p_{t}^2)p_{t}^{2}dt +
\sigma_{2}p_{t}^{2}dW_{t}^{2}
\end{split}
\end{equation*}
In a more compact form
\begin{equation}\label{eq: SDE compact}
dp_{t} = a(p_{t})dt + b(p_{t})dW_{t},
\end{equation}
where $p_t = (p_{t}^1, p_{t}^2)$, $a(p_{t})$ is a 2-dimensional
vector function $a = (Q_{1}(.)p_{t}^{1}, Q_{2}(.)p_{t}^{2})':
\mathbb{R}^{2} \rightarrow \mathbb{R}^{2}$ and
$$b(p_{t}) = \begin{pmatrix}
  \sigma_{1}p_{t}^{1} & 0 \\
  0 & \sigma_{2}p_{t}^{2}
\end{pmatrix}$$
is a $2\times 2$ matrix function $b(t): \mathbb{R}^{2} \rightarrow
\mathbb{R}^{2\times 2}$. Some regularity conditions on $a(x)$,
$b(x)$, and the initial condition $p_{t_{0}}$ should be imposed
for the existence and uniqueness of a strong solution of the SDE,
meaning that the solution $p_{t}$ is a measurable functional of
$p_{t_{0}}$ and the Wiener process $W_{u},\ u\in [t_{0}, t]$ (see
Theorem 4.5.3 p. 131 and Theorem 4.5.6, p. 139 in \cite{Klo92}).
The classical conditions given, for example, in \cite{Klo92}
cannot be applied in our case because the function $a(x)$ violates
the Lipschitz condition. There is a result due to Zvonkin which
guarantees existence and uniqueness of a strong solution while
imposing weaker assumptions on $a(x)$, see Theorem 6.13, p. 152 in
\cite{Kleb98}. According to it, a strong solution of the
one-dimensional version of the SDE in (\ref{eq: SDE compact})
exists and is unique if $a(x)$ is a bounded function and $b(x)$ is
Lipschitz and bounded away from zero. While we are not aware of a
multi-dimensional extension of Zvonkin's theorem, we hypothesize
that a similar result holds. Under this hypothesis, a strong
solution of our SDE exists and is unique.

\begin{figure}[ht]
\begin{center}
\includegraphics[width=7 cm, height=8 cm]{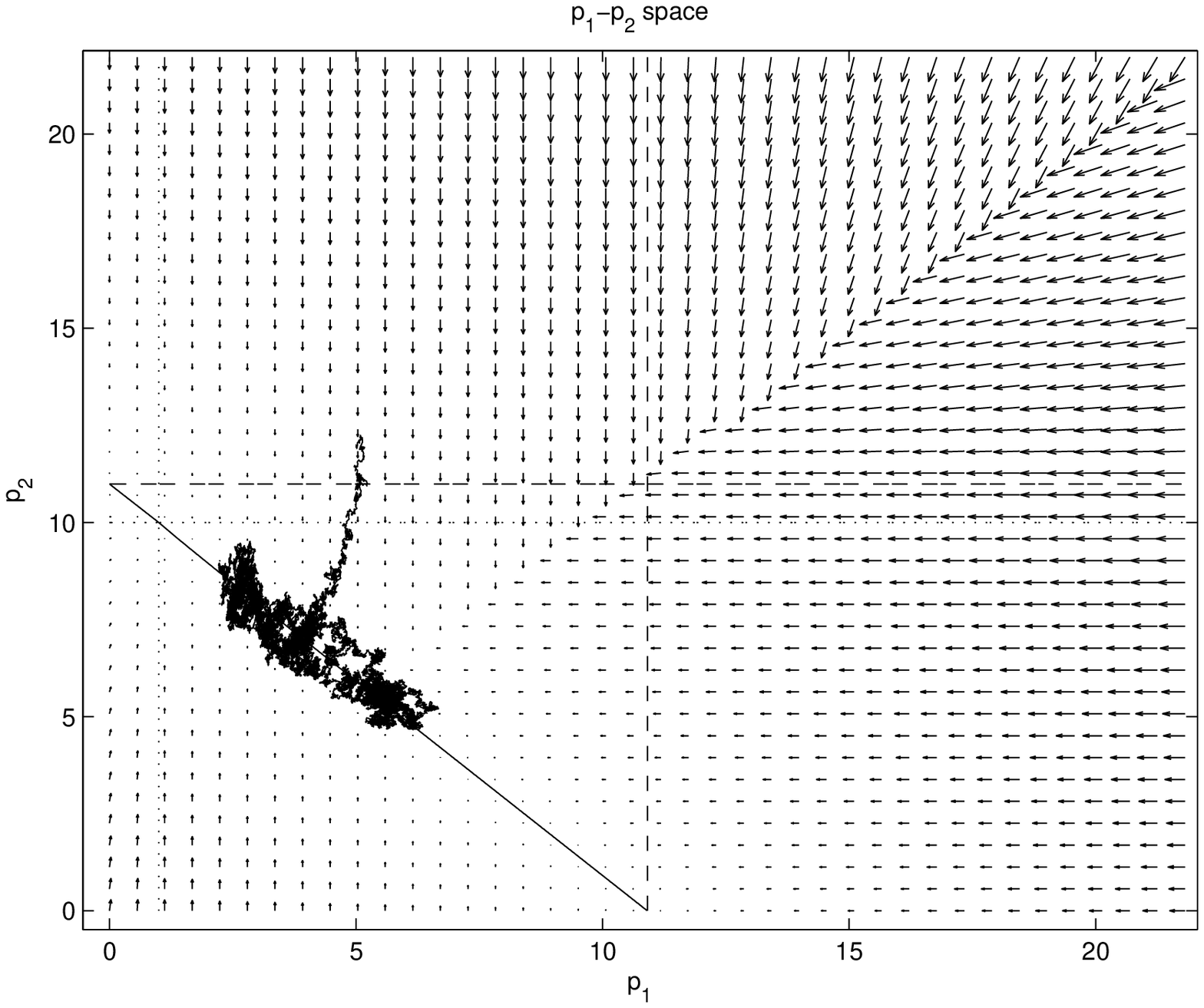}
\includegraphics[width=7 cm, height=8 cm]{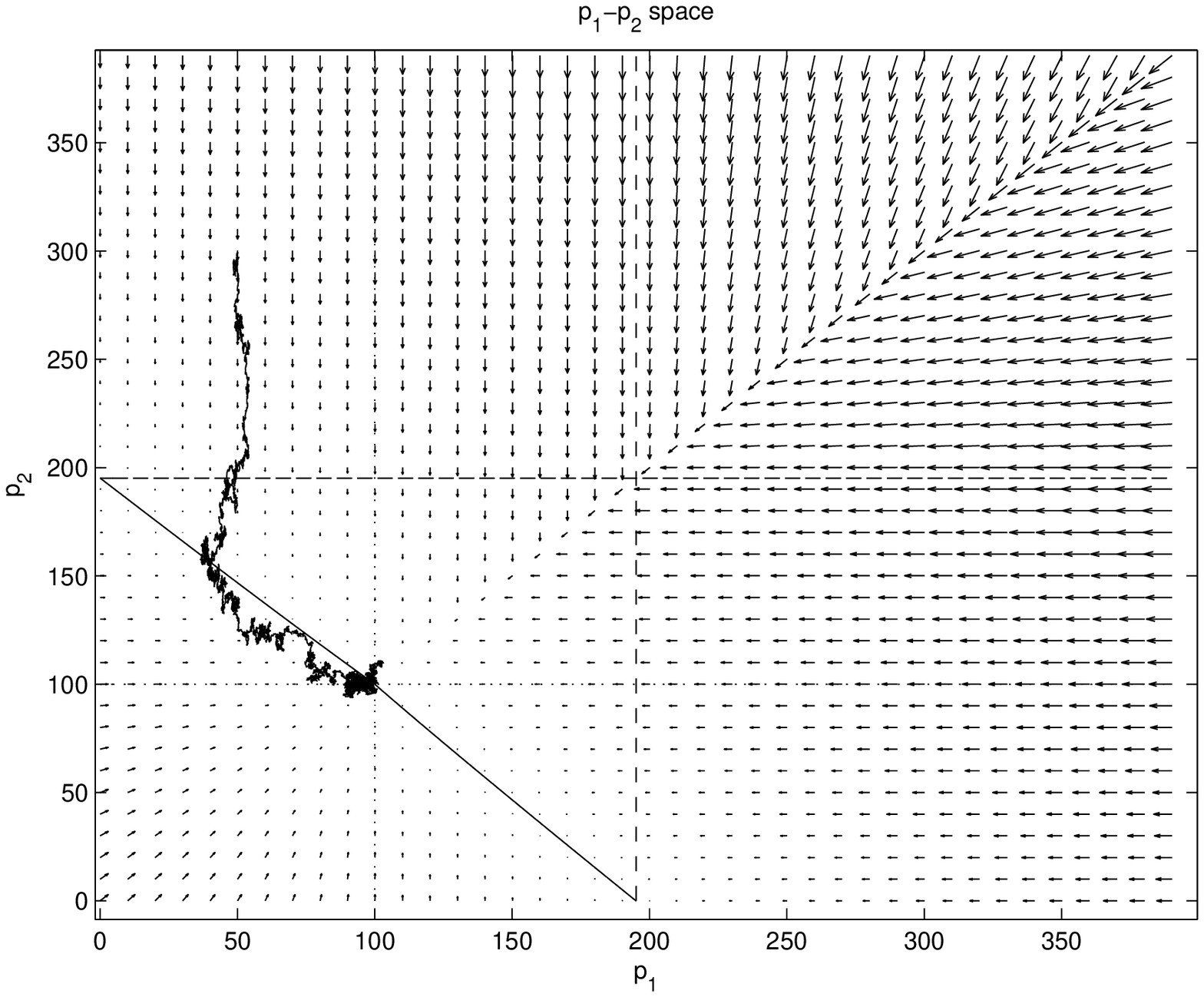}
\caption{Behaviour of the differential system under stochastic
shocks for different initial prices and parameterizations of the
problem (1)} \label{fig:stochsim1}
\end{center}
\end{figure}

\begin{figure}[ht]
\begin{center}
\includegraphics[width=7 cm, height=8 cm]{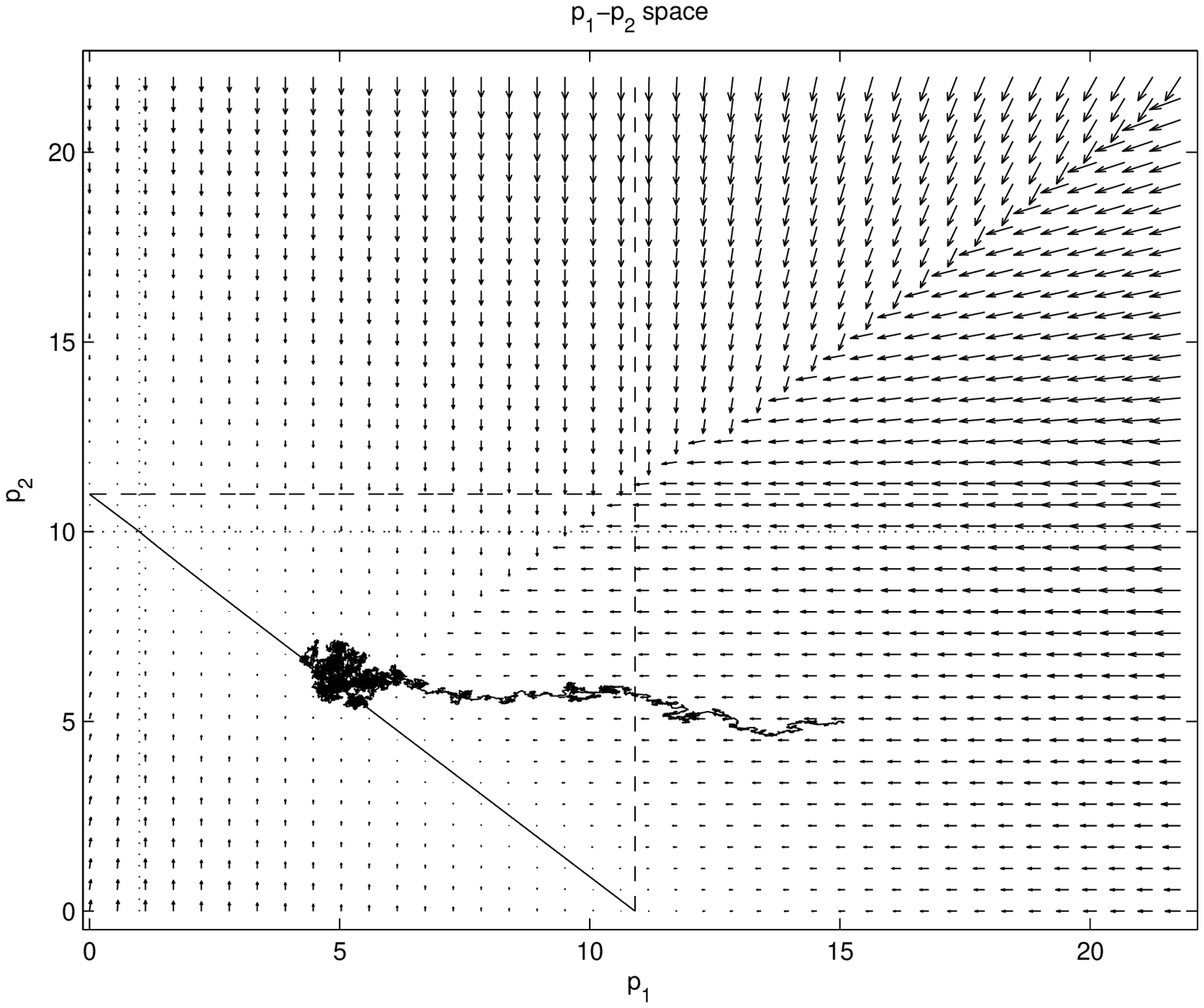}
\includegraphics[width=7 cm, height=8 cm]{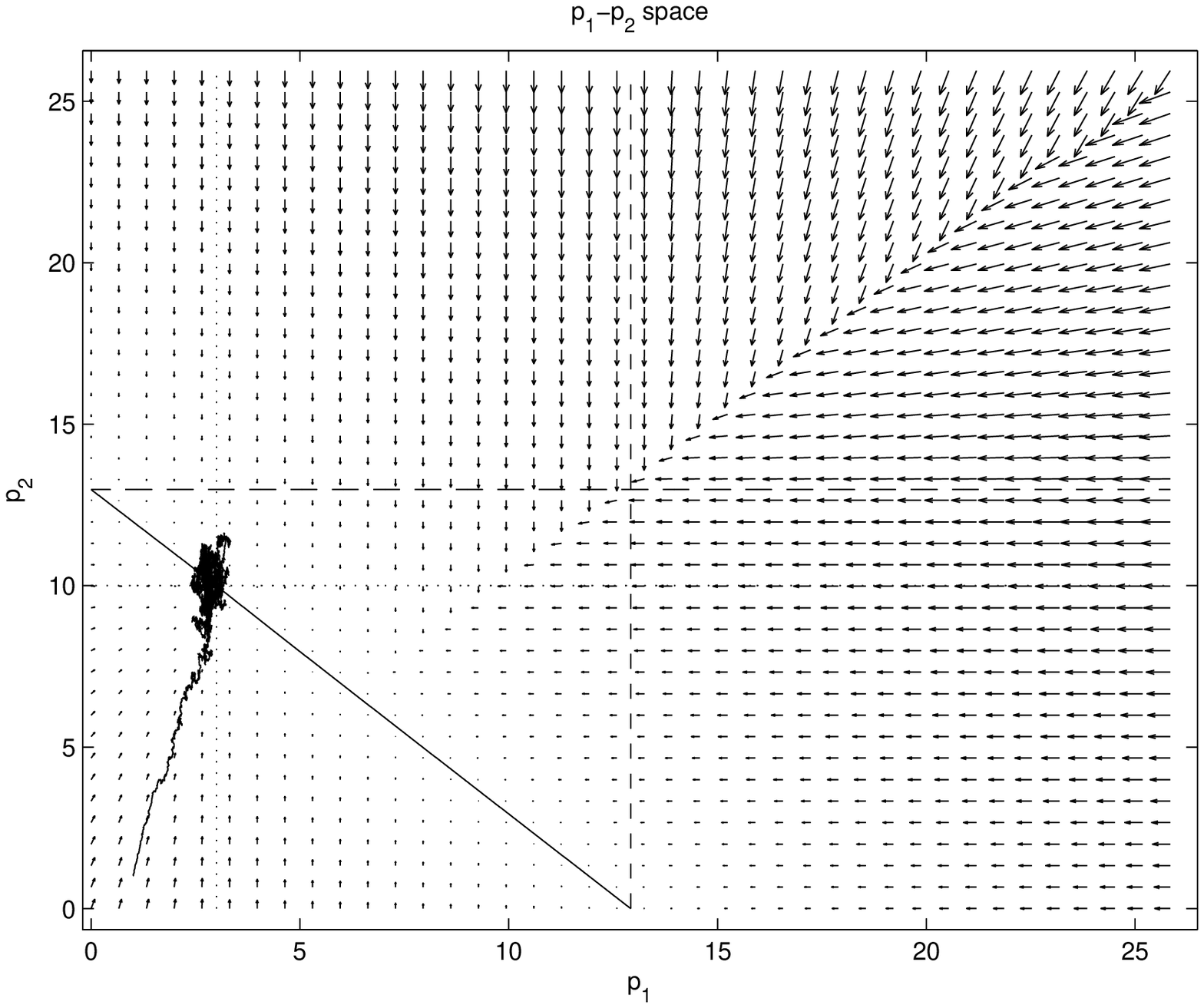}
\caption{Behaviour of the differential system under stochastic
shocks for different initial prices and parameterizations of the
problem (2)} \label{fig:stochsim2}
\end{center}
\end{figure}

We explore the sample paths of the solution in the phase space
employing the Euler scheme to solve the stochastic differential
system numerically. Figures \ref{fig:stochsim1} and
\ref{fig:stochsim2} illustrate the behaviour of the stochastic
differential system for different starting values of the prices.
In a fashion similar to the deterministic case, the sample path
approaches a stationary point depending on the initial condition.
If a stationary point has been reached, the random shocks perturb
the system away from it in a small neighborhood of the stationary
point. The simulation studies illustrate that provided the scales
$\sigma_{1}$ and $\sigma_{2}$ are small enough, the solution
remains in a small neighborhood of a stationary point.

\section*{Appendix}
\renewcommand{\theequation}{A.\arabic{equation}}
\setcounter{equation}{0}\numberwithin{equation}{subsection}
\renewcommand{\thesubsection}{A.\arabic{subsection}}
\setcounter{subsection}{0}
\renewcommand{\thefigure}{A.\arabic{figure}}
\setcounter{figure}{0}
\numberwithin{theorem}{subsection}
\setcounter{theorem}{0}

Here we sketch the proofs of Propositions
\ref{prop:zoneII}-\ref{prop:zoneIagain}. Although the proof of
Proposition \ref{prop:zoneIII} is contained in the main body of
the paper, we will provide a sketch for it as well in order to
illustrate the expository style adopted in this appendix.

We remind the reader that we always start with initial prices
$p_{i,0}$. These define, through the lines $\ell_i$, the zone in
the income space partition that the point $(Y_1,Y_2)$ belongs to.
We also assume that $(\alpha,\beta,\gamma,\delta)$ is a NE. After
determining the specific values of $(\alpha,\beta,\gamma,\delta)$,
we check to what extent financial resources have been used and
perform the necessary price adjustments. This leads to shifts in
the lines $\ell_i$, thus redefining the partition and changing the
position of $(Y_1,Y_2)$ with respect to the newly obtained zones.
Ultimately, we seek to find the respective p.e.s. Whenever the use
of more precise notation is called for, we write $\ell_{i,t}$,
$t=0,1,2,\ldots$ The current coordinates in the equations of the
respective lines are denoted by $(\tilde{Y}_1,\tilde{Y}_2)$.

\subsection{Sketch of the proof of Proposition \ref{prop:zoneIII} for zone
III}\label{sec:sketchZoneIII}

\textbf{1)} $q_2-\delta \leq 0
\xrightarrow{\textrm{1-I-}A}$$\alpha=q_1,\beta=0
\xrightarrow{\textrm{2-I-}A}\gamma=0,\delta=q_2$. We obtain
\textbf{NE} $\mathbf{(q_1,0,0,q_2)}$.

\textbf{2)} $0<q_2-\delta\leq Y_1/p'_2 \longrightarrow$ 1-II-$A_1$
or 1-II-$A_2$

- for 1-II-$A_1$: $\alpha=q_1 \xrightarrow{\textrm{2-I-}A}
\delta=q_2$ (impossible in case \textbf{2)}),

- for 1-II-$A_2$: $q_1-\alpha > 0$ - only for $\Delta p>\rho$
$\xrightarrow{\textrm{(3,1) 2-II-}A_2}$ for $\Delta p>-\rho$:
$\delta=q_2$ (impossible).

\textbf{3)} $Y_1/p'_2 < q_2-\delta \xrightarrow{\textrm{1-III-}A}$
$q_1-\alpha > 0$ only for $\Delta p>\rho$
$\xrightarrow{\textrm{(3,1)}}$ for $\Delta p>-\rho$:
2-II-$A_1~(A_2)$ and 2-III-$A$ $\longrightarrow$ $q_2-\delta = 0$,
which is impossible.

For the NE $(q_1,0,0,q_2)$ obtained, the supply of goods is
exhausted. If $Y_i=p_{i,0}q_i,~i=1,2$, the financial resources are
also exhausted, i.e. we are at a p.e. $p_{i,t}\equiv
p_{i,0},~\forall t \geq 0,~i=1,2$. If for some $i$, $i=1,2$, we
have $Y_i>p_{i,0}q_i$, the respective price $p_{i,0}$ is adjusted
to the level $p_{i,1}$, defined by the condition $Y_i=p_{i,1}q_i$.
Thus, we reach a p.e. $p_{i,t}\equiv p_{i,1},~\forall t \geq 1$.

\subsection{Sketch of the proof of Proposition \ref{prop:zoneII} for zone
II}\label{sec:sketchZoneII}

\subsubsection{Zone II-3 (see \eqref{eq:zoneII-3init})}\label{sec:sketchZoneII-3}

\textbf{1)} $q_2-\delta \leq 0 \xrightarrow{\textrm{1-I-}B}$
$\alpha = Y_1/p_1 (<q_1),\beta=0 \longrightarrow$ \textbf{1a)} or
\textbf{1b)}

\textbf{1a)} $0<q_1-\alpha<Y_2/p'_1
\xrightarrow{\textrm{2-II-}A_1}$ $\gamma=q_1-\alpha$,
$\delta=q_2$. We get \textbf{NE}
$\mathbf{(Y_1/p_1,0,q_1-Y_1/p_1,q_2)}$.

\textbf{1b)} $Y_2/p'_1 < q_1-\alpha$, which is impossible, since
for $\alpha=Y_1/p_1$ we obtain that $(Y_1,Y_2)$ is below $\ell_1$
and so below $\ell_4$.

\textbf{2)} $0<q_2-\delta \leq Y_1/p'_2
\xrightarrow{\textrm{1-II-}B}$ $\left\{
  \begin{array}{l}
   \Delta p\leq \rho:~ \alpha=Y_1/p_1(<q_1),~\beta=0  \\
   \Delta p> \rho:~\alpha=\frac{Y_1-p'_2(q_2-\delta)}{p_1}(< q_1),~\beta=q_2-\delta
  \end{array}
 \right\}$ $\longrightarrow$ \textbf{2a)} or
\textbf{2b)}

\textbf{2a)} $0<q_1-\alpha \leq Y_2/p'_1$

- for 2-II-$A_1$: $\delta=q_2$, which is impossible in case
\textbf{2)}.

- for 2-II-$A_2$: i) or ii)

i) for $\Delta p \geq -\rho$ (cases (1,1), (1,2), (2,1), (3,1))
$\longrightarrow$ $\delta = q_2$, which is impossible.

ii) for $\Delta p < -\rho$ (case (1,3)) $\longrightarrow$
$\delta=\frac{Y_2-p'_1 \left(q_1-\frac{Y_1}{p_1}\right)}{p_2}$
which implies, in view of the first inequality in \textbf{2)},
that $(Y_1,Y_2)$ is strictly below $\ell_4$.

\textbf{2b)} $Y_2/p'_1 < q_1-\alpha$
$\xrightarrow{\textrm{2-III-}A}$ i) or ii)

i) for $\Delta p \geq -\rho$: $\delta = q_2$ -- impossible in
\textbf{2)}

ii) for $\Delta p < -\rho$: case (1,3), which is impossible, since
inequality \textbf{2b)} for $\alpha=Y_1/p_1$ implies that
$(Y_1,Y_2)$ is below $\ell_1$.

\textbf{3)} $Y_1/p'_2 < q_2-\delta \xrightarrow{\textrm{1-III-}B}$
$\left\{
  \begin{array}{l}
   \Delta p\leq \rho:~ \alpha=Y_1/p_1,~\beta=0  \\
   \Delta p> \rho:~\alpha=0,~\beta=Y_1/p'_2
  \end{array}
 \right\}$ $\longrightarrow$ \textbf{3a)} or \textbf{3b)}

\textbf{3a)} $0<q_1-\alpha \leq Y_2/p'_1$ -- impossible, see
\textbf{2a)}

\textbf{3b)} $Y_2/p'_1 < q_1-\alpha$ -- impossible, see
\textbf{2b)}

For the unique NE, obtained in \textbf{1a)}, the quantities $q_1$,
$q_2$ and $Y_1$ are exhausted. The condition that $(Y_1,Y_2)$ is
above $\ell_4$, i.e. \beq
p'_{1,0}\left(q_1-\frac{Y_1}{p_{1,0}}\right)+p_{2,0}q_2 \leq Y_2
\label{eq:A_Yaboveell4}\eeq leads to two cases.

\textbf{Case I.} The condition \eqref{eq:A_Yaboveell4} holds with
equality. Then $p_{2,0}$ also remains unchanged, i.e. the points
on $\ell_4$ in zone II are p.e.s.

\textbf{Case II.} If there is a strict inequality in
\eqref{eq:A_Yaboveell4}, then
$$Y_2^{res}=Y_2-p'_{1,0}\left( q_1 - \frac{Y_1}{p_{1,0}} \right) > p_{2,0}q_2
.$$ We increase $p_{2,0}$ to $p_{2,1}$, for which $p_{2,1}q_2 =
Y_2^{res}$. With the new prices $p_{1,1}=p_{1,0}$ and $p_{2,1}$,
the point $(Y_1,Y_2)$ falls on the line $$\ell_{4,1} : p_{2,1}q_2
+ p'_{1,0}q_1 = \tilde{Y}_2+\frac{p'_{1,0}}{p_{1,0}}\tilde{Y}_1,$$
i.e. the p.e. is reached in one adjustment step.

\subsubsection{Zone II-2 (see \eqref{eq:zoneII-2init})}\label{sec:sketchZoneII-2}

\textbf{1)} $q_2-\delta \leq 0 \xrightarrow{\textrm{1-I-}B}$
$\alpha = Y_1/p_1 (<q_1),\beta=0 \longrightarrow$ \textbf{1a)} or
\textbf{1b)}

\textbf{1a)} $0<q_1-\alpha\leq Y_2/p'_1$

- case 2-II-$A_1$ is impossible, since $(Y_1,Y_2)$ is below
$\ell_4$

- case 2-II-$A_2$ $\longrightarrow$ \textbf{i)} or \textbf{ii)}

\textbf{i)} for $\Delta p \geq -\rho$: $\gamma=\frac{Y_2-p_2
q_2}{p'_1},~\delta=q_2$, which leads for $\bs{\Delta p \geq
-\rho}$ to $$\textrm{\textbf{NE}}\mathbf{\left(
\frac{Y_1}{p_1},0,\frac{Y_2-p_2q_2}{p'_1},q_2 \right).}$$

\textbf{ii)} for $\Delta p < -\rho$ $\longrightarrow$
$\delta<q_2$, which is impossible in \textbf{1)}.

\textbf{1b)} $q_1-\alpha > Y_2/p'_1$, which is impossible for
$\alpha=Y_1/p_1$.

\textbf{2)} $0<q_2-\delta \leq Y_1/p'_2
\xrightarrow{\textrm{1-II-}B}$ $\left\{
  \begin{array}{l}
  \Delta p \leq \rho:~ \alpha=Y_1/p_1 (<q_1),
\beta=0 \\
\Delta p > \rho: \alpha=\frac{Y_1-p'_2(q_2-\delta)}{p_1} (<q_1),
\beta=q_2-\delta
  \end{array}
 \right\}$ $\longrightarrow$ \textbf{2a)} or \textbf{2b)}

\textbf{2a)} $0<q_1-\alpha \leq Y_2/p'_1$

- case 2-II-$A_1$: $\delta=q_2$, which is impossible in
\textbf{2)}

- case 2-II-$A_2$ $\longrightarrow$ \textbf{i)} or \textbf{ii)}

\textbf{i)} for $\Delta p \geq -\rho$ $\longrightarrow$
$\delta=q_2$, which is impossible in \textbf{2)}

\textbf{ii)} for $\bs{\Delta p < -\rho}$: only (1,3)
$\longrightarrow$ We have
$$\textrm{\textbf{NE}}\mathbf{\left(\frac{Y_1}{p_1},0,q_1-\frac{Y_1}{p_1},\frac{Y_2-p'_1(q_1-Y_1/p_1)}{p_2}
\right)}.$$ (In this case the condition $\delta<q_2$ is obviously
satisfied. A comment on the second condition, $q_2-\delta \leq
Y_2/p'_1$, is offered following case \textbf{3a)}.)

\textbf{2b)} $Y_2/p'_1 < q_1-\alpha$ $\longrightarrow$ \textbf{i)}
or \textbf{ii)}

\textbf{i)} for $\alpha=Y_1/p_1$ (i.e. for $\Delta p \leq \rho$):
$(Y_1,Y_2)$ is strictly below $\ell_1$, which is impossible.

\textbf{ii)} for $\Delta p > \rho$: only (3,1)
$\xrightarrow{\textrm{2-III-}A}$ $\delta=q_2$, which is impossible
in \textbf{2)}.

\textbf{3)} $Y_1/p'_2 < q_2-\delta$
$\xrightarrow{\textrm{1-III-}B}$ $\left\{
  \begin{array}{l}
   \Delta p \leq \rho:~ \alpha=Y_1/p_1 (< q_1),~ \beta=0  \\
   \Delta p > \rho:~ \alpha=0 (<q_1),~ \beta=Y_1/p'_2
  \end{array}
 \right\}$ $\longrightarrow$ \textbf{3a)} or \textbf{3b)}

\textbf{3a)} $0<q_1-\alpha \leq Y_2/p'_1$ $\longrightarrow$
\textbf{i)} or \textbf{ii)}

\textbf{i)} for 2-II-$A_1$: $\delta=q_2$ -- impossible in
\textbf{3)}

\textbf{ii)} for 2-II-$A_2$ $\longrightarrow$ $\Delta p \geq
-\rho$ or $\Delta p < -\rho$

- for $\Delta p \geq -\rho$: $\delta = q_2$ -- impossible

- for $\Delta p < -\rho$: only (1,3) $\longrightarrow$
$\gamma=q_1-Y_1/p_1,~\delta = \frac{Y_2 -
p'_1(q_1-Y_1/p_1)}{p_2}$, which leads to the same NE as in case
2a)-(1,3). \newline (It turns out that whether $q_2-\delta$ is
greater than or less than $Y_1/p'_2$ is irrelevant, since we
obtain the same NE.)

\textbf{3b)} $Y_2/p'_1 < q_1-\alpha$
$\xrightarrow{\textrm{2-III-}A}$  case $\Delta p \geq -\rho$ or
case $\Delta p < -\rho$

- case $\Delta p \geq -\rho$: $\delta=q_2$ -- impossible

- case $\Delta p < -\rho$: only (1,3) $\longrightarrow$ condition
\textbf{3b)} would imply that $(Y_1,Y_2)$ is below $\ell_1$, which
is impossible.

We now turn to the study of the price dynamics, starting from the
NE obtained above.

\textbf{I)} \textbf{Analysis of the case NE}
$\mathbf{\left(\frac{Y_1}{p_{1,0}},0,q_1-\frac{Y_1}{p_{1,0}},\frac{Y_2-p'_{1,0}(q_1-Y_1/p_{1,0})}{p_{2,0}}
\right)}$ \textbf{(see 2a) or 3a) for $\bs{\Delta p < -\rho}$)}

The quantities $Y_1$, $q_1$ and $Y_2$ are depleted. Since
$(Y_1,Y_2)$ is strictly below $\ell_4$, we have $\delta< q_2$,
i.e. $q_2$ is not used up completely. Additionally, \beq
q_{2,0}^{cons}=\delta=\frac{1}{p_{2,0}}\left[
Y_2+\frac{p'_{1,0}}{p_{1,0}}Y_1 - p'_{1,0}q_1 \right]\geq 0,
\label{eq:A_q20consgeq0}\eeq since $(Y_1,Y_2)$ is on or above
$\ell_1$.

\textbf{Case I,i)}: In \eqref{eq:A_q20consgeq0} we have
$q_{2,0}^{cons} > 0$, i.e. $(Y_1,Y_2)$ is strictly above $\ell_1$.
Now $p_{2,0}$ decreases to $p_{2,1}$, which is defined by
$$p_{2,1} q_2 = p_{2,0}q_{2,0}^{cons}=Y_2-p'_{1,0}\left( q_1-\frac{Y_1}{p_{1,0}}
\right).$$ Since $p_{1,1}=p_{1,0}$, the point $(Y_1,Y_2)$ turns
out to be on the line $\ell_{4,1}$, which is parallel to
$\ell_{4,0}$, whose points are all equilibria, i.e. we reach a
p.e. in one adjustment step.

\begin{figure}[ht]
\centering
\input{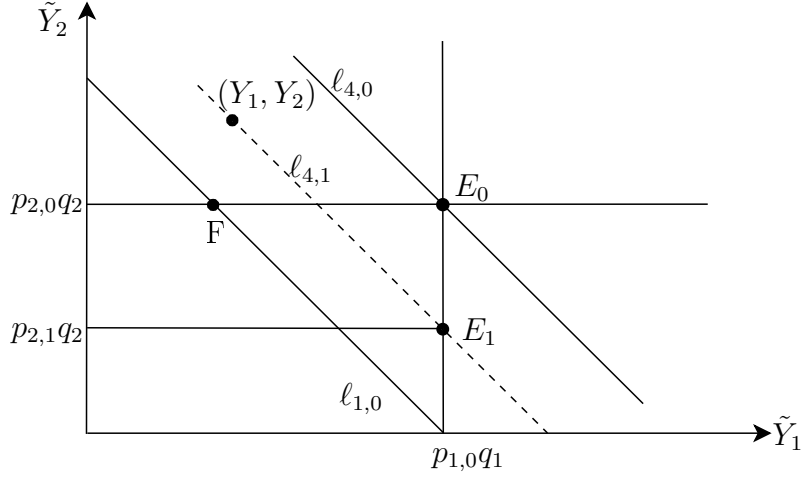}
\caption{Price adjustment in case \textbf{I,i)}.} \label{fig:A1}
\end{figure}

The above case is graphically illustrated in Figure \ref{fig:A1}.

\textbf{Case I,ii)}: In \eqref{eq:A_q20consgeq0} we have
$q_{2,0}^{cons}=0$, which is possible for the points in
$\ell_{1,0} \bigcap \{\textrm{Zone II}\}$, where
$q_1-Y_1/p_{1,0}=Y_2/p'_{1,0}$. Consequently, $p_{2,0}$ is reduced
to $p_{2,1}=0$ (while $p_{1,1}=p_{1,0}$) and we reach a degenerate
case. Writing for brevity $(p_1,0)$ instead of
$(p_{1,1},p_{2,1})$, this case is described by Table
\ref{tab:breply1} for $p'_2=\rho,p_1>0$, and Table
\ref{tab:breply2prime}. We now have the problem of finding NE
$(\alpha,\beta,\gamma,\delta)$ subject to the constraints \beq
0<Y_1<p_1q_1, 0<Y_2 \label{eq:Y1>p1q1Y2>0} \eeq and the additional
condition \beq \frac{Y_1}{p_1}+\frac{Y_2}{p'_1} = q_1.
\label{eq:A_AddCond} \eeq

We solve the problem through the already familiar approach:

\textbf{1)} $q_2-\delta \leq 0$ $\xrightarrow{\textrm{1-I-}B}$
$\alpha= Y_1/p_1 (<q_1),~\beta = 0$ $\longrightarrow$ \textbf{1a)}
or \textbf{1b)}

\textbf{1a)} $0<q_1-\alpha = Y_2/p'_1$
$\xrightarrow{\textrm{4-II}}$ $\gamma = q_1-\alpha,~\delta=q_2$.

This leads to \textbf{NE}
$\mathbf{\left(\frac{Y_1}{p_1},0,q_1-\frac{Y_1}{p_1},q_2 \right)}
= \mathbf{\left(\frac{Y_1}{p_{1,0}},0,\frac{Y_2}{p'_{1,0}},q_2
\right)}$

Since all resources are depleted, we reach a \textbf{degenerate
p.e.}, for which $$p_{1,t}=p_{1,0},~p_{2,t}=0, ~ \forall t \geq
1.$$ Obviously only the first economy is fully functioning, while
in the second economy local output becomes irrelevant as its
market price is zero.

\textbf{1b)} $Y_2/p'_1 < q_1-\alpha$, which is impossible since
$(Y_1,Y_2)$ is above $\ell_1$.

In the degenerate case there are no other NE, since for all
possible cases, namely

\textbf{2)} $0< q_2-\delta\leq Y_1/\rho$ \newline or

\textbf{3)} $Y_1/\rho < q_2-\delta$, \newline after applying Table
\ref{tab:breply2prime}, we obtain $\delta=q_2$, which leads to a
contradiction.

This completes the analysis of case \textbf{I)}.

\textbf{II)} \textbf{Analysis of the case NE}
$\mathbf{\left(\frac{Y_1}{p_{1,0}},0,\frac{Y_2-p_{2,0}q_2}{p'_{1,0}},q_2
\right)}$ \textbf{(see 1a) for $\bs{\Delta p \geq -\rho}$)}

The quantities $Y_1$, $Y_2$ and $q_2$ are depleted and the
condition that $(Y_1,Y_2)$ is below $\ell_4$ is equivalent to
$$q_1>\left[ \frac{Y_1}{p_{1,0}}+\frac{1}{p'_{1,0}}\left( Y_2 -p_{2,0}q_2 \right)
\right]=q_{1,0}^{cons},$$ i.e. $q_1$ is not depleted. At the same
time, $$q_{1,0}^{cons} =
\frac{Y_1}{p_{1,0}}+\frac{Y_2}{p'_{1,0}}-\frac{p_{2,0}q_2}{p'_{1,0}}
\geq q_1-\frac{p_{2,0}q_2}{p'_{1,0}} = \frac{\rho
q_1+p_{1,0}q_1-p_{2,0}q_2}{p'_{1,0}} > 0
$$ according to \eqref{eq:assump}. Consequently, $p_{1,0}$ is
reduced to $p_{1,1}$, where \beq
p_{1,1}q_1=p_{1,0}q_{1,0}^{cons}=Y_1+\frac{p_{1,0}}{p'_{1,0}}\left(
Y_2-p_{2,0}q_2 \right). \label{eq:A_p11determ}\eeq

The following subcases are possible:

\textbf{II,i)} $Y_2=p_{2,0}q_2$, i.e. $(Y_1,Y_2)$ is a point on
the segment $F E_0$ in Figure \ref{fig:A1}. Now the NE under
consideration takes the form $(Y_1/p_1,0,0,q_2)$. After the above
adjustment of $p_{1,0}$, $(Y_1,Y_2)$ turns out to be at $E_1$,
i.e. we reach a p.e. for which $$p_{1,t}=p_{1,1},~p_{2,t}=p_{2,0},
~ \forall t \geq 1.$$

\textbf{II,ii)} $Y_2>p_{2,0}q_2$. First we find the location of
$(Y_1,Y_2)$ with respect to the new position of $\ell_{4,0}$
(after the adjustment \eqref{eq:A_p11determ}) i.e. with respect to
$$\ell_{4,1}: p_{2,0}q_2+p'_{1,1}q_1 =
\frac{p'_{1,1}}{p_{1,1}}\tilde{Y}_1+\tilde{Y}_2.$$ We compare
$p_{2,0}q_2-Y_2$ and $$\frac{p'_{1,1}}{p_{1,1}}\left(
Y_1-p_{1,1}q_1 \right)=\frac{p'_{1,1}}{p_{1,1}}\left[
Y_1-Y_1-\frac{p_{1,0}}{p'_{1,0}}\left( Y_2-p_{2,0}q_2 \right)
\right] = \left( p_{2,0}q_2-Y_2 \right)\frac{p'_{1,1}}{p_{1,1}}
\frac{p_{1,0}}{p'_{1,0}}.$$ Since $p_{1,1}=p_{1,0}-\Delta$,
$\Delta>0$, it is easy to check that $$\frac{p'_{1,1}}{p_{1,1}}
\frac{p_{1,0}}{p'_{1,0}}>1,$$ which, after multiplication by
$p_{2,0}q_2-Y_2(<0)$, yields $$p_{2,0}q_2-Y_2 >
\frac{p'_{1,1}}{p_{1,1}} (Y_1-p_{1,1}q_1).$$ Consequently,
$(Y_1,Y_2)$ turns out to be \emph{below} $\ell_{4,1}$.

The results obtained are illustrated graphically in Figure
\ref{fig:A2}. Let $\tilde{\ell}$ denote a line through the point
$(Y_1,Y_2)$, which is parallel to $\ell_{4,0}$, i.e.
$$\tilde{\ell} : \frac{p'_{1,0}}{p_{1,0}}\tilde{Y}_1+\tilde{Y}_2= const \left(= \frac{p'_{1,0}Y_1}{p_{1,0}}+Y_2 \right).$$

Obviously, the point M$(p_{1,1}q_1,p_{2,0}q_2)$ lies on
$\tilde{\ell}$, as well as on $\ell_{4,1}$. Also, since $$\tan
\theta_0 = \frac{p'_{1,0}}{p_{1,0}} <
\frac{p'_{1,1}}{p_{1,1}}=\tan \theta_1 ,$$ it follows that
$$\theta_0<\theta_1 .$$

Thus, the line $\ell_{4,1}$ must turn in the negative direction
around the point M to coincide with $\tilde{\ell}$. Since
$(Y_1,Y_2)$ lies on $\tilde{\ell}$, it is located \emph{below}
$\ell_{4,1}$. At the same time, as $\ell_{1,1}$ is below
$\ell_{1,0}$, the point $(Y_1,Y_2)$ remains above $\ell_{1,1}$.

\begin{figure}[ht]
\centering
\input{FigA2.tex}
\caption{Price adjustment in case \textbf{II,ii)}} \label{fig:A2}
\end{figure}

Obviously, $p'_{1,1} < p'_{1,0}$ and $p_{2,0}q_2 < p'_{1,0}q_1$,
yet it is possible for $p'_{1,1}q_1$ to be greater or smaller than
$p_{2,0}q_2$ (see below).

We shall study separately the cases

\textbf{II,ii-1)} $\Delta p = -\rho$, i.e. $p'_{1,0}=p_{2,0}$,

\textbf{II,ii-2)} $\Delta p > -\rho$, i.e. $p'_{1,0}>p_{2,0}$.

In case \textbf{II,ii-1)}, after a downward adjustment of
$p_{1,0}$ to $p_{1,1}$ (see \eqref{eq:A_p11determ}) we obtain (in
the new zone II-2 -- see Figure \ref{fig:A2}, between $\ell_{4,1}$
and $\ell_{1,1}$) the case leading to NE for $\Delta p < -\rho$ of
the type in case \textbf{I,i)}. Consequently, after an adjustment
of $p_{2,1}=p_{2,0}$ to a smaller positive value $p_{2,2}$, the
point $(Y_1,Y_2)$ lies on $\ell_{4,2}$, whose points are
equilibria. Note also that even if $(Y_1,Y_2)\in \ell_{1,0}$, this
point will be \emph{strictly above} $\ell_{1,1}$ and so a
degenerate equilibrium cannot be obtained.

In case \textbf{II,ii-2)} there are many possibilities, which we
describe below. Suppose that, after the first adjustment of
$p_{1,0}$ as per \eqref{eq:A_p11determ} down to $p_{1,1}$, we
obtain the condition $\Delta p \leq -\rho$, i.e. $p'_{1,1} \leq
p_{2,1}$. (We have $p_{2,1}=p_{2,0}$, since only $p_{1,0}$ has
been changed.) Then, when consumption in the next period is
carried out ($t=2$), the above described adjustment according to
the NE of type \textbf{I,i)} obtains.

To find out whether such points exist at all, we write the
condition (which is the converse of the one mentioned above) \beq
p'_{1,1}>p_{2,0} \label{eq:A_p'11>p20} \eeq in the equivalent form
\beq Y_1+\frac{p_{1,0}}{p'_{1,0}}Y_2 >
(p_{2,0}-\rho)q_1+\frac{p_{1,0}}{p'_{1,0}}p_{2,0}q_2 .
\label{eq:A_equivp'11>p20}\eeq

Consequently, the condition \eqref{eq:A_p'11>p20} means that the
point $(Y_1,Y_2)$ is above the line $$\bar{\bar{\ell}} :
\tilde{Y}_1 + \frac{p_{1,0}}{p'_{1,0}} \tilde{Y}_2 =
(p_{2,0}-\rho)q_1+\frac{p_{1,0}}{p'_{1,0}}p_{2,0}q_2 ,$$ which is
parallel to $\ell_{4,0}$ and $\ell_{1,0}$. Since the abscissa of
the intersection point of $\bar{\bar{\ell}}$ with the
$\tilde{Y}_1$ axis is smaller than the abscissa of the
intersection point of $\ell_{4,0}$ (in view of \textbf{II,ii-2)}),
there exist points $(Y_1,Y_2)$ in zone II-2 with the property
\eqref{eq:A_p'11>p20}. Respectively, in the case when
\eqref{eq:A_p'11>p20} does not hold, the relevant points belong to
the closed area in zone II-2 enclosed between $\ell_{1,0}$ and
$\bar{\bar{\ell}}$ (when $\bar{\bar{\ell}}$ is between
$\ell_{1,0}$ and $\ell_{4,0}$), or the segment of
$\ell_{1,0}=\bar{\bar{\ell}}$ belonging to zone II-2 (when the
last two line coincide), and for then the adjustment process from
\textbf{II,ii-1)} obtains. When $\bar{\bar{\ell}}$ is strictly
below $\ell_{1,0}$, no such points exist.

For all points strictly above $\bar{\bar{\ell}}$ the conditions
$$p'_{1,0}>p_{2,0}\textrm{ and }p'_{1,1}>p_{2,0}$$ are
simultaneously valid.

\begin{figure}[ht]
\centering
\input{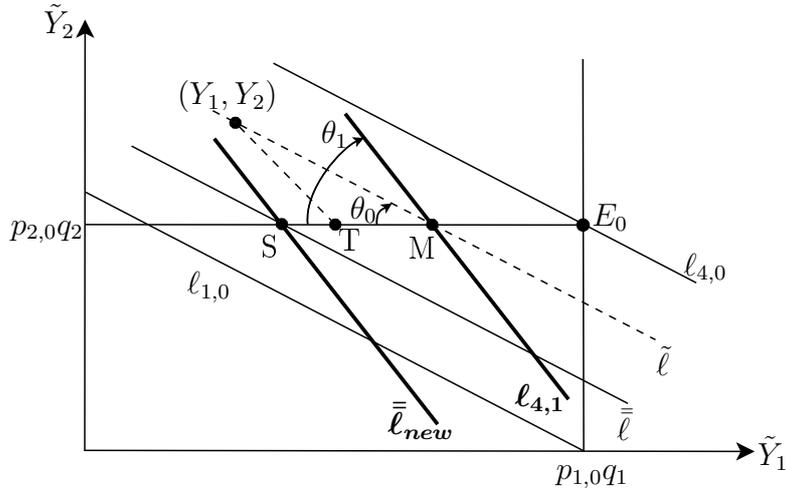}
\caption{Price adjustment in case \textbf{II,ii-2)}}
\label{fig:A3}
\end{figure}

Figure \ref{fig:A3} illustrates this case, with $\bar{\bar{\ell}}$
taken to lie between $\ell_{1,0}$ and $\ell_{4,0}$ for
determinacy.

In the new zone II-2, defined by $\ell_{1,1}$ and $\ell_{4,1}$, we
again obtain a NE of type II. However, because of the rotation at
an angle $\theta_1-\theta_0$ in the negative direction (see
above), it is not certain whether after the adjustment in
$p_{1,1}$ (in which a reduction to $p_{1,2}$ occurs), we can have
\beq p'_{1,2}>p_{2,2} ~(=p_{2,0}), \label{eq:A_p'12>p22} \eeq i.e.
whether $(Y_1,Y_2)$ would be above $\bar{\bar{\ell}}_{new}$ (see
Figure \ref{fig:A3}).

To describe all possible cases, we introduce the function \beq
g(x)=\frac{1}{q_1}\left[ Y_1 + \left( Y_2-p_{2,0}q_2
\right)\frac{x}{x+\rho} \right], \label{eq:A_gx} \eeq for which
$g(0)>0$, $g'(x)>0$, $\lim_{x\rightarrow \infty}g(x)>g(0)>0$ (see
Figure \ref{fig:A4}).

\begin{figure}[ht]
\centering
\input{FigA4.tex}
\caption{The function $g(x)$} \label{fig:A4}
\end{figure}

With the aid of the function $g(x)$, the condition that the point
$(Y_1,Y_2)$ is below the line $\ell_{4,0}$ becomes equivalent to
\beq p_{1,0}>g(p_{1,0}).\label{eq:A_p10>gp10}\eeq This, in
particular, implies that $p_{1,0}>k$, where $k$ is the only
positive number defined by $k=g(k)$.

Respectively, the condition \eqref{eq:A_p11determ} determining
$p_{1,1}$ can be written as \beq
p_{1,1}=g(p_{1,0}).\label{eq:A_p11equalsgp10} \eeq

If \eqref{eq:A_p'12>p22} were valid, in the next step we would set
$$p_{1,2}=g(p_{1,1})=g^2(p_{1,0})$$ and so on.

There are two possible cases:

\textbf{II,ii-2a)}: $p_{2,0}-\rho \leq k (<p_{1,0})$,

\textbf{II,ii-2b)}: $k<p_{2,0}-\rho ~(<p_{1,0})$.

\textbf{Case II,ii-2a).} In this case $$p_{1,1}=g(p_{1,0})>g(k)=k
\geq p_2-\rho .$$ We reach a NE of the type \textbf{II,ii)} and
set
$$p_{1,2}=g(p_{1,1})>k\geq p_2-\rho$$ etc. $$p_{1,t+1}=g(p_{1,t}),~\forall t > 0 .$$

Obviously the sequence $\{ p_{1,t} \}$ is convergent and tends to
$k$. The points S, T and M in Figure \ref{fig:A3} have abscissas
respectively $(p_{2,0}-\rho)q_1$, $kq_1$ and $p_{1,1}q_1$, and the
line through T and $(Y_1,Y_2)$ is the limit position of $\ell_4$
after infinitely many adjustments of the first price.

\textbf{Case II,ii-2b).} Let $$p_{2,0}-\rho \in
[g^s(p_{1,0}),g^{s-1}(p_{1,0)}],$$ where $s$ is a natural number.
In this case, after $s$ adjustments of the first price, we reach
an equilibrium for which the respective NE is of type \textbf{I)}
and the price adjustment process evolves accordingly.

\begin{remark} It is possible, as a result of the price reduction in the first
market, to reach for some $t\in\mathbb{N}$ the situation
$$p'_{1,t}\leq p_{2,0}q_2,$$ i.e. zone II-1 disappears. (As a matter of fact, this
is the case of zone  IV-2, with the roles of the two economies
reversed.) We can directly see that if $p'_{1,0}q_1 \leq
p_{2,0}q_2$, in the ``expanded'' zone II-2
\begin{equation*}\left\{
  \begin{array}{l}
    Y_1<p_1q_1,~Y_2\geq p_2q_2  \\
    (Y_1,Y_2)\textrm{ is stricly below }\ell_4
  \end{array}
\right.\end{equation*} one obtains the NEs of type \textbf{I)} and
\textbf{II)} derived above. The only qualitative difference here
is that no degenerate equilibria exist.
\endprf \label{rem:A_zoneII1disappears}\end{remark}

\subsubsection{Zone II-1 (see \eqref{eq:zoneII-1init})}\label{sec:sketchZoneII-1}

We note that this zone is characterized by relatively low
financial resources in both economies. We have $Y_1<Y_0$, where
the point $(Y_0,p_2q_2)=\ell_1 \bigcap \{ \tilde{Y}_2 = p_2 q_2
\}$ and $Y_2<p'_1 q_1$. To find the NE one proceeds as follows.

\textbf{1)} $q_2-\delta \leq 0 \xrightarrow{\textrm{1-I-}B}$
$\alpha = Y_1/p_1 (<q_1),\beta=0 \longrightarrow$ \textbf{1a)} or
\textbf{1b)}

\textbf{1a)} $0< q_1-\alpha \leq Y_2/p'_1$ is impossible, since
for $\alpha=Y_1/p_1$ the point $(Y_1,Y_2)$ would be above
$\ell_1$.

\textbf{1b)} $Y_2/p'_1 < q_1-\alpha$
$\xrightarrow{\textrm{2-III-}A}$ i) or ii)

i) for $\bs{\Delta p \geq -\rho}$: $\gamma =
\frac{Y_2-p_2q_2}{p'_1},\delta = q_2$, which leads to the
\textbf{NE} $\mathbf{\left( \frac{Y_1}{p_1},0,\frac{Y_2 -
p_1q_2}{p'_1},q_2 \right)}$.

ii) for $\Delta p < -\rho \longrightarrow $ $\delta=0
\longrightarrow q_2-\delta > 0$, which is impossible in
\textbf{1)}.

\textbf{2)} $0< q_2-\delta \leq Y_1/p'_2$
$\xrightarrow{\textrm{1-III-}B}$ $\left\{
  \begin{array}{l}
    \Delta p \leq \rho :~ \alpha=\frac{Y_1}{p_1}(<q_1),~\beta=0 \\
    \Delta p > \rho :~
    \alpha=\frac{Y_1-p'_2(q_2-\delta)}{p_1}(<q_1),~\beta=q_2-\delta
  \end{array}
 \right\}$ $ \longrightarrow$ \textbf{2a)} or \textbf{2b)}

\textbf{2a)} $0<q_1-\alpha \leq Y_2/p'_1$

- for $\Delta p \leq \rho$: impossible (see \textbf{1a)})

- for $\Delta p > \rho$: only (3,1) $\longrightarrow$ 2-II-$A_1$
and 2-II-$A_2$ ($\Delta p > -\rho$) $\longrightarrow$ $\delta =
q_2$ -- impossible.

\textbf{2b)} $q_1-\alpha > Y_2/p'_1$
$\xrightarrow{\textrm{2-III-}A}$ case $\Delta p \geq -\rho$ or
case $\Delta p < -\rho$

- for $\Delta p \geq -\rho$: $\delta = q_2$ -- impossible.

- for $\bs{\Delta p < -\rho}$: only in (1,3), $\gamma =
Y_2/p'_1,\delta = 0$, which leads to the \textbf{NE}
$\mathbf{\left( \frac{Y_1}{p_1},0,\frac{Y_2}{p'_1},0 \right)}$.
(See the comment after \textbf{3b)} for a check of the condition
$q_2-0 \leq Y_1/p'_2$.)

\textbf{3)} $Y_1/p'_2 < q_2-\delta $
$\xrightarrow{\textrm{1-III-}B}$ $\left\{
  \begin{array}{l}
    \Delta p \leq \rho :~ \alpha=Y_1/p_1(<q_1),~\beta=0 \\
    \Delta p > \rho :~
    \alpha=0(<q_1),~\beta=Y_2/p'_2
  \end{array}
 \right\}$ $ \longrightarrow$ \textbf{3a)} or \textbf{3b)}

\textbf{3a)} $0<q_1-\alpha \leq Y_2/p'_1$

- for $\Delta p \leq \rho$ -- impossible (see \textbf{1a)})

- for $\Delta p > \rho$ -- impossible (see \textbf{1a)})

\textbf{3b)} $q_1-\alpha > Y_2/p'_1$
$\xrightarrow{\textrm{2-III-}A}$ case $\Delta p \geq -\rho$ or
case $\Delta p < -\rho$

- for $\Delta p \geq -\rho$: $\delta=q_2$ -- impossible

- for $\Delta p < -\rho$: only in (1,3), $\gamma = Y_2/p'_1,\delta
= 0$, which leads to the same NE as in \textbf{2b)}. (It follows
that the check whether $Y_1/p'_2$ is less than or greater than
$q_2$ is unnecessary.)

We now turn to the study of the price dynamics, starting from the
NE obtained above.

\textbf{I)} \textbf{Analysis of the case NE}
$\mathbf{\left(\frac{Y_1}{p_{1,0}},0,\frac{Y_2}{p'_{1,0}},0
\right)}$ \textbf{for $\bs{\Delta p < -\rho}$ (see 2a) or 3b))}

The financial resources $Y_i$ are depleted, $q_1$ is only consumed
in part (since $(Y_1,Y_2)$ is strictly below $\ell_1$ and thus
$q_{1,0}^{cons}\in (0,q_1)$) and $q_2$ is unchanged
($q_{2,0}^{cons}=0$). Consequently, $p_{2,1}=0$ and
$p_{1,1}<p_{1,0}$ is determined by $$p_{1,1}q_1 =
p_{1,0}q_{1,0}^{cons} = Y_1 + \frac{p_{1,0}}{p'_{1,0}}Y_2 .$$ It
is immediately seen that the point $(Y_1,Y_2)$ remains below the
line
$$\ell_{1,1} : p_{1,1}q_1 =
\tilde{Y}_1+p_{1,1}\frac{\tilde{Y}_2}{p'_{1,1}}$$ and (omitting
the index $t=1$) this point lies in the following set (degenerate
zone II-1): \beq \left\{
  \begin{array}{l}
    \frac{Y_1}{p_1}+\frac{Y_2}{p'_1} < q_1 , \\
    0<Y_1<p_1q_1,~0<Y_2 .
  \end{array}
 \right. \label{eq:A_zoneII-1degen}\eeq

\begin{figure}[ht]
\centering
\input{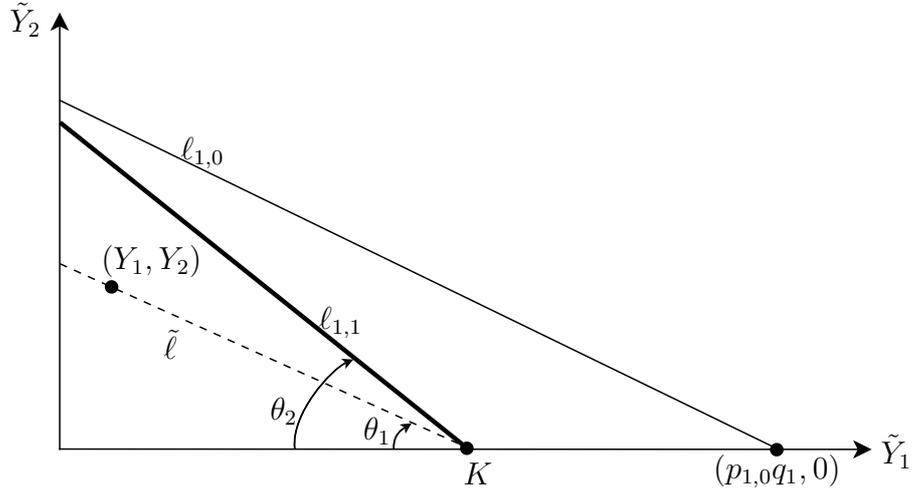}
\caption{The adjustment of the line $\ell_1$ when $p_1$ is
reduced} \label{fig:A4a}
\end{figure}

Figure \ref{fig:A4a} provides a geometric illustration of the
adjustment of the line $\ell_1$ in the case when the price $p_1$
is reduced.

The line $\tilde{\ell}$, passing through the point $(Y_1,Y_2)$, is
parallel to $\ell_{1,0}$ and has the equation $$\tilde{\ell} :
\tilde{Y}_1+\frac{p_{1,0}}{p'_{1,0}}\tilde{Y}_2=Y_1+\frac{p_{1,0}}{p'_{1,0}}Y_2
= p_{1,1}q_1 .$$ It intersects the line $$\ell_{1,1} :
\tilde{Y}_1+\frac{p_{1,1}}{p'_{1,1}}\tilde{Y}_2 = p_{1,1}q_1$$ at
the point $K (p_{1,1}q_1)$. Since (see Figure \ref{fig:A4a} for
the notation) $$\tan \theta_1 = \frac{p'_{1,0}}{p_{1,0}} <
\frac{p'_{1,1}}{p_{1,1}} = \tan \theta_2 \Rightarrow
\theta_1<\theta_2 ,$$ it is obvious that $(Y_1,Y_2)$ is below
$\ell_{1,1}$.

Using Table \ref{tab:breply1} (for $p'_2=\rho$ and $p_1>0$) and
Table \ref{tab:breply2prime}, we find the \textbf{NE}
$\mathbf{\left( \frac{Y_1}{p_{1,0}},0,\frac{Y_2}{p'_{1,0}},0
\right)}$ for the set defined in \eqref{eq:A_zoneII-1degen}. For
this NE, $Y_1$, $Y_2$ and $q_2$ are depleted, so the price
$p_1~(=p_{1,1})$ is reduced, as above, to
$p_{1,2}=h(p_{1,1})\equiv \frac{1}{q_1}\left[ Y_1 +
\frac{p_{1,1}}{p_{1,1}+\rho}Y_2 \right]$.

This adjustment process for the price is infinite and in the limit
we reach $$p_{1,\infty}=\lim_{t\rightarrow \infty}p_{1,t},$$ where
$p_{1,\infty}$ is the positive solution of the equation  $k=h(k)$.
In general, the prices tend to (although they never reach it) a
degenerate ``equilibrium'' $$p_{1,t}\rightarrow
p_{1,\infty},~p_{2,t}=0,~\forall t \geq 0 .$$ In this situation,
the limiting position of $(Y_1,Y_2)$ is on the line \beq
\ell_{1,\infty} :
\frac{\tilde{Y}_1}{p_{1,\infty}}+\frac{\tilde{Y}_2}{p'_{1,\infty}}=q_1
. \label{eq:A_ell1infty}\eeq From this one can easily obtain the
number \beq p_{1,\infty} = \frac{1}{2 q_1}\left( Y_1+Y_2-\rho q_1
+ \sqrt{(Y_1+Y_2-\rho q_1)^2+4\rho q_1 Y_1} \right).
\label{eq:A_numberp1infty}\eeq

Returning to the situation shown in Figure \ref{fig:A4a}, we note
that in the adjustment process for $p_1$ described above, the
points $(p_{1,t}q_1,0)$, which are counterparts to the point $K$,
tend to the limit point $(p_{1,\infty}q_1,0)$, while the lines
$\ell_{1,t}$ converge to the limit position $\ell_{1,\infty}$
(with the latter line passing through $(Y_1,Y_2)$).

\textbf{II)} \textbf{Analysis of the case NE}
$\mathbf{\left(\frac{Y_1}{p_{1,0}},0,\frac{Y_2-p_{2,0}q_2}{p'_{1,0}},q_2
\right)}$ \textbf{for $\bs{\Delta p \geq -\rho}$ (see 1b))}

The analysis and results in this case coincide with those for case
II) from b) from Proposition \ref{prop:zoneII} (when the
constraint coming from $\ell_1$ is not binding).

\subsection{Sketch of the proof of Proposition \ref{prop:zoneI} for zone
I}\label{sec:sketchZoneI}

The financial resources are smaller than the supply in both
economies, which technically means that we shall use parts $B$ in
Tables \ref{tab:breply1} and \ref{tab:breply2}. For the same
reason, the relationship between initial prices and transportation
costs plays an important role for the evolution of prices here.

\subsubsection{Zone I-1 (see \eqref{eq:zoneI-1init})}\label{sec:sketchZoneI-1}

\textbf{1)} $q_2-\delta \leq 0$ $\xrightarrow{\textrm{1-I-}B}$
$\alpha = Y_1/p_1 (<q_1),\beta=0$ $\longrightarrow$ \textbf{1a)}
or \textbf{1b)}

\textbf{1a)} $0<q_1-\alpha \leq Y_2/p'_1$
$\xrightarrow{\textrm{2-II-}B}$ $q_2-\delta \geq q_2-Y_2/p_2>0$,
which is impossible in \textbf{1)}.

\textbf{1b)} $Y_2/p'_1<q_1-\alpha $ -- impossible for $\alpha =
Y_1/p_1$, since it would imply that $(Y_1,Y_2)$ is strictly below
$\ell_1$.

\textbf{2)} $0<q_2-\delta \leq Y_2/p'_2$
$\xrightarrow{\textrm{1-II-}B}$ $\left\{
  \begin{array}{l}
   \Delta p \leq \rho:~ \alpha=Y_1/p_1 (<q_1),~\beta=0  \\
   \Delta p > \rho:~ \alpha=\frac{Y_1-p'_2(q_2-\delta)}{p_1}(<q_1),~\beta=q_2-\delta
  \end{array}
 \right\}$ $\longrightarrow$ \textbf{2a)} or \textbf{2b)}

\textbf{2a)} $0<q_1-\alpha \leq Y_2/p'_1$
$\xrightarrow{\textrm{2-II-}B}$ i) or ii)

i) for $\Delta p \geq -\rho$: $\gamma = 0,\delta = Y_2/p_2$

ii) for $\Delta p < -\rho$: $\gamma = q_1-\alpha,\delta =
\frac{Y_2-p'_1(q_1-\alpha)}{p_2}$

We obtain respectively:

- in cases (1,1), (1,2) and (2,1), i.e. for $\bs{\Delta p\in
[-\rho,\rho]}$: \textbf{NE}
$\mathbf{\left(\frac{Y_1}{p_1},0,0,\frac{Y_2}{p_2} \right)}$

- in case (1,3), $\bs{\Delta p < -\rho}$: \textbf{NE}
$\mathbf{\left(\frac{Y_1}{p_1},0,q_1-\frac{Y_1}{p_1},\frac{Y_2-p'_1\left(q_1-\frac{Y_1}{p_1}\right)}{p_2}
\right)}$

- in case (3,1), $\bs{\Delta p > \rho}$: \textbf{NE}
$\mathbf{\left(\frac{Y_1-p'_2\left(q_2-\frac{Y_2}{p_2}\right)}{p_1},q_2-\frac{Y_2}{p_2},0,
\frac{Y_2}{p_2} \right)}$.

(The condition from \textbf{2)} holds for $\delta=Y_2/p_2$, since
$(Y_1,Y_2)$ is above $\ell_2$. For a check of the condition when
$\delta$ is as in the NE for $\Delta p < -\rho$, see \textbf{3a)}.
The condition $q_1 \leq \alpha-Y_2/p'_1$ holds for $\alpha =
Y_1/p_1$, since $(Y_1,Y_2)$ is above $\ell_1$. For a check of the
condition when $\alpha$ is as in the NE for $\Delta p > \rho$, see
\textbf{2b)}.)

\textbf{2b)} $Y_2/p'_1<q_1-\alpha$ -- impossible for
$\alpha=Y_1/p_1$, since $(Y_1,Y_2)$ is above $\ell_1$
$\longrightarrow$ we have only case (3,1): $\Delta p > \rho$, for
which we find the same NE as in the respective case in
\textbf{2a)}. Therefore, it is unnecessary to compare $q_1$ and
$\alpha+Y_2/p'_1$.

\textbf{3)} $Y_1/p'_2 < q_2-\delta $
$\xrightarrow{\textrm{1-III-}B}$ $\left\{
  \begin{array}{l}
    \Delta p \leq \rho :~ \alpha=Y_1/p_1(<q_1),~\beta=0 \\
    \Delta p > \rho :~
    \alpha=0(<q_1),~\beta=Y_1/p'_2
  \end{array}
 \right\}$ $ \longrightarrow$ \textbf{3a)} or \textbf{3b)}

\textbf{3a)} $0<q_1-\alpha \leq Y_2/p'_1$

- for $\Delta p \leq \rho$ -- impossible, since $\delta=Y_2/p_2$
will violate the condition that $(Y_1,Y_2)$ is above $\ell_2$.

- we have only (1,3) for $\Delta p < - \rho$, in which case we
again arrive at the NE from \textbf{2a)}. (It follows that it is
unnecessary to compare $Y_1/p'_2$ and $q_2-\delta$.)

\textbf{3b)} $q_1-\alpha > Y_2/p'_1$

- for $\Delta p \leq \rho$ -- impossible, since $(Y_1,Y_2)$ is
above $\ell_1$

- for $\Delta p >\rho$: only in (3,1)
$\xrightarrow{\textrm{2-III-}B}$ $\delta = Y_2/p_2$, which is
impossible (see \textbf{3a)}).

\textbf{I)} \textbf{Analysis of the case NE}
$\mathbf{\left(\frac{Y_1}{p_{1,0}},0,0,\frac{Y_2}{p_{2,0}}
\right)}$ \textbf{for $\bs{\Delta p \in [-\rho,\rho ]}$ (see 2a))}

Since both quantities $q_i$ are not consumed completely
($q_{i,0}^{cons}=Y_i/p_i<q_i$), the prices are reduced once to
$$p_{i,1}=\frac{Y_i}{q_i}.$$ Then the point $(Y_1,Y_2)$ coincides
with $E_1$, the new position of the point $E_0$, which is a p.e.
(In a sense, the situation here is the exact opposite of that in
zone III, where the new point $E_1$ is reached after one upward
adjustment.)

\textbf{II)} \textbf{Analysis of the case NE}
$\mathbf{\left(\frac{Y_1}{p_{1,0}},0,q_1-\frac{Y_1}{p_{1,0}},\frac{Y_2-p'_{1,0}\left(q_1-\frac{Y_1}{p_{1,0}}\right)}{p_{2,0}}
\right)}$ \textbf{for $\bs{\Delta p < -\rho}$ (see 2a))}

Obviously $Y_1$, $Y_2$ and $q_1$ are depleted. Since $(Y_1,Y_2)$
is strictly below $\ell_{4,0}$,
$$q_{2,0}^{cons}=\frac{Y_2-p'_{1,0}\left(q_1-\frac{Y_1}{p_{1,0}}\right)}{p_{2,0}}<q_2
.$$ Consequently, $p_{2,0}$ is reduced to $p_{2,1}$: \beq
p_{2,1}q_2=Y_2-p'_{1,0}\left(q_1-\frac{Y_1}{p_{1,0}}\right).
\label{eq:A_p21determZoneI}\eeq

\textbf{II-1)}: $p_{2,1}>0$, i.e. $(Y_1,Y_2)$ is strictly above
$\ell_{1,0}$. Since $p_{1,0}$ is unchanged, i.e.
$p_{1,1}=p_{1,0}$, equation \eqref{eq:A_p21determZoneI} shows that
$(Y_1,Y_2)$ lies on $\ell_{4,1}$. Moreover,
$$Y_2-p_{2,1}q_2=Y_2-\left[Y_2-p'_{1,0}\left(q_1-\frac{Y_1}{p_{1,0}}\right)\right]>0,$$
i.e.  $(Y_1,Y_2)$ lies on the part of $\ell_{4,1}$ in the new zone
II-3, whose points are p.e.s

\textbf{II-2)} $p_{2,1}=0$, i.e. $(Y_1,Y_2)$ lies on $\ell_{1,0}$.
Thus, we reach a degenerate case ($p_1>0,p_2=0$), which was
analyzed for zone II-2, case \textbf{I-ii)}. It leads to a
degenerate $\ell_1$-equilibrium: $p_{2,t}=0,~\forall t \geq 1 .$

\textbf{III)} \textbf{Analysis of the case NE}
$\mathbf{\left(\frac{Y_1-p'_{2,0}\left(q_2-\frac{Y_2}{p_{2,0}}\right)}{p_{1,0}},q_2-\frac{Y_2}{p_{2,0}},0,\frac{Y_2}{p_{2,0}}
\right)}$ \textbf{for $\bs{\Delta p > \rho}$ (see 2a))}

The quantities $Y_1$, $Y_2$ and $q_2$ are depleted, and
$q_{1,0}^{cons}<q_1$, since $(Y_1,Y_2)$ is strictly below the line
$\ell_3$. After a reduction of $p_{1,0}$ to $p_{1,1}$, where
$$p_{1,1}q_1=Y_1-p'_{2,0}\left(q_2-\frac{Y_2}{p_{2,0}}\right),$$
the point $(Y_1,Y_2)$ lies on the new line $\ell_3$, whose points
in zone IV are p.e.s.

\subsubsection{Zone I-2 (see \eqref{eq:zoneI-2init})}\label{sec:sketchZoneI-2}

\textbf{1)} $q_2-\delta \leq 0$ $\xrightarrow{\textrm{1-I-}B}$
$\alpha = Y_1/p_1 (<q_1),\beta=0$ $\longrightarrow$ \textbf{1a)}
or \textbf{1b)}

\textbf{1a)} $0<q_1-\alpha \leq Y_2/p'_1$ -- impossible, since for
$\alpha=Y_1/p_1$ it contradicts the assumption that $(Y_1,Y_2)$ is
strictly below $\ell_1$.

\textbf{1b)} $Y_2/p'_1<q_1-\alpha $
$\xrightarrow{\textrm{2-III-}B}$ case $\Delta p \geq - \rho$ or
case $\Delta p < - \rho$

- case $\Delta p \geq - \rho$: $\delta=Y_2/p_2<q_2$, which
contradicts \textbf{1)}.

- case $\Delta p < - \rho$: $\delta=0$, which in \textbf{1)}
implies $q_2 \leq 0$ (impossible).

\textbf{2)} $0<q_2-\delta \leq Y_1/p'_2$
$\xrightarrow{\textrm{1-II-}B}$ for $\Delta p \leq \rho$:
$\alpha=Y_1/p_1 (<q_1),\beta=0$ or, for $\Delta p > \rho$:
$\alpha=\frac{Y_1-p'_2(q_2-\delta)}{p_1} ,\beta=q_2-\delta$
$\longrightarrow$ \textbf{2a)} or \textbf{2b)}

\textbf{2a)} $0<q_1-\alpha \leq Y_2/p'_1$, which is impossible,
since it would imply, for $\alpha \leq Y_1/p_1$, that $(Y_1,Y_2)$
is above $\ell_1$.

\textbf{2b)} $Y_2/p'_1 < q_1-\alpha$
$\xrightarrow{\textrm{2-III-}B}$ $\left\{
  \begin{array}{l}
    \Delta p \geq -\rho :~ \gamma = 0,~\delta=Y_2/p_2 \\
    \Delta p < - \rho :~ \gamma=Y_2/p'_1,~\delta=0
  \end{array}
 \right\}$ $ \longrightarrow$ three alternatives:

- for $\bs{\Delta p \in [-\rho,\rho]}$ (cases (1,1), (1,2) and
(2,1)): \textbf{NE}
$\mathbf{\left(\frac{Y_1}{p_1},0,0,\frac{Y_2}{p_2} \right)}$;

- for $\bs{\Delta p < -\rho}$ (case (1,3)): \textbf{NE}
$\mathbf{\left(\frac{Y_1}{p_1},0,\frac{Y_2}{p'_1},0 \right)}$;

- for $\bs{\Delta p > \rho}$ (case (3,1)): \textbf{NE}
$\mathbf{\left(\frac{Y_2-p'_2\left(q_2-\frac{Y_2}{p_2} \right
)}{p_1},q_2-\frac{Y_2}{p_2},0,\frac{Y_2}{p_2} \right)}$.

\textbf{3)} $Y_1/p'_2 < q_2-\delta$
$\xrightarrow{\textrm{1-III-}B}$ $\left\{
  \begin{array}{l}
    \Delta p \leq \rho :~ \alpha = Y_1/p_1 (<q_1),~\beta = 0 \\
    \Delta p > \rho :~ \alpha = 0,~\beta = Y_1/p'_2
  \end{array}
 \right\}$ $ \longrightarrow$ \textbf{3a)} or \textbf{3b)}

\textbf{3a)} $0<q_1-\alpha \leq Y_2/p'_1$, which is impossible
(see \textbf{2a)})

\textbf{3b)} $Y_2/p'_1<q_1-\alpha$
$\xrightarrow{\textrm{2-III-}B}$ i) or ii)

i) for $\Delta p \geq -\rho$: $\delta=Y_2/p_2$, which is
impossible, since \textbf{3)} implies that $(Y_1,Y_2)$ is strictly
below $\ell_2$.

ii) for $\Delta p < -\rho$, we only have case (1,3):
$\gamma=Y_2/p'_1,\delta=0$, which leads to the NE arising in
\textbf{2b)}. (Therefore, a comparison of $q_2$ and $Y_1/p'_2$ to
check the feasibility of the NE is unnecessary.)

\textbf{I)} \textbf{Analysis of the case NE}
$\mathbf{\left(\frac{Y_1}{p_{1,0}},0,0,\frac{Y_2}{p_{2,0}}
\right)}$ \textbf{for $\bs{\Delta p \in [ -\rho,\rho]}$ (see 2b))}

This case coincides with case \textbf{I)} in the analysis of zone
I-1.

\textbf{II)} \textbf{Analysis of the case NE}
$\mathbf{\left(\frac{Y_1}{p_{1,0}},0,\frac{Y_2}{p'_{1,0}},0
\right)}$ \textbf{for $\bs{\Delta p < -\rho}$ (see 2b))}

Here the  financial resources are completely spent;
$q_{1,0}^{cons}=Y_1/p_{1,0}+Y_2/p'_{1,0}<q_1$, since $(Y_1,Y_2)$
is strictly below $\ell_{1,0}$, $q_{2,0}^{cons}=0$. Then,
$p_{1,0}$ adjusts downward to $p_{1,1}$, where
$$p_{1,1}q_1=Y_1+\frac{p_{1,0}}{p'_{1,0}}Y_2,$$ while $p_{2,0}$
is adjusted downwards to $p_{2,1}=0$. It is immediately verified
that $(Y_1,Y_2)$ remains strictly below the line $$\ell_{1,1} :
p_{1,1}q_1=\tilde{Y}_1+\frac{p_{1,1}}{p'_{1,1}}\tilde{Y}_2 .$$
Thus, we obtain a degenerate zone II-1.

\textbf{III)} \textbf{Analysis of the case NE}
$\mathbf{\left(\frac{Y_1-p'_{2,0}\left(q_2-\frac{Y_2}{p_{2,0}}
\right)}{p_{1,0}}, q_2-\frac{Y_2}{p_{2,0}},0 ,\frac{Y_2}{p_{2,0}}
\right)}$ \textbf{for $\bs{\Delta p
>\rho}$ (see 2b))}

See case \textbf{III)} from zone I-1.

\subsubsection{Zone I-3 (see \eqref{eq:zoneI-3init})}\label{sec:sketchZoneI-3}

\textbf{1)} $q_2-\delta \leq 0$ $\xrightarrow{\textrm{1-I-}B}$
$\alpha = Y_1/p_1 (<q_1),\beta=0$ $\longrightarrow$ \textbf{1a)}
or \textbf{1b)}

\textbf{1a)} $0<q_1-\alpha \leq Y_2/p'_1$ -- impossible, since it
contradicts the assumption that $(Y_1,Y_2)$ is strictly below
$\ell_{2,0}$.

\textbf{1b)} $Y_2/p'_1<q_1-\alpha $
$\xrightarrow{\textrm{2-III-}B}$ case $\Delta p \geq - \rho$ or
case $\Delta p < - \rho$

- case $\Delta p \geq - \rho$: $\delta=Y_2/p_2<q_2$, which
contradicts \textbf{1)}.

- case $\Delta p < - \rho$: $\delta=0<q_2$, which is incompatible
with \textbf{1)}.

\textbf{2)} $0<q_2-\delta \leq Y_1/p'_2$
$\xrightarrow{\textrm{1-II-}B}$ $\left\{
  \begin{array}{l}
  \Delta p \leq \rho:~ \alpha=Y_1/p_1 (<q_1),~\beta=0  \\
  \Delta p > \rho:~ \alpha=\frac{Y_1-p'_2(q_2-\delta)}{p_1},~\beta=q_2-\delta
  \end{array}
 \right\}$ $\longrightarrow$ \textbf{2a)} or \textbf{2b)}

\textbf{2a)} $0<q_1-\alpha \leq Y_2/p'_1$, which is impossible,
since for $\alpha\leq Y_1/p_1$ it would imply that $(Y_1,Y_2)$ is
above $\ell_1$.

\textbf{2b)} $Y_2/p'_1 < q_1-\alpha$
$\xrightarrow{\textrm{2-III-}B}$ $\delta=Y_2/p_2$ or $0$, and by
the second inequality in \textbf{2)} it would imply that
$(Y_1,Y_2)$ is above $\ell_2$, which is impossible.

\textbf{3)} $Y_1/p'_2 < q_2-\delta$
$\xrightarrow{\textrm{1-III-}B}$ $\left\{
  \begin{array}{l}
    \Delta p \leq \rho :~ \alpha = Y_1/p_1 (<q_1),~\beta = 0 \\
    \Delta p > \rho :~ \alpha = 0 (<q_1),~\beta = Y_1/p'_2
  \end{array}
 \right\}$ $ \longrightarrow$ \textbf{3a)} or \textbf{3b)}

\textbf{3a)} $0<q_1-\alpha \leq Y_2/p'_1$, which is impossible
(see \textbf{2a)})

\textbf{3b)} $Y_2/p'_1<q_1-\alpha$
$\xrightarrow{\textrm{2-III-}B}$ i) or ii)

i) for $\Delta p \geq -\rho$: $\gamma=0,\delta=Y_2/p_2$;

ii) for $\Delta p < -\rho$: $\gamma=Y_2/p_2,\delta=0$.

From this we find:

- for $\bs{\Delta p \in [-\rho,\rho]}$ (cases (1,1), (1,2) and
(1,3)): \textbf{NE}
$\mathbf{\left(\frac{Y_1}{p_1},0,0,\frac{Y_2}{p_2} \right)}$;

- for $\bs{\Delta p < -\rho}$ (case (1,3)): \textbf{NE}
$\mathbf{\left(\frac{Y_1}{p_1},0,\frac{Y_2}{p'_1},0 \right)}$;

- for $\bs{\Delta p > \rho}$ (case (3,1)): \textbf{NE}
$\mathbf{\left(0,\frac{Y_1}{p'_2},0,\frac{Y_2}{p_2} \right)}$.

\textbf{I)} \textbf{Analysis of the case NE}
$\mathbf{\left(\frac{Y_1}{p_{1,0}},0,0,\frac{Y_2}{p_{2,0}}
\right)}$ \textbf{for $\bs{\Delta p \in [ -\rho,\rho]}$}

See case \textbf{I)} in the analysis of zone I-1.

\textbf{II)} \textbf{Analysis of the case NE}
$\mathbf{\left(\frac{Y_1}{p_{1,0}},0,\frac{Y_2}{p'_{1,0}},0
\right)}$ \textbf{for $\bs{\Delta p < -\rho}$ }

See case \textbf{II)} in the analysis of zone I-2.

\textbf{III)} \textbf{Analysis of the case NE}
$\mathbf{\left(0,\frac{Y_1}{p'_{2,0}},0,\frac{Y_2}{p_{2,0}}
\right)}$ \textbf{for $\bs{\Delta p > \rho}$ }

This case is symmetric (with respect to a change of roles of the
two economies) with case \textbf{II)}. The financial resources
$Y_1$, $Y_2$ are entirely spent, $q_{2,0}^{cons}=Y_1/p'_2+Y_2/p_2
< q_2$ (since $(Y_1,Y_2)$ is strictly below $\ell_2$), and $q_1$
is not consumed at all. Consequently, $p_{2,0}$ is reduced to
$p_{2,1}$: $$p_{2,1}q_2 = \frac{p_{2,0}}{p'_{2,0}}Y_1+Y_2 ,$$ and
$p_{1,1}=0$. We thus reach a degenerate case $p_1=0,~p_2>0$:
$$\left\{
  \begin{array}{l}
    0<Y_1,~0<Y_2<p_2q_2 , \\
    Y_2+\frac{p_2}{p'_2}Y_1<q_2 .
  \end{array}
 \right.$$

In the standard way, using Tables \ref{tab:breply1prime} and
\ref{tab:breply2} (for $p'_1=\rho,~p_2>0$), we find the next NE
$\left( q_1,\frac{Y_1}{p'_2},0,\frac{Y_2}{p_2} \right)$. After
that one obtains an infinite price adjustment process for the
price $p_2$, $\{p_{2,t}\}$, for which $$p_{2,t+1}q_2 =
\frac{p_{2,t}}{p'_{2,t}}Y_1+Y_2 ,$$ i.e. the system of two
economies tends to a degenerate $\ell_2$-equilibrium.

\subsection{Sketch of the proof of Proposition \ref{prop:zoneIV} for zone
IV}\label{sec:sketchZoneIV}

\subsubsection{Zone IV-1 (see \eqref{eq:zoneIV-1init})}\label{sec:sketchZoneIV-1}

\textbf{1)} $q_2-\delta \leq 0$ $\xrightarrow{\textrm{1-I-}A}$
$\alpha = q_1$ $\xrightarrow{\textrm{2-I-}B}$
$\delta=Y_2/p_2<q_2$, which is incompatible with \textbf{1)}.

\textbf{2)} $0<q_2-\delta \leq Y_1/p'_2$

- for 1-II-$A_1$: $\alpha=q_1,\beta=q_2-\delta$
$\xrightarrow{\textrm{2-I-}B}$ $\gamma=0,\delta=Y_2/p_2$, which
leads to \newline \textbf{NE}
$\mathbf{\left(q_1,q_2-\frac{Y_2}{p_2},0,\frac{Y_2}{p_2}\right)}$.
(The inequalities in \textbf{2)} hold. In particular, the second
one holds since $(Y_1,Y_2)$ is above $\ell_2$.)

- for 1-II-$A_2$: i) for $\Delta p \leq \rho$:
$\alpha=q_1,~\beta=\frac{Y_1-p_1q_1}{p'_2}$;

ii) for $\Delta p > \rho$:
$\alpha=\frac{Y_1-p'_2(q_2-\delta)}{p_1},~\beta=q_2-\delta$.

Then:

i) $\xrightarrow{\textrm{2-I-}B}$ $\gamma=0,\delta=Y_2/p_2$.
However, for this value of $\delta$ the condition from 1-II-$A_2$
together with the condition that $(Y_1,Y_2)$ is on or above
$\ell_3$ yield $q_1=\frac{Y_1-p'_2(q_2-Y_2/p_2)}{p_1}$, i.e.
$\frac{Y_1-p_1q_1}{p'_2}=q_2-\frac{Y_2}{p_2}$, which shows that in
this case we do not obtain a different NE from the one above.

ii) $\longrightarrow$ ii-1), ii-2) or ii-3)

ii-1) for $q_1\leq \alpha$ $\xrightarrow{\textrm{2-I-}B}$
$\gamma=0,\delta=Y_2/p_2$ and again the condition 1-II-$A_2$ and
the assumption that $(Y_1,Y_2)$ is on or above $\ell_3$ imply
$\alpha=q_1$, so that the familiar NE obtains.

ii-2) for $0<q_1-\alpha \leq Y_2/p'_2$
$\xrightarrow{\textrm{2-II-}B}$ (only (3,1)) for $\Delta p >
-\rho$: $\gamma=0,\delta=Y_2/p_2$, which is impossible in view of
the first inequality in this case ($\alpha<q_1$) and the condition
that $(Y_1,Y_2)$ is above $\ell_3$.

ii-3) for $Y_2/p'_2 < q_1-\alpha$ $\xrightarrow{\textrm{2-III-}B}$
(only (3,1)) $\delta=Y_2/p_2$, which is impossible in this case
(see ii-2)).

\textbf{3)} $Y_1/p'_2 < q_2-\delta$
$\xrightarrow{\textrm{1-III-}A}$ i) or ii)

i) for $\Delta p \leq \rho$:
$\alpha=q_1,\beta=\frac{Y_1-p_1q_1}{p'_2}$

ii) for $\Delta p > \rho$: $\alpha=0,\beta=\frac{Y_1}{p'_2}$

Respectively, we have:

i) $\xrightarrow{\textrm{2-I-}B}$ $\gamma=0,\delta=Y_2/p_2$, which
is a contradiction, since \textbf{3)} would imply that $(Y_1,Y_2)$
is below $\ell_2$ (impossible in zone IV-1).

ii) $\xrightarrow{\textrm{2-II-}B\textrm{ and }\textrm{2-III-}B}$
(only (3,1)) $\gamma=0,\delta=Y_2/p_2$ (impossible, as just shown
in i))

\textbf{Analysis of the case NE}
$\mathbf{\left(q_1,q_2-\frac{Y_2}{p_{2,0}},0,\frac{Y_2}{p_{2,0}}
\right)}$

The analysis and the results are symmetric (with respect to a
change of roles of the two economies) to those for zone II-3.

\subsubsection{Zone IV-2 (see \eqref{eq:zoneIV-2init})}\label{sec:sketchZoneIV-2}

\textbf{1)} $q_2-\delta \leq 0$ $\xrightarrow{\textrm{1-I-}A}$
$\alpha = q_1,\beta=0$ $\xrightarrow{\textrm{2-I-}B}$
$\gamma=0,\delta=Y_2/p_2$, which is incompatible with \textbf{1)}.

\textbf{2)} $0<q_2-\delta \leq Y_1/p'_2$ $\longrightarrow$ i) or
ii)

i) 1-II-$A_1$ $\longrightarrow$ $\alpha=q_1,\beta=q_2-\delta$
$\xrightarrow{\textrm{2-I-}B}$ $\delta=Y_2/p_2$, for which the
condition from 1-II-$A_1$ does not hold, since $(Y_1,Y_2)$ is
below $\ell_3$.

ii) 1-II-$A_2$ $\longrightarrow$ ii-1) or ii-2)

ii-1) for $\Delta p \leq \rho$:
$\alpha=q_1,\beta=\frac{Y_1-p_1q_1}{p'_2}$
$\xrightarrow{\textrm{2-I-}B}$ $\gamma=0,\delta=Y_2/p_2$ and for
$\bs{\Delta p \leq \rho}$ we get
$\textrm{\textbf{NE}}\mathbf{\left(q_1,\frac{Y_1-p_1q_1}{p'_2},0,\frac{Y_2}{p_2}
\right)}.$

ii-2) for $\Delta p > \rho$:
$\alpha=\frac{Y_1-p'_2(q_2-\delta)}{p_1},\beta=q_2-\delta$, so
that:

- If $q_1-\alpha \leq 0$ $\xrightarrow{\textrm{2-I-}B}$
$\delta=Y_2/p_2$, which is impossible, since $(Y_1,Y_2)$ is
strictly below $\ell_3$.

- If $0<q_1-\alpha \leq Y_2/p'_1$ $\xrightarrow{\textrm{2-II-}B}$
for (3,1): $\gamma=0,\delta=Y_2/p_2$ and we obtain
\newline$\textrm{\textbf{NE}}\mathbf{\left( \frac{Y_1-p'_2\left(
q_2-\frac{Y_2}{p_2}
\right)}{p_1},q_2-\frac{Y_2}{p_2},0,\frac{Y_2}{p_2} \right)}$. (As
the same NE arises under the assumption $Y_2/p'_1\leq q_1-\alpha$
(see 2-III-$B$ for $\Delta p > -\rho$), it is unnecessary to
compare $Y_2/p'_1$ and $q_1-\alpha$.)

\textbf{3)} $Y_1/p'_2 < q_2-\delta$
$\xrightarrow{\textrm{1-III-}A}$ i) or ii)

i) for $\Delta p \leq \rho$:
$\alpha=q_1,\beta=\frac{Y_1-p_1q_1}{p'_2}$ -- impossible, since by
2-I-$B$ we have $\gamma=0,\delta=Y_2/p_2$ and \textbf{3)} would
imply that $(Y_1,Y_2)$ is below $\ell_2$.

ii) for $\Delta p > \rho$: $\alpha=0,\beta=\frac{Y_1}{p'_2}$
$\longrightarrow$ $q_1-\alpha>0$:

- if $q_1-\alpha \leq Y_2/p'_1$ $\xrightarrow{\textrm{2-II-}B}$
(only (3,1)) $\gamma=0,\delta=Y_2/p_2$ and we obtain the same
contradiction from \textbf{3)}.

- if $Y_2/p'_1 < q_1-\alpha$ $\xrightarrow{\textrm{2-III-}B}$
(only (3,1)) $\gamma=0,\delta=Y_2/p_2$ and \textbf{3)} leads to a
contradiction.

\textbf{The analysis of the price adjustment in zone IV} is
analogous to the one in zone II-2, as the situations obtain are
symmetric as regards a change of roles of the two economies.

\subsection{Sketch of the proof of Proposition \ref{prop:zoneIagain} for case \eqref{eq:4.4iii} and \mbox{zone
I (1)} }\label{sec:sketchZoneIagain}

We first note that under the assumption made, in zone IV one
obtains the situation in zone II that was discussed under the
condition \eqref{eq:4.4i} (and, respectively, with interchanged
roles of the two economies).

\begin{figure}[ht]
\centering
\input{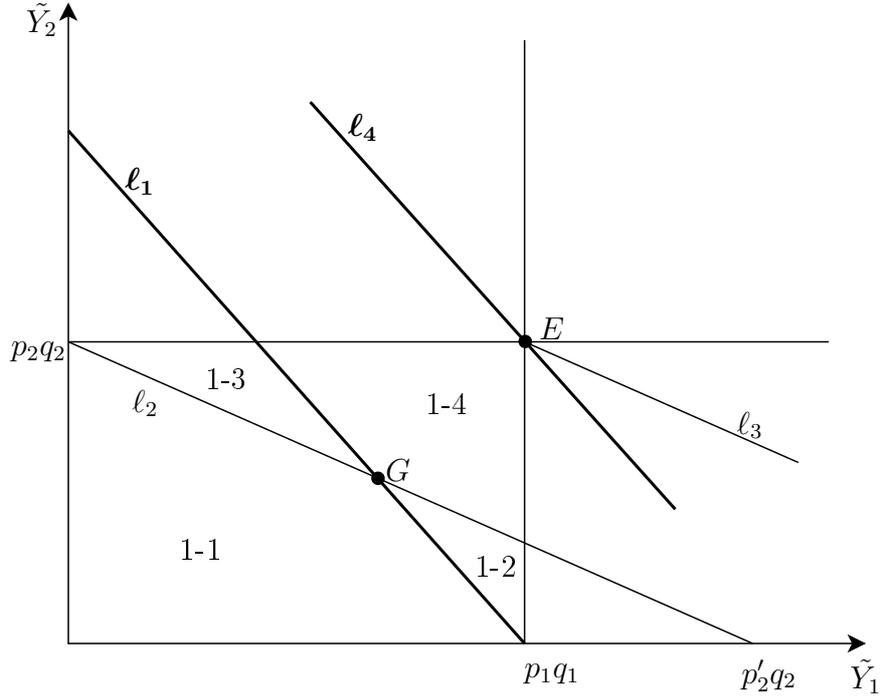}
\caption{The income space partition featuring zone 1}
\label{fig:A5}
\end{figure}

In this case zone 1 is divided into four subzones (see Figure
\ref{fig:A5}). (We draw the reader's attention to the fact that we
use Arabic numerals to denote the zone in the present setup.)

The case $p'_2q_2=p_1q_1$ is not qualitatively different from the
case \eqref{eq:4.4i}.

\subsubsection{Zone 1-4 (see \eqref{eq:zone1-4init})}\label{sec:sketchZone1-4}

\textbf{1)} $q_2-\delta \leq 0$ $\xrightarrow{\textrm{1-I-}B}$
$\alpha=Y_1/p_1 (<q_1),\beta=0$ $\longrightarrow$ \textbf{1a)} or
\textbf{1b)}

\textbf{1a)} $0<q_1-\alpha \leq Y_2/p'_1$
$\xrightarrow{\textrm{2-II-}B}$ i) or ii)

i) for $\Delta p \geq \rho $: $\delta = Y_2/p_2 <q_2$
(incompatible with \textbf{1)})

ii) for $\Delta p < \rho$:
$\delta=\frac{Y_2-p'_1(q_1-\alpha)}{p_2}$ -- impossible, since for
$\alpha=Y_1/p_1$ \textbf{1)} implies that $(Y_1,Y_2)$ is above
$\ell_4$.

\textbf{1b)} $Y_2/p'_1 < q_1-\alpha$ -- impossible, since for
$\alpha=Y_1/p_1$ one would get that $(Y_1,Y_2)$ is above $\ell_1$.

\textbf{2)} $0<q_2-\delta \leq Y_1/p'_2$
$\xrightarrow{\textrm{1-II-}B}$ i) or ii)

i) for $\Delta p \leq \rho $: $\alpha = Y_1/p_1 (<q_1),\beta=0$

ii) for $\Delta p > \rho $: $\alpha =
\frac{Y_1-p'_2(q_2-\delta)}{p_1} (<q_1),\beta=q_2-\delta$,

so we have \textbf{2a)} or \textbf{2b)}

\textbf{2a)} $0<q_1-\alpha \leq Y_2/p'_1$
$\xrightarrow{\textrm{2-II-}B}$ i) or ii)

i) for $\Delta p \geq -\rho $: $\gamma=0,\delta=Y_2/p_2$,

ii) for $\Delta p < -\rho $:
$\gamma=q_1-\alpha,\delta=\frac{Y_2-p'_1(q_1-\alpha)}{p_2}$,

From this we find:

- for $\bs{\Delta p \in [-\rho,\rho]}$ (cases (1,1), (1,2) and
(2,1)): \textbf{NE}
$\mathbf{\left(\frac{Y_1}{p_1},0,0,\frac{Y_2}{p_2} \right)}$;

- for $\bs{\Delta p > \rho}$ (case (3,1)): \textbf{NE}
$\mathbf{\left( \frac{Y_1-p'_2\left( q_2-\frac{Y_2}{p_2}
\right)}{p_1},q_2-\frac{Y_2}{p_2},0,\frac{Y_2}{p_2} \right)}$.

- for $\bs{\Delta p < -\rho}$ (case (1,3)): \textbf{NE}
$\mathbf{\left(\frac{Y_1}{p_1},0,q_1-\frac{Y_1}{p_1},\frac{Y_2-p'_1\left(
q_1-\frac{Y_1}{p_1} \right)}{p_2} \right)}$.

(For the above NEs conditions \textbf{2)} and \textbf{2a)}
obviously hold. We present more details on the feasibility of the
NEs after the analysis of cases \textbf{2b)} and \textbf{3)}.)

\textbf{2b)} $Y_2/p'_1 < q_1-\alpha$, which is impossible for
$\alpha=Y_1/p_1$, so that only case (3,1) for $\Delta p > \rho$ is
left. This case (by 2-III-$B$) leads to $\gamma=0,\delta=Y_2/p_2$
and we find the same NE as in \textbf{2a)}. (Thus, it is not
necessary to check inequalities \textbf{2b)} and the second
inequality in \textbf{2a)}.)

\textbf{3)} $Y_1/p'_2<q_2-\delta$ $\xrightarrow{\textrm{1-III-}B}$
i) or ii)

i) for $\Delta p \leq \rho $: $\alpha = Y_1/p_1 (<q_1),\beta=0$

ii) for $\Delta p > \rho $: $\alpha = 0 (<q_1),\beta=Y_1/p'_2$, so
that we have \textbf{3a)} or \textbf{3b)}

\textbf{3a)} $0<q_1-\alpha \leq Y_2/p'_1$, which for $\Delta p
\geq -\rho $ (by 2-II-$B$) leads to $\delta=Y_2/p_2$, so that
\textbf{3)} is impossible, while for $\Delta p < -\rho $ (case
(1,3)), we obtain the same result as in \textbf{2a)}.

\textbf{3b)} $Y_2/p'_1 <q_1-\alpha$, which is impossible for
$\alpha=Y_1/p_1$ (as $(Y_1,Y_2)$ is above $\ell_1$), so only (3,1)
is left and by 2-III-$B$ we find $\delta=Y_2/p_2$, for which
\textbf{3)} cannot hold, as $(Y_1,Y_2)$ is above $\ell_2$.

The price adjustment process for the NEs in question is the same
as in zone I-1 (basic case). The points on $(\ell_1 \bigcup
\ell_2)\bigcap \{\textrm{zone1-4}\}$ lead to a degenerate
equilibrium for which one of the prices becomes zero.

As a special illustration for the point $G$ we list the possible
cases:

a) for $\Delta p < -\rho $: the initial NE is $\left(
\frac{Y_1}{p_1},0,q_1- \frac{Y_1}{p_1},0 \right)$ and for prices
$p_1>0,p_2=0$ we reach the equilibrium $\left(
\frac{Y_1}{p_1},0,q_1- \frac{Y_1}{p_1},q_2 \right)$,

b) for $\Delta p > \rho $: the initial NE is $\left( 0,q_2-
\frac{Y_2}{p_2},0,\frac{Y_2}{p_2} \right)$ and for prices
$p_1=0,p_2>0$ we reach the equilibrium $\left( q_1,q_2-
\frac{Y_2}{p_2},0,\frac{Y_2}{p_2} \right)$,

c) for $\Delta p \in [-\rho,\rho] $: the initial NE is $\left(
\frac{Y_1}{p_1},0,0,\frac{Y_2}{p_2} \right)$ and, after a downward
adjustment of both prices, we reach a regular equilibrium that
coincides with the new position of the point $E_0$.

\subsubsection{Zone 1-3 (see \eqref{eq:zone1-3init})}\label{sec:sketchZone1-3}

\textbf{1)} $q_2-\delta \leq 0$ $\xrightarrow{\textrm{1-I-}B}$
$\alpha=Y_1/p_1,\beta=0$ $\longrightarrow$ \textbf{1a)} or
\textbf{1b)}

\textbf{1a)} $0<q_1-\alpha \leq Y_2/p'_1$, which is impossible,
since $(Y_1,Y_2)$ is below $\ell_1$.

\textbf{1b)} $Y_2/p'_1 < q_1-\alpha$
$\xrightarrow{\textrm{2-III-}B}$ $\delta \leq Y_2/p_2$, which
contradicts the assumption $Y_2<p_2q_2$.

\textbf{2)} $0<q_2-\delta \leq Y_1/p'_2$
$\xrightarrow{\textrm{1-II-}B}$ the same situation as in case
\textbf{2)} for zone 1-4.

Here, however, the case

\textbf{2a)} $0<q_1-\alpha \leq Y_2/p'_1$ is impossible, as
$\alpha \leq Y_1/p_1$ $\longrightarrow$ $(Y_1,Y_2)$ is above
$\ell_1$,

and in the case

\textbf{2b)} $Y_2/p'_1 <q_1-\alpha $ by 2-III-$B$ we find

- for $\Delta p \in [-\rho,\rho]$: NE
$\left(\frac{Y_1}{p_1},0,0,\frac{Y_2}{p_2} \right)$;

- for $\Delta p > \rho$: NE $\left( \frac{Y_1-p'_1\left(
q_2-\frac{Y_2}{p_2}\right)}{p_2},q_2-\frac{Y_2}{p_2},0,\frac{Y_2}{p_2}
\right)$, just as in the respective subcases from \textbf{2a)} in
zone 1-4;

- for $\Delta p < -\rho$: NE
$\left(\frac{Y_1}{p_1},0,\frac{Y_2}{p'_1},0 \right)$, which is
different from the equilibrium computed in \textbf{2a)} for zone
1-4.

\textbf{3)} $Y_1/p'_2 < q_2-\delta$
$\xrightarrow{\textrm{1-III-}B}$ $\left\{
  \begin{array}{l}
   \Delta p \leq \rho :~ \alpha = Y_1/p_1 ,~\beta = 0 \\
   \Delta p > \rho :~
    \alpha = 0 ,~\beta = Y_1/p'_2
  \end{array}
 \right\}$ $ \longrightarrow$ \textbf{3a)} or \textbf{3b)}

\textbf{3a)} $0<q_1-\alpha \leq Y_2/p'_1$, which is impossible,
since it would imply either that $(Y_1,Y_2)$ is above $\ell_1$, or
that $(Y_1,Y_2)$ is below $\ell_2$, both of which are wrong here.

\textbf{3b)} $Y_2/p'_1<q_1-\alpha$, for which, after eliminating
the impossible cases, we reach the NE from \textbf{2b)} for
$\Delta p < -\rho$.

The price adjustment process for $\Delta p \in [-\rho,+\infty )$
is the same as in zone 1-4 (i.e. as in the basic case for zone
I-1).

The price adjustment process for $\Delta p < -\rho$ with initial
NE $\left(\frac{Y_1}{p_1},0,\frac{Y_2}{p'_1},0 \right)$ is the
same as in the counterpart case for zone I-2 (basic case).

\subsubsection{Zone 1-2 (see \eqref{eq:zone1-2init})}\label{sec:sketchZone1-2}

\textbf{1)} $q_2-\delta \leq 0$ $\xrightarrow{\textrm{1-I-}B}$
$\alpha=Y_1/p_1(<q_1),\beta=0$. In this case neither

\textbf{1a)} $0<q_1-\alpha \leq Y_2/p'_2$, nor

\textbf{1b)} $Y_2/p'_2 < q_1-\alpha$ are possible, since they
contradict the inequality $Y_2<p_2q_2$ (by \textbf{1)}) or the
condition that $(Y_1,Y_2)$ is above $\ell_1$.

\textbf{2)} $0< q_2-\delta < Y_1/p'_2 $
$\xrightarrow{\textrm{1-II-}B}$ $\left\{
  \begin{array}{l}
    \Delta p \leq \rho :~ \alpha = Y_1/p_1 (<q_1),~\beta = 0 \\
    \Delta p > \rho :~ \alpha = \frac{Y_1-p'_2(q_2-\delta)}{p_1},~\beta = q_2-\delta
  \end{array}
 \right\}$ $ \longrightarrow$ \textbf{2a)} or \textbf{2b)}

\textbf{2a)} $0<q_1-\alpha \leq Y_2/p'_1$

- for $\Delta p \geq -\rho$ we obtain (by 2-II-$B$)
$\delta=Y_2/p_2$ which, together with \textbf{2)}, implies that
$(Y_1,Y_2)$ is above $\ell_2$, which is impossible. Then we have
only case (1,3), where for $\bs{\Delta p < -\rho}$ we find
\textbf{NE}
$\mathbf{\left(\frac{Y_1}{p_1},0,q_1-\frac{Y_1}{p_1},\frac{Y_2-p'_1\left(
q_1-\frac{Y_1}{p_1} \right)}{p_2} \right)}$.

\textbf{2b)} $Y_2/p'_1 < q_1-\alpha$ -- impossible, since for
$\alpha = Y_1/p_1$ one obtains that $(Y_1,Y_2)$ is above $\ell_1$,
and for (3,1) (by 2-III-$B$ and the inequality \textbf{2)}) one
finds that $(Y_1,Y_2)$ is above $\ell_2$.

\textbf{3)} $Y_1/p'_2 < q_2-\delta$
$\xrightarrow{\textrm{1-III-}B}$ $\left\{
  \begin{array}{l}
    \Delta p \leq \rho :~ \alpha = Y_1/p_1 ,~\beta = 0 \\
    \Delta p > \rho :~ \alpha = 0 ,~\beta = Y_1/p'_2
  \end{array}
 \right\}$ $ \longrightarrow$ \textbf{3a)} or \textbf{3b)}

\textbf{3a)} $0<q_1-\alpha \leq Y_2/p'_1$
$\xrightarrow{\textrm{2-II-}B}$ $\left\{
  \begin{array}{l}
    \Delta p \geq -\rho :~ \gamma = 0 ,~\delta = Y_2/p_2 , \\
    \Delta p < -\rho :~ \gamma = q_1-\alpha ,~\delta = \frac{Y_2-p'_1(q_1-\alpha)}{p_2}
  \end{array}
 \right.$

We obtain:

- for $\bs{\Delta p \in [-\rho,\rho]}$ (cases (1,1), (1,2) and
(2,1)): \textbf{NE}
$\mathbf{\left(\frac{Y_1}{p_1},0,0,\frac{Y_2}{p_2} \right)}$;

- for $\bs{\Delta p < -\rho}$ (case (1,3)): \textbf{NE}
$\mathbf{\left(\frac{Y_1}{p_1},0,q_1-\frac{Y_1}{p_1},\frac{Y_2-p'_1\left(q_1-\frac{Y_1}{p_1}\right)}{p_2}
\right)}$;

- for $\bs{\Delta p > \rho}$ (case (3,1)): \textbf{NE}
$\mathbf{\left(0,\frac{Y_1}{p'_2},0,\frac{Y_2}{p_2} \right)}$.

\textbf{3b)} $Y_2/p'_1<q_1-\alpha$
$\xrightarrow{\textrm{2-III-}B}$ (3,1): we reach the NE from
\textbf{3a)} for $\Delta p > \rho$.

\textbf{The price adjustment process for $\bs{\Delta p \in
(-\infty,\rho]}$ is the same as in the respective cases from zones
1-4, while for $\bs{\Delta p > \rho}$ it is as in the
corresponding case from zone I-3 (basic case).}

\subsubsection{Zone 1-1 (see \eqref{eq:zone1-1init})}\label{sec:sketchZone1-1}

In view of the definition of this zone, the analysis and the
results obtained fully coincide with those for the basic case in
zone I-3.


\begin{thebibliography}{1}

\bibitem{Fil85}
Alexei~F. Filippov.
\newblock {\em Differential Equations with Discontinuous Right-hand Side}.
\newblock Nauka, Moscow, 1985.
\newblock (in Russian).

\bibitem{Fri90}
James~W. Friedman.
\newblock {\em Game Theory with Applications to Economics}.
\newblock Oxford University Press, New York, 2 edition, 1990.

\bibitem{Fud91}
Drew Fudenberg and Jean Tirole.
\newblock {\em Game Theory}.
\newblock Cambridge, MIT Press, 1991.

\bibitem{Ior04}
Iordan~V. Iordanov, Stoyan~V. Stoyanov, and Andrey~A. Vassilev.
\newblock Price dynamics in a two region model with strategic interaction.
\newblock In {\em Proceedings of the {$33^{\textrm{rd}}$} Spring Conference of
  the Union of {B}ulgarian Mathematicians}, pages 144--149, 2004.

\bibitem{Kleb98}
Fima~C. Klebaner.
\newblock {\em Introduction to Stochastic Calculus with Applications}.
\newblock Imperial College Press, London, 1998.

\bibitem{Klo92}
Peter~E. Kloeden and Eckhard Platen.
\newblock {\em Numerical Solution of Stochastic Differential Equations}.
\newblock Springer-Verlag, New York, 1992.

\bibitem{Vas05}
Andrey~A. Vassilev, Iordan~V. Iordanov, and Stoyan~V. Stoyanov.
\newblock A strategic model of trade between two regions with continuous-time
  price dynamics.
\newblock {\em Comptes rendus de l'Acad bul. Sci.}, 58(4):361--366, 2005.

\end{thebibliography}
\end{document}